%October 2, 2008
\magnification=1200
\input amstex
\documentstyle{amsppt}
\hoffset=-0.5pc
\nologo
\vsize=57.2truepc
\hsize=38.5truepc
\spaceskip=.5em plus.25em minus.20em

\define\Nflat{\bold N}
\define\Mflat{\bold M}
\define\SSS{\roman S'}

\define\LLL{\roman B\Lambda_{\partial}[s\fra g]}

\define\chain{\Cal C}
\define\Cat{\Cal D}
\define\VV{\bold V} 
\define\VVV{\bold V}
\define\hgg{\frak g}
\define\Sigm{\roman S}

\define\fra{\frak}
\define\bbar{\roman B}
\define\rbar{\overline{\roman B}}
\define\rcob{\overline{\Omega}}

\define\Bobb{\Bbb}

\define\allpuone{1}
\define\almei{2}
\define\barrone{3}
\define\bottone{4}
\define\botshust{5}
\define\cartanon{6}
\define\cartantw{7}
\define\cartanse{8}
\define\cartanei{9}
\define\cheveile{10}
\define\doldpupp{11}
\define\duponone{12}
\define\duskinon{13}
\define\eilmothr{14}
\define\franzone{15}
\define\franztwo{16}
\define\godebook{17}
\define\gorkomac{18}
\define\gugenhtw{19}
\define\gugenmay{20}
\define\gugenmun{21}
\define\hochsone{22}
\define\hochmost{23}
\define\habili{24}
\define\perturba{25}
\define\cohomolo{26}
\define\modpcoho{27}
\define\intecoho{28}
\define\kan{29}
\define\poiscoho{30}
\define\extensta{31}
\define\lradq{32}

\define\minimult{33}

\define\duaone{34}

\define\pertlie{35}

\define\pertltwo{36}

\define\huebkade{37}

\define\huebstas{38}

\define\husmosta{39}

\define\kamtonfo{40}
\define\kamtontw{41}
\define\kamtonfi{42}
\define\kostaeig{43}
\define\maclafiv{44}
\define\maclaboo{45}
\define\maclbotw{46}
\define\maszwebe{47}
\define\mooretwo{48}
\define\moorefiv{49}
\define\quilltwo{50}
\define\rinehone{51}
\define\gsegatwo{52}
\define\shulmone{53}
\define\stashsev{54}
\define\stashalp{55}
\define\vanestwo{56}
\define\vanesthr{57}
\define\Nsddata#1#2#3#4#5{
 \left( #4
\text{
\vbox
%to 1.15 pc
to 0.8 pc
{
  \hbox {$@>{\,\,\, #3\,\,\,}>>$}
  \vskip-1.2pc
  \hbox {$@<<{\,\,\, #2\,\,\,}<$}
                }
      }
 #1,#5 \right)
}

\define\Mmsddata#1#2#3#4#5{
 ( #1
\text{
\raise 1.5pt
\vbox 
to 1.15 pc
{
  \hbox {$@>{\,\,\, #2\,\,\,}>>$}
  \vskip-1.2pc
  \hbox {$@<<{\,\,\, #3\,\,\,}<$}
                }
      }
 #4,#5 )
}

\topmatter
\title  Relative
 homological algebra,
equivariant de Rham cohomology and Koszul duality
\endtitle
\author Johannes Huebschmann
\endauthor
\affil
Universit\'e des Sciences et Technologies de Lille
\\
UFR de Math\'ematiques
\\
CNRS-UMR 8524
\\
F-59 655 VILLENEUVE D'ASCQ C\'edex, France
\\
Johannes.Huebschmann\@math.univ-lille1.fr
\endaffil
\date{September 29, 2008}
\enddate

\abstract{Let $G$ be a general (not necessarily finite dimensional compact) Lie group,
let $\fra g$ be its Lie algebra, let $C\fra g$ be the cone on $\fra g$ in the 
category of differential graded Lie algebras, and let $\Cal G$ be the functor
which assigns to a chain complex $V$ the $V$-valued total de Rham complex of $G$.
We describe the $G$-equivariant de Rham cohomology in terms of a suitable relative 
differential graded Ext, defined on the appropriate category of $(G,C \fra g)$-modules.
The meaning of \lq\lq relative\rq\rq\ is made precise via the dual standard 
construction associated with the monad involving the functor $\Cal G$ and the 
associated forgetful functor. The corresponding infinitesimal equivariant cohomology 
is the relative differential Ext over $C\fra g$  relative to $\fra g$.
The functor $\Cal G$ decomposes into two functors, the functor which determines 
differentiable cohomology in the sense of Hochschild-Mostow and the functor
which determines the infinitesimal equivariant theory, suitably interpreted.
This functor decomposition, in turn, entails an extension of a Decomposition Lemma 
due to Bott. Appropriate models for the differential graded Ext involving a 
comparison between a suitably defined simplicial Weil coalgebra and the Weil 
coalgebra dual to the familiar ordinary Weil algebra yield small models for
equivariant de Rham cohomology including the standard Weil and Cartan models
for the special case where the group $G$ is compact and connected. 
Koszul duality in de Rham theory results from these considerations
in a straightforward manner.}
\endabstract
\keywords {Relative differential 
derived functors, comonads and standard constructions,
monads and dual standard constructions,
equivariant de Rham cohomology, 
Eilenberg-Moore type descriptions,  
Bott's decomposition lemma,
Weil coalgebra,
simplicial Weil coalgebra,
Eilenberg-Zilber theorem, 
small models for geometric bundles, 
homological perturbation theory, 
Cartan model, Cartan construction, relative Cartan construction,
Weil model, 
Koszul duality}
\endkeywords
\subjclass \nofrills{{\rm 2000} {\it Mathematics Subject
Classification}. \usualspace } {Primary 18G25 55N91; Secondary
16S37 16E45 18G15 18G55 55N10 55N33 55T20 57T30 57U10 58A12 }
\endsubjclass
\toc
\widestnumber\subhead{3.2.}
\specialhead{ } Introduction \endspecialhead
\head 1. The CCE complex of a Lie algebra\endhead
\subhead 1.1. The cone construction\endsubhead
\subhead 1.2. The CCE complex\endsubhead
\subhead 1.3. The $(C\fra h)$-module structure\endsubhead
\subhead 1.4. The category of $(G,C\fra g)$-modules \endsubhead
\subhead 1.5. The de Rham complex of a Lie group\endsubhead
\subhead 1.6.  The de Rham complex of a fiber 
bundle over a homogeneous space \endsubhead
\subhead 1.7.  The diagonal structure on the de Rham complex 
with values in a  $(G,C\fra g)$-module \endsubhead 

\head 2. Relative differential homological algebra\endhead
\subhead 2.1. Adjunctions and (co)monads \endsubhead
\subhead 2.2. Differentiable cohomology \endsubhead
\subhead 2.3. The relative differential $\roman{Ext}_{(C\fra g, \fra g)}$
\endsubhead
\subhead 2.4.  The relative differential 
$\roman{Ext}_{((G,C\fra g);\chain)}$ \endsubhead
\subhead 2.5. The Borel construction\endsubhead
\subhead 2.6. Extension of Bott's decomposition lemma 
\endsubhead
\subhead 2.7. Equivariant de Rham theory as a differential 
Ext \endsubhead
\subhead 2.8. The relative differential 
$\roman{Tor}^{(C\fra g, \fra g)}$\endsubhead
\subhead 2.9. Lie-Rinehart algebras and Lie algebroids\endsubhead
\subhead 2.10. Rational cohomology of algebraic groups\endsubhead
\subhead 2.11. Equivariant de Rham cohomology for algebraic varieties\endsubhead

\head 3. Infinitesimal equivariant (co)homology\endhead
\subhead 3.1.  The Weil coalgebra\endsubhead
\subhead 3.2. The relative bar resolution\endsubhead
\subhead 3.3. Comparison between the Weil coalgebra and the 
relative bar resolution \endsubhead

\subhead 3.4. The Weil coalgebra $W'[\fra g]$ as a relative
$\roman U[\fra g]$-contractible construction in the
reductive case\endsubhead
\subhead 3.5. The Cartan model\endsubhead
\subhead 3.6. Cutting the Cartan model to size\endsubhead
\head 4. The simplicial Weil coalgebra\endhead
\head 5. Cutting the defining object for 
$\roman{Ext}_{((G,C\fra g);\chain)}$ to size \endhead
\subhead 5.1.  The model for 
$\roman{Ext}_{((G,C\fra g);\chain)}$ arising from the Weil coalgebra
\endsubhead
\subhead 5.2. A small object for 
$\roman{Ext}_{((G,C\fra g);\chain)}$ 
in the strictly exterior case\endsubhead
\subhead 5.3. The Weil and Cartan models for compact $G$\endsubhead
\head 6. Small models in equivariant de Rham theory\endhead
\head 7. Duality\endhead
\head{ } References\endhead
\endtoc

\endtopmatter
\document

\leftheadtext{Johannes Huebschmann} \rightheadtext{Homological
algebra and equivariant de Rham cohomology}

\medskip\noindent {\bf Introduction\/}
\smallskip\noindent
The main result of this paper describes 
equivariant 
de Rham theory, in the
spirit of Eilenberg and Moore, in terms of 
a differential
graded Ext, defined
on an appropriate category:
Let $G$ be a Lie group,
let $\fra g$ be its Lie algebra, let $C\fra g$ be the cone
on $\fra g$ in the category of differential graded Lie algebras,
let $\roman{Mod}_{(G,C \fra g)}$ be the category of right
$(G,C \fra g)$-modules where the $G$- and $(C\fra g)$-actions
intertwine in the obvious way,
and let $\chain$ be the category of chain complexes. Given the
$(G,C \fra g)$-module $W$, we define 
$\roman{Ext}_{((G,C \fra g);\chain)}(W,\,\cdot \,)$ to be the {\it right 
derived\/} functor 
(we view the collection of the various Ext as a single functor)
of the functor $\roman{Mod}_{(G,C \fra g)}\to \chain$
which assigns to the $(G,C \fra g)$-module $V$ the
chain complex 
$$
\roman{Hom}(W,V)^{(G,C \fra g)} = 
\roman{Hom}_{(G,C \fra g)}(\Bbb R,\roman{Hom}(W,V))
$$
of invariants,
the real numbers $\Bbb R$ being viewed as a trivial 
$(G,C \fra g)$-module in the obvious way; here 
the convention is to write $\roman{Ext}_{((G,C \fra g);\chain)}$
rather than $\roman{Ext}^*_{((G,C \fra g);\chain)}$,
the term {\it right derived\/} is interpreted in a suitable {\it relative\/}
sense,  the term \lq\lq relative\rq\rq\ 
being made precise by means of  the notions of
{\it monad\/} and {\it dual standard construction\/}.
The requisite categorical language was developed by S. Mac Lane, indeed,
the underlying ideas go back to \cite\maclafiv\ (\S 3).
Thus, let $\Cal G$ be the functor
which assigns to a chain complex $V$ the familiar
$V$-valued (totalized) de Rham complex of $G$.
We define the  differential graded 
$\roman{Ext}_{((G,C \fra g);\chain)}(\, \cdot \, ,\,\cdot \,)$ 
via  the dual standard construction associated
with the monad involving the functor
$\Cal G$ and the corresponding forgetful functor.
The de Rham algebra $\Cal A(X)$ of any smooth $G$-manifold $X$ inherits a
$(G,C\fra g)$-module structure in an obvious manner via the
$G$-action and the
operations of contraction and Lie derivative.
Theorem 2.7.1 below includes the statement that {\sl the $G$-equivariant de Rham
cohomology of a $G$-manifold $X$ is given by the
differential graded\/} 
$\roman{Ext}_{((G,C\fra g);\chain)}(\Bobb R,\Cal A(X))$. 
The infinitesimal version of this differential Ext 
is a suitably defined {\it relative differential\/} 
$\roman{Ext}_{(C \fra g,\fra g)}(\, \cdot \, ,\,\cdot \,)$ over  
$C\fra g$  relative to $\fra g$. 
In (2.3) below,
we introduce this infinitesimal theory accordingly
via the appropriate monad,
and in (2.8) below we spell out
the corresponding comonad.
Occasionally we refer to $\roman{Ext}_{(C \fra g,\fra g)}(\Bbb R ,\,\cdot \,)$
as the {\it infinitesimal\/} equivariant cohomology (relative to $G$
or relative to $\fra g$).
The exploration of the infinitesimal theory
involves notions of {\it Weil coalgebra\/}
and {\it simplicial
Weil coalgebra\/} associated with a Lie algebra. 
The Weil coalgebra of the Lie algebra $\fra g$
is dual to the familiar Weil algebra of $\fra g$;
this coalgebra arises as the ordinary 
differential graded 
Cartan-Chevalley-Eilenberg (CCE) coalgebra
$\Lambda'_{\partial}[sC\fra g]$
of $C\fra g$.
It is well known that the functor which assigns to a 
vector space $V$ the space 
of smooth $V$-valued maps on $G$
and the associated forgetful functor combine to a monad which defines, via the corresponding dual standard construction,
the differentiable cohomology of $G$
in the sense
of Hochschild-Mostow.
For our purposes, a crucial observation is then to the effect that the functor
$\Cal G$ can be written as the {\it composite\/}  of the 
functor which determines {\it differentiable cohomology\/} with the functor
which determines the {\it infinitesimal equivariant theory\/}, suitably interpreted.
This functor decomposition, in turn, leads to an extension,
given as Theorem 2.6.1 below, of Bott's
Decomposition Lemma \cite\bottone.
Appropriate models for the relative differential graded Ext 
involving a comparison between the simplicial 
Weil coalgebra and the Weil coalgebra 
similar to the classical comparison
between the CCE resolution and the
bar complex
yield small models for
equivariant de Rham cohomology including the familiar {\it Weil\/} 
and {\it Cartan\/} models
for the special case where the group is compact and connected. 
Koszul duality
in equivariant de Rham theory then results from these considerations 
in a straightforward manner.
The present paper generalizes in 
particular a result of Bott's
\cite\bottone\ relating the 
Chern-Weil construction
with differentiable cohomology via a 
certain spectral sequence.
Indeed, our approach recovers
equivariant de Rham cohomology in terms of a
suitable higher homotopies construction
having a spectral sequence of the kind considered by Bott as an invariant
and thereby yields in particular, at least in principle, complete
information about the higher differentials
in Bott's spectral sequence.
See Remark 5.1.17 below for details.

Equivariant cohomology is usually defined by means of the 
{\it Borel\/} construction. In the de Rham setting,  given 
the smooth $G$-manifold $X$,
the appropriate way to realize this construction
is to apply the de Rham functor $\Cal A$ to the simplicial Borel
construction
$N(G,X)$ so that the cosimplicial differential graded algebra $\Cal A
(N(G,X))$ results; totalization and normalization then
yield the chain complex $|\Cal A(N(G,X))|$ defining the
$G$-equivariant de Rham theory of $X$ \cite\bottone, \cite\botshust,
\cite\shulmone.
For the special case where the Lie group $G$ is compact,
 older 
constructions of $G$-equivariant cohomology 
in the literature proceed via the 
{\it Weil\/} and {\it Cartan\/} models.
According to folk lore, 
the resulting equivariant cohomology is the
same as that coming from the Borel construction in that particular case; 
indeed, in the literature,
there are various comparison maps between
the Cartan and Weil models and the Borel construction.
These comparison maps establish the equivalence between the various
theories in the compact case but do {\it not explain why\/} these theories
are then equivalent.
Our description of equivariant cohomology in terms of the
aforementioned differential Ext
entails an explanation of the relationship between
the Cartan and Weil models and the Borel construction in a 
{\it conceptual\/} manner:
this relationship  results as a comparison map for
various objects calculating the same derived functor.
The ordinary Weil coalgebra then leads to what we refer to as
the {\it Weil\/} and {\it Cartan\/} models
for the relative differential graded Ext under discussion.
When the group $G$ is compact, any differentiable $G$-module is differentiably
injective, and the differentiable cohomology is non-zero only in degree 
zero and boils down to the $G$-invariants whence, in view of the
aforementioned functor decomposition,
the Weyl and Cartan models 
indeed calculate
the $G$-equivariant cohomology.

In Section 3 below we shall 
explore the infinitesimal equivariant cohomology per se.
In \cite\hochsone, 
for a pair $(\fra a, \fra b)$ of ordinary Lie algebras,
Hochschild has introduced
an acyclic relatively projective CCE complex
which yields the {\it relative Lie algebra cohomology of the pair
$(\fra a, \fra b)$ in the sense of Chevalley-Eilenberg\/} \cite\cheveile.
This CCE complex arises by abstraction from the situation of
the invariant de Rham complex of a homogeneous space of compact 
connected Lie groups.
We shall show that
the literal translation of that CCE construction,
to the pair $(C\fra g,\fra g)$ of differential graded Lie algebras,
yields the Weil coalgebra; see Proposition 3.1.7 below.
In the situation where the smaller Lie algebra
is reductive in the ambient one,
Hochschild's chain complex is
actually  a relatively projective resolution of the ground ring 
\cite\hochsone.
In our case, a similar result holds.
To clarify the situation, extending the idea of a {\it construction\/}
which goes back to H. Cartan \cite\cartanse\ (expos\'e 3), in Section 3 below,
we introduce the notion of {\it relative\/} construction.
Proposition 3.1.7 actually says that the {\it Weil coalgebra $W'[\fra g]$ is 
a  construction for $\roman U[C \fra g]$  relative to $\roman U[\fra g]$
that is $R$-acyclic, even $R$-contractible\/}.
In Theorem 3.4.1 we then show that 
a result similar to that of Hochschild's quoted above holds:
when $\fra g$ is reductive,
the Weil coalgebra $W'[\fra g]$ admits a $\fra g$-equivariant 
contracting homotopy.

In the paper \cite\bottone,
Bott communicates a formula
which he indicates was inspired by some work of Hochschild,
one of the creators of relative 
homological algebra. 
Thus our approach 
explains in particular equivariant cohomology
in terms
of relative homological algebra
and thus closes, perhaps, a circle of ideas.
In Section 4 we shall introduce the already mentioned
simplicial Weil coalgebra.
We will then sometimes refer to the 
Weil coalgebra as the {\it ordinary Weil coalgebra\/},
in particular
when there is a need to distinguish it from 
the simplicial Weil coalgebra.
The ordinary Weil coalgebra leads to a small object calculating
the (relative cohomology which yields the) equivariant cohomology.
The 
canonical comparison between the (normalized
chain complex of the) simplicial Weil coalgebra
and the ordinary
Weil coalgebra, cf. Corollary 4.7 below,
then induces a comparison between the object defining
equivariant cohomology and a small object calculating 
this cohomology, interpreted as the relative derived
functor in the sense explained before;
when a compact connected Lie group is behind,
as already hinted at above,
this procedure leads 
eventually to the Weil and Cartan models.
Thereby the canonical comparison
between the Weil coalgebra and the 
simplicial
Weil coalgebra is formally exactly of the 
same kind
as the classical
comparison, spelled out in detail
in \cite\cartanei\ (chap. XIII),
between the CCE complex for an ordinary Lie algebra
and the bar complex for its universal enveloping algebra.
This relies on the fact,
to be established in Theorem 4.5
below, that 
the (normalized
chain complex of the) simplicial Weil coalgebra
of the Lie algebra $\fra g$
is precisely the (homogeneous form of the) relative bar resolution
for the pair $(C\fra g,\fra g)$.

In Section 5, by means of various HPT techniques
which, in \cite\duaone, we have 
used to construct small models for ordinary singular equivariant
(co)homology, 
we shall cut to size the defining objects for the various
derived functors under discussion.
Using the small objects we shall then show 
in Section 6 that, for a
finite dimensional compact group, 
the ordinary
Weil and Cartan models for equivariant cohomology
result as special cases.
In particular, for a compact group $G$,
 the $G$-{\sl equivariant de
Rham cohomology of $X$ is given by the
invariants of the
 relative differential graded\/}
$\roman{Ext}_{(C\fra g,\fra g)}(\Bobb R,\Cal A(X))$
with respect to the group $\pi_0(G)$ of connected components of $G$, 
and the standard object
calculating this differential graded Ext 
 contracts onto the
Cartan model. Pushing the
HPT-procedure a bit further, we obtain another 
(familiar) model which is even
smaller than the Cartan model for equivariant cohomology. 

In Section 7 we shall exploit
the models constructed in the present paper to introduce, 
via the procedure explained at the
end of \cite\duaone, a certain algebraic duality involving the
object which defines the  differential graded Ext; what is
referred to in the literature as {\it Koszul duality\/},  cf. e.~g.
\cite{\gorkomac}, is an immediate consequence thereof. This yields
a conceptual explanation of 
Koszul duality for de Rham theory in terms of the
extended functoriality of the relevant differential derived
functors and places 
this kind of
Koszul duality in the sh-context. 
The idea behind this extended
functoriality goes back to \cite\stashalp\ and was pushed further
in \cite\gugenmun. For our purposes, the
categories of sh-modules and sh-comodules serve as {\it
replacements for various derived categories\/} exploited in
\cite\gorkomac\ and elsewhere. In particular, when a Lie group $G$
acts on a smooth manifold $X$, even
when the induced action of $\roman H_*G$ on $\roman H^*X$ lifts to
an action on $\Cal A(X)$, in general only an sh-action of $\roman
H_*G$ on $\Cal A(X)$ will recover the geometry of the original
action.

Given a topological group $G$, in ordinary (singular)
(co)homology, the $G$-equivariant (co)homology can be described
via suitable differential Tor- and Ext-functors
in the sense of Eilenberg and Moore over
the chain algebra $C_*G$. When $G$ is an algebraic group,
from the group
multiplication, the algebraic de Rham algebra of $G$ inherits a differential
graded coalgebra, in fact Hopf algebra structure, and the
$G$-equivariant de Rham theory of a nonsingular algebraic variety
is then given by a {\it
differential graded\/} Cotor with respect to this differential
graded coalgebra structure; a similar observation leads
to a description of
rational cohomology of algebraic groups. In the {\it smooth
setting\/}, such a description is of no avail since the smooth
de Rham algebra on a Lie group $G$ does {\it not\/} inherit a diagonal
map (in the usual algebraic sense) turning the de Rham algebra 
into a Hopf algebra.
Our approach in terms of the relative differential Ext
explained above
entails that
a replacement for the
non-existent category of comodules over the de Rham complex of $G$
is provided by the category of $(G,C\fra g)$-modules.
Thus our description of equivariant de Rham cohomology in terms
of a differential Ext {\sl can be seen as a result of the Eilenberg-Moore type.\/} 
To complete the story we shall show,
in (2.10) and (2.11) below, 
that the rational cohomology of algebraic groups and 
the algebraic equivariant de Rham theory of nonsingular
algebraic varieties relative to an algebraic group can likewise
be subsumed under the formalism of monads and dual
standard constructions.

We view the present paper as belonging to a certain
differential homological algebra tradition which started  with
Eilenberg-Mac Lane and H. Cartan and was developed further by J. Moore and
his school. Within this tradition, the theory  
takes care of itself and formulas drop out more or less automatically.
A typical example is the notion of twisting cochain;
once isolated, it explains, in a conceptual way, 
all sorts of perturbed operators
and explicit formulas can then always be derived from structural insight.
For example, given an ordinary Lie algebra $\fra h$, 
in terms of (i) the exterior coalgebra $\Lambda_{\partial}'[s \fra h]$
on the suspension $s\fra h$ of $\fra h$, endowed with the differential
determined by the Lie bracket on $\fra h$ 
and of (ii) the universal algebra $\roman U[\fra h]$ of $\fra h$,
the CCE resolution
can be written in the form
$
\Lambda_{\partial}'[s \fra h] \otimes_{\tau_{\fra h}} \roman U[\fra h]
$
with respect to the corresponding universal twisting cochain
$\tau_{\fra h}$ from $\Lambda_{\partial}'[s \fra h]$ to $\roman U[\fra h]$;
see (1.2.1) below.
Another typical example is the idea of a (co)monad. Yet another example
is given by the description of 
the formalism of contraction and Lie derivative
in terms of an action of the cone on the corresponding Lie algebra;
we shall heavily use this observation in the paper.
The familiar Gerstenhaber algebra structure on the CCE complex calculating the
homology of a Lie algebra actually amounts to a module 
structure over the cone on that Lie algebra; see (1.3) below.
Despite its flexibility and vast range of possible
applications, this differential homological algebra
technology has so far hardly been used in differential geometry.

The reader is assumed to be familiar with the
notation, terminology, and preliminary material in \cite\duaone; 
this material will
not be repeated here. In particular, we will 
use the HPT-techniques explained
in \cite\duaone\ without further explanation. 
As usual, the group of connected components
of a topological group $G$ is written as $\pi_0(G)$.
The ground ring is denoted
by $R$ and the (real) de Rham functor by $\Cal A$.
We treat
chain complexes and cochain complexes on equal footing:
We  consider 
a cochain complex $(C^*,d)$ as a chain complex
$(C_*,d)$ by letting $C_j = C^{-j}$ for $j \in \Bobb Z$.
An ordinary cochain complex, concentrated in non-negative
degrees as a cochain complex, is then a chain complex which
is {\it concentrated in non-positive degrees\/}.
The identity morphism of an object will occasionally 
be denoted by the same symbol as that object and
the operation of suspension will be written as $s$. 
For any smooth manifold $N$,
we write the tangent bundle as $\roman TN \to N$.

I am much indebted to J. Stasheff for a number of comments
on various drafts of the manuscript. I had posted an earlier version
to the arxiv under math.DG/0401161.
Since then, the article \cite\almei\ has appeared,
posted to the arxiv as math.DG/0406350;
that article contains, for the special case
where the group $G$ under discussion is compact, material related to
Subsection 3.6, cf. Remark 3.6.10,  and to Section 6 
below. Publication of our paper
has been delayed for personal 
(non-mathematical)
reasons.

The results presented here can be generalized
to equivariant Lie-Rinehart cohomology
arising from a group acting on an arbitrary Lie-Rinehart algebra.
This is interesting not only in its own right
since this kind of equivariant Lie-Rinehart cohomology arises in arithmetic geometry,
equivariant sheaf theory and algebraic K-theory.
In a different direction, the theory can, perhaps, be extended to 
cover equivariant cohomology relative to actions of
Lie groupoids rather than just Lie groups. We hope 
to return to these issues elsewhere.

\medskip\noindent {\bf 1. The CCE complex of a Lie algebra}
\smallskip\noindent
For later reference we describe various pieces of structure 
on the CCE complex of a Lie
algebra which are most easily explained in terms of twisting cochains.

Let $\fra h$ be an $R$-Lie algebra which we suppose to be free 
or at least projective as an $R$-module
when $R$ is not a field, so that the CCE complex 
then has the desired features,
cf. \cite\barrone.
Here is an example of the kind of Lie algebra we have in mind:
Let $G$ be a Lie group,
with Lie algebra $\fra g$, let
$\xi \colon P \to M$ be a principal $G$-bundle, and let
$R= C^{\infty}(M)$, the algebra of smooth functions on $M$.
Then the space of sections $\fra g(\xi)$ of the {\it adjoint bundle\/}
$\fra g \times_G P \to M$ acquires in an obvious way an $R$-Lie algebra
structure. As an $R$-module,  $\fra g(\xi)$ is projective.
This example justifies building the theory over a ground ring more 
general than a field.

\smallskip
\noindent
{\smc 1.1. The cone construction.\/}
The {\it cone\/} $C\fra h$ on $\fra h$ in the
category of differential graded Lie algebras is the contractible differential
graded Lie algebra $C\fra h$ characterized as follows:
$(C\fra h)_0= \fra h$,
$(C\fra h)_1= s\fra h$, the differential $d\colon
(C\fra h)_1 \to (C\fra h)_0$ is determined by the identity 
$ds=\fra h \ (=\roman{Id}_{\fra h})$,
 the degree 1
constituent $s\fra h$ is abelian, and the action
of $\fra h$ on $s\fra h$ is induced from the adjoint action.
Thus, as a graded 
Lie algebra (i.~e. when the differential is ignored), 
$C\fra h$ can be written as the semi-direct
product $C\fra h = s\fra h \rtimes \fra h$.
The universal
enveloping algebra $\roman U[C\fra h]$ of $C\fra h$ is contractible. As a
graded algebra,
$\roman U[C\fra h]$ decomposes as a crossed product algebra
$\Lambda[s\fra h]\odot \roman U[\fra h]$
relative to the obvious action of the Hopf algebra  $\roman U[\fra h]$
on $\Lambda[s\fra h]$. In particular,
$\Lambda[s\fra h]$ embeds into
$\roman U[C \fra h]$ as a graded subalgebra and $\roman U[\fra h]$
 embeds into
$\roman U[C \fra h]$ as a differential graded subalgebra.

Occasionally we will also use the cone
$\overline C\fra h$ whose underlying graded $R$-module is the same as that of
$C\fra h$ but whose differential is the negative of the differential
of $C\fra h$. The obvious map which is the identity in 
degree zero and multiplication by $-1$ in degree 1 plainly 
identifies the two cones as differential graded Lie algebras. 

Let $V$ be a projective graded $R$-module, concentrated in odd degrees,
and consider the graded exterior algebra $\Lambda[V]$ on $V$.
The diagonal map $V\to V \oplus V$ is well known to induce a diagonal
map for $\Lambda[V]$ turning the latter into a graded Hopf algebra.
We then denote the resulting graded coalgebra by
$\Lambda'[V]$ and, as usual, refer to it as the {\it exterior coalgebra\/}.
Whenever a graded exterior coalgebra of the kind  $\Lambda'[V]$
is under discussion,
we will suppose throughout that the resulting coalgebra is the
graded symmetric coalgebra $\SSS [V]$ on $V$, that is, that the
canonical morphism of coalgebras from 
$\Lambda'[V]$ to $\SSS [V]$ (induced by the canonical projection
from $\Lambda'[V]$ to $V$)
is an isomorphism of graded coalgebras.
This excludes the prime 2 being a zero divisor in the ground ring $R$.
In particular, a field of characteristic 2 is not admitted as ground ring.

\smallskip
\noindent
{\smc 1.2. The CCE complex.\/}
The algebra $\roman U[C\fra h]$ has the CCE resolution $\roman K(\fra h)$
of $R$ (in the category of left $\roman U[\fra h]$-modules
and, suitably modified, in that of  right
$\roman U[\fra h]$-modules, see (1.3) below) 
as its underlying  differential graded
$\roman U[\fra h]$-module, cf. \cite\cartanei\ 
(Ex. XIII.14 
where the ground ring is written as $K$), and we
will identify $\roman U[C\fra h]$ and $\roman K(\fra h)$ 
in notation.
In particular, as a graded
coalgebra, $\roman K(\fra h)$ 
 amounts to the tensor product
$\Lambda'[s \fra h] \otimes \Sigm_{\Delta}[\fra h]$ 
of the graded exterior coalgebra $\Lambda'[s \fra h]$
on $s \fra h$ with the (cocommutative) coalgebra  $\Sigm_{\Delta}[\fra h]$ 
underlying the obvious Hopf algebra structure on the symmetric algebra
$\Sigm[\fra h]$ on $\fra h$, with the
tensor product diagonal.
When $R$ contains the rational numbers as a subring,
$\roman K(\fra h)$ is actually a primitively generated differential
graded Hopf algebra having $C\fra h$ as its space of primitives;
furthermore, as a graded coalgebra,
the constituent 
$\Sigm_{\Delta}[\fra h]$ is isomorphic to
the symmetric coalgebra $\SSS[\fra h]$  on $\fra h$, and
$\roman K(\fra h)$
is the graded symmetric coalgebra cogenerated
by $C\fra h$.

The quotient $\roman K(\fra h)\otimes_{\roman U[\fra h]} R$  calculates
the Lie algebra homology of $\fra h$. This quotient inherits,
furthermore, a differential graded coalgebra
structure having 
the ordinary exterior coalgebra
$\Lambda'[s\fra h]$ on $s \fra h$ 
as its underlying graded coalgebra and
having as differential  
the coderivation $\partial$ corresponding to the
Lie bracket of $\fra h$; we denote this differential graded coalgebra
by $\Lambda'_{\partial}[s \fra h]$ and refer to it
as the CCE {\it coalgebra\/} of $\fra h$. With respect to
the coaugmentation filtration, $\partial$ is a {\it perturbation\/} of the
trivial differential. As a differential graded
left
$(\Lambda_{\partial}'[s \fra h])$-comodule and right
$(\roman U[\fra h])$-module,
$\roman K(\fra h)$ is isomorphic to the twisted tensor product
$$
\Lambda_{\partial}'[s \fra h] \otimes_{\tau_{\fra h}} \roman U[\fra h]
\tag1.2.1
$$
where  $\tau_{\fra h} \colon \Lambda'[s\fra h] \to \roman U [\fra h]$ is the
twisting cochain induced by the differential in $C\fra h$
and in this manner
$\roman K(\fra h)$ appears as a free resolution of $R$ in the category
of right $(\roman U[\fra h])$-modules.
It may, of course, also be rewritten as a
free resolution of $R$ in the category
of {\it left\/} $(\roman U[\fra h])$-modules.

For later use we recall 
some of the technical details
for the case of a general differential graded Lie algebra
where, for simplicity, we suppose that the prime 2 
is invertible in the ground ring.
This is all we need since later in the paper we
shall exclusively work over the reals;
see \cite\pertlie\ and \cite\pertltwo\ 
for the general case:
Let
$C$ be a coaugmented differential graded 
cocommutative coalgebra and $\fra g$
a differential graded Lie algebra which we suppose to be
projective as a graded $R$-module. 
We denote the differential of $C$ and that of $\fra g$ by $d$.
Since $\fra g$ is 
$R$-projective, the 
symmetric 
coalgebra $\SSS[s\fra g]$ on the suspension $s\fra g$ exists;
indeed, this is the cofree
coaugmented differential graded cocommutative 
coalgebra on $s\fra g$.
Let $\tau_{\fra g}\colon \SSS[s\fra g] \to \fra g$ be the homogeneous degree 
$-1$
morphism (of the underlying graded $R$-modules)
which is the desuspension
$ \SSS_{1}[s\fra g]=s\fra g \to \fra g$ from the homogeneous
degree 1 constituent of $\SSS[s\fra g]$ to $\fra g$
and which is zero on the higher degree constituents of
 $\SSS[s\fra g]$.
Given homogeneous morphisms $a,b
\colon C \to \fra g$, with a slight abuse of the bracket notation
$[\, \cdot \, , \, \cdot \, ]$, their {\it cup bracket\/} $[a, b]$
is given by the composite
$$
C @>{\Delta}>> C\otimes C @>{a\otimes b}>>\fra g
\otimes\fra g @> {[\cdot,\cdot]}>> \fra g.
$$
The cup
bracket $[\, \cdot \, , \, \cdot \, ]$ turns $\roman{Hom}(C,\fra g)$ 
into a differential
graded Lie algebra. Define the coderivation
$$
\partial\colon\SSS[s\fra g] \longrightarrow 
\SSS[s\fra g] 
$$
on $\SSS[s\fra g]$  by the requirement that the identity
$$
\tau_{\fra g} \partial 
= \frac 12 [\tau_{\fra g}, \tau_{\fra
g}]\colon \SSS_2[s\fra g] \to \fra g
$$
hold in $\roman{Hom}(\SSS[s\fra g],\fra g)$.
Then $D\partial\  (=d\partial + \partial d) = 0$ since the Lie
algebra structure on $\fra g$ 
is supposed to be compatible with
the differential on $\fra g$. Moreover, the 
property that the bracket $[\,
\cdot \, , \, \cdot \, ]$ on $\fra g$ satisfies the graded Jacobi
identity
is equivalent to the vanishing of $\partial\partial$, that is, to
$\partial$ being a coalgebra 
perturbation of the differential $d$
on $\SSS[s\fra g]$, cf. \cite\pertlie,
\cite\huebstas. 
The resulting
differential graded coalgebra $\SSS_{\partial}[s\fra
g]$ is the 
CCE or {\it classifying\/} coalgebra for $\fra g$;
cf. e.~g. \cite\quilltwo\ (p.~291) for the case where $R$ is the field 
of rational numbers.
When the prime 2 is not invertible in the ground ring, by means of
suitable squaring operations
on $\fra g$ and
 $\roman{Hom}(\SSS[s\fra g],\fra g)$, 
the theory can still be set up but we spare
the reader and ourselves
these added troubles here; see e.~g. \cite\pertlie.

A {\it Lie algebra twisting cochain\/} $t \colon C \to \fra g$ is
a homogeneous morphism of degree $-1$ whose 
composite with the
coaugmentation map of $C$ is zero and which satisfies the equation
$$
Dt = \frac 12 [t,t],
$$
cf. \cite\moorefiv,
\cite\quilltwo,
referred to nowadays
in the literature as {\it deformation equation\/}
or {\it master equation\/}.
When the canonical morphism from $\fra g$ to
$\roman U[\fra g]$ is injective,
the homogeneous degree $-1$ 
morphism $t \colon C \to \fra g$ is a Lie algebra twisting cochain
if and only if the composite of
$t$ with the injection into $\roman U[\fra g]$
is an ordinary twisting cochain.
In particular,
$\tau_{\fra g}\colon \SSS_{\partial}[s\fra
g] \to \fra g$ is a Lie algebra twisting cochain.
When $\fra h$
is an ordinary Lie algebra
(concentrated in degree zero),
$ \SSS_{\partial}[s\fra
h]$ comes down to the ordinary CCE coalgebra of $\fra h$
and, maintaining
notation 
established earlier, we write
$\Lambda'_{\partial}[s \fra h]$
for the CCE coalgebra.
To illustrate our sign conventions we note that, 
given $x_1,x_2 \in \fra h$,
$$
\frac 12 [\tau_{\fra h},\tau_{\fra h}](sx_1 sx_2) =
-[\tau_{\fra h}(sx_1),\tau_{\fra h}(sx_2)]= [x_2,x_1]
$$
whence $\partial(sx_1 sx_2) = [x_2,x_1]$ etc.

For intelligibility we recall the notion of twisted Hom-object,
cf. \cite\duaone\ (2.4.1).
Let $A$ be an
augmented differential graded algebra, 
$C$ a coaugmented differential graded coalgebra,
and $\tau \colon C \to A$ a twisting cochain. Given a
differential graded right
$A$-module $N$ 
let $\delta^{\tau}$ be the operator on $\roman{Hom}(C,N)$
given, for homogeneous $f$, by $\delta^{\tau}(f) = (-1)^{|f|}f \cup \tau$.
With reference to the filtration induced by the coaugmentation
filtration of $C$, 
the operator $\delta^{\tau}$ is a {\it perturbation\/}
of the differential $d$ on $\roman{Hom}(C,N)$, 
and we write the perturbed
differential on 
$\roman{Hom}(C,N)$ as $d^{\tau}= d + \delta^{\tau}$.
Likewise, 
 given a
differential graded left
$A$-module $M$, 
the operator $-\tau \cup\,\cdot\,$ on $\roman{Hom}(C,M)$
is a {\it perturbation\/}
of the differential $d$ on $\roman{Hom}(C,M)$, 
and we write the perturbed
differential on 
$\roman{Hom}(C,M)$ as $d^{\tau}= d - \tau \cup\,\cdot\,$.
We refer to $\roman{Hom}^{\tau}(C,N)=(\roman{Hom}(C,N),d^{\tau})$
and $\roman{Hom}^{\tau}(C,M)=(\roman{Hom}(C,M),d^{\tau})$
as  {\it twisted\/} Hom-{\it objects\/}, cf. \cite\duaone. 

With this preparation out of the way,
let  $N$ be a right $\fra h$-module and $M$ a left
$\fra h$-module. 
The cohomology of $\fra h$ with
values in $N$ ($M$) is calculated  as the homology of the chain complex
$\roman{Hom}_{\roman U[\fra h]}(\roman K(\fra h),N)$
($\roman{Hom}_{\roman U[\fra h]}(\roman K(\fra h),M)$), the 
(differential graded)
subspace of $\roman{Hom}(\roman K(\fra h),N)$ 
($\roman{Hom}(\roman K(\fra h),M)$)
which consists of
$(\roman U[\fra h])$-linear morphisms from 
$\roman K(\fra h)$ to $N$ (to $M$).
The assignment to
$\alpha \in \roman{Hom}(\Lambda'[s\fra h]_{\partial},N)$ of
$$
\Phi_{\alpha} \colon
\Lambda_{\partial}'[s \fra h] \otimes_{\tau_{\fra h}} \roman U[\fra h]
@>>> N,
\quad
\Phi_{\alpha} (w\otimes a)= \alpha(w) a,\ w \in
\Lambda_{\partial}'[s \fra h] ,\ a \in
\roman U[\fra h],
\tag1.2.2
$$
yields an injective chain map
$$
\roman{Hom}^{\tau_{\fra h}}(\Lambda_{\partial}'[s\fra h],N)
@>>>
\roman{Hom}(\roman K(\fra h),N)
\tag1.2.3
$$
which identifies the twisted Hom-object
$\roman{Hom}^{\tau_{\fra h}}(\Lambda_{\partial}'[s\fra h],N)$
with  $\roman{Hom}_{\roman U[\fra h]}(\roman K(\fra h),N)$.
The same kind of association
 identifies the twisted Hom-object
$\roman{Hom}^{\tau_{\fra h}}(\Lambda_{\partial}'[s\fra h],M)$
with  $\roman{Hom}_{\roman U[\fra h]}(\roman K(\fra h),M)$.
The chain complex $\roman{Alt}(\fra h,N)$ of
$N$-valued alternating forms on $\fra h$ with the CCE
differential is exactly 
the source $\roman{Hom}^{\tau_{\fra h}}(\Lambda_{\partial}'[s\fra h],N)$
of (1.2.3).
For $N=R$, we refer
to 
the differential graded algebra $\roman{Alt}(\fra h,R)$ of
$R$-valued alternating forms on $\fra h$ as the 
CCE algebra of $\fra h$ or,
following \cite\vanesthr, as the
{\it Maurer-Cartan\/} algebra
of $\fra h$.

\smallskip\noindent
{\smc 1.3. The $(C\fra h)$-module structures.\/} 
The Lie algebra $\fra h$ acts on
the CCE coalgebra $\Lambda'_{\partial}[s \fra h]$ of $\fra h$
via the action induced by the adjoint action of $\fra h$ on itself.
This 
action is well known to be trivial on (co)homology.
For later reference, we will now refine this observation.

There is a canonical isomorphism
$$
\roman U[\overline C \fra h]\cong 
\Lambda_{\partial}'[s \fra h] \otimes_{\tau_{\fra h}}\roman U[\fra h]
$$
of differential graded
left
$(\Lambda_{\partial}'[s \fra h])$-comodules and right
$(\roman U[\fra h])$-modules and, likewise,
 a canonical isomorphism
$$
\roman U[C \fra h]\cong 
\roman U[\fra h] \otimes_{\tau_{\fra h}}\Lambda_{\partial}'[s \fra h]
$$
of differential graded
right
$(\Lambda_{\partial}'[s \fra h])$-comodules and left
$(\roman U[\fra h])$-modules.
The point here is that, for {\it both\/} isomorphisms,
$\roman U[\fra h]$, ${\tau_{\fra h}}$, and $\Lambda_{\partial}'[s \fra h]$
are the {\it same\/} constituents.
The above isomorphisms entail canonical isomorphisms
$$
\roman U[\overline C \fra h]\otimes_{\roman U[\fra h]}R \cong 
\Lambda_{\partial}'[s \fra h],\quad
R \otimes_{\roman U[\fra h]}\roman U[C \fra h]\cong 
\Lambda_{\partial}'[s \fra h].
$$
Consequently
the CCE coalgebra $\Lambda'_{\partial}[s \fra h]$ of $\fra h$
acquires a differential graded left
$\roman U[\overline C \fra h]$-module structure
and a 
differential graded right
$\roman U[C \fra h]$-module structure.
We write these structures as
$$
\align
\roman U[\overline C\fra h] 
\times \Lambda'_{\partial}[s \fra h] &\longrightarrow 
\Lambda'_{\partial}[s \fra h],\ (a,b) \longmapsto a \cdot b,
\\
\Lambda'_{\partial}[s \fra h] \times \roman U[C\fra h]  
&\longrightarrow 
\Lambda'_{\partial}[s \fra h],\ (b,a) \longmapsto b \cdot a.
\endalign
$$
It is immediate that, given $Y\in \fra h$ and $b \in \Lambda[s \fra h]$,
$$
\gather
Y \cdot b = \roman{ad}_Y(b),\ b\cdot Y  = -\roman{ad}_Y(b)
\\ sY \cdot b = (sY) b, \ b\cdot sY = b(sY)
\ (\text {exterior\ multiplication)}.
\endgather
$$
As a side remark we note that
these  $(C \fra h)$-module and $(\overline C \fra h)$-module  structures
are actually equivalent to the familiar fact that
{\sl the Lie algebra homology
operator $\partial$ generates the Gerstenhaber 
bracket 
$[[\,\cdot \, , \, \cdot \,]]$ 
on\/}
$\Lambda[s \fra h]$,
that is,
for homogeneous $a,b \in \Lambda[s\fra h]$,
$$
\partial (ab)=  (\partial a)b +
(-1)^{|a|} a \partial b  +(-1)^{|a|}[[a,b]].
$$
The ground ring being viewed
as a trivial $(C \fra h)$-module in the obvious way, the
induced $(C \fra h)$-module structure 
($(\overline C \fra h)$-module structure) 
on  $\roman{Alt}(\fra h, R)$
is the familiar action via the operations 
of Lie derivative $\lambda$ and contraction $i$;
thus, given $Y \in \fra h$ and
$\alpha \in \roman{Alt}(\fra h, R)$,
$$
Y(\alpha) = \lambda_Y(\alpha),\ (sY)(\alpha)=i_Y(\alpha),
$$
and the  $(C \fra h)$-action 
on  $\roman{Alt}(\fra h, R)$
is well known to be compatible with the multiplicative
structure. 

\smallskip\noindent
{\smc 1.4. The category of $(G,C\fra g)$-modules.\/}
Let the ground ring $R$ to be that of the real numbers $\Bbb R$.
Let $G$ be a Lie group and let $\fra g$ be its Lie algebra.
We will use the notion of {\it differentiable\/} $G$-module
in the sense of \cite\hochmost.
Henceforth \lq\lq $G$-module\rq\rq\  will mean 
\lq\lq differentiable $G$-module\rq\rq.
Given a differentiable right $G$-module $V$, we will occasionally write
the induced $\fra g$-action as
$$
[\,\cdot \, , \, \cdot \,]
\colon V \times \fra g \longrightarrow V.
$$
In particular, 
for $Y\in\fra g$ 
and $b \in V$,
$$
\frac{\roman d}{\roman dt}(b\,\roman{exp}(tY))\big|_{t=0} = [b,Y] ,
\tag1.4.1
$$
and the actions intertwine in the sense that, given $x \in G$,
$$
[bx,Y] = [b,\roman{Ad}_xY]x .
\tag1.4.2
$$
We use right $G$-modules rather than left ones
since we will eventually apply the theory to left actions of $G$ on smooth 
manifolds; the induced 
$G$-action and 
infinitesimal $\fra g$-action on the functions etc. are then
right actions.

The group $G$ acts on $\roman U[C\fra g]$ compatibly with
the differential graded algebra structure in the obvious way;
indeed $\roman U[C\fra g]$ is a differentiable $G$-module
(in the appropriate category).
Let
$\roman{Mod}_{(G,C \fra g)}$ be the {\it category of right
differential $(G, C\fra g)$-modules 
where the $G$-action is differentiable and where
the actions of $G$ and $C\fra g$ intertwine in
the obvious manner\/}, that is, the obvious extension of (2.4.2)
is satisfied. 
Thus, a $(G, C\fra g)$-module is a chain complex endowed with
a differential right
differentiable $G$-module structure
and a (differential graded) right $(C\fra g)$-module structure
which, restricted to $\fra g$, amounts to the infinitesimal
$\fra g$-module structure induced by the differential
$G$-module structure, and the actions intertwine.
Notice the usage of the adjective \lq\lq differential\rq\rq\
vs that of the adjective \lq\lq differentiable\rq\rq.

Let $\VV$ be a $(G, C\fra g)$-module.
Given $Y\in \fra g$, for convenience, we will occasionally write
the degree zero operator on $\VV$ induced by $Y$ as
a Lie derivative operator
$\lambda_Y\colon \VV \to \VV$ and the degree one operator on
$\VV$ induced by $sY$ as a contraction operator
$i_Y\colon \VV \to \VV$;
the intertwining of the $G$- and $(C\fra g)$-actions
then means that, given $x \in G$, $Y \in \fra g$, and a homogeneous
member $v$ of $\VV$,
$$
(\lambda_Y(v))x= \lambda_{\roman{Ad}_{x^{-1}}Y}(vx),\ 
(i_Y(v))x= i_{\roman{Ad}_{x^{-1}}Y}(vx).
$$

The crucial example of  a  $(G, C\fra g)$-module is the de Rham complex
$\Cal A(X)$ of a smooth $G$-manifold $X$:
In this situation,
the left $G$-action on $X$ induces an action of
the differential graded algebra $\roman U[C\fra g]$ on the
differential graded de Rham algebra $\Cal A(X)$ from the right via
the operations of contraction and Lie-derivative, evaluated
through the infinitesimal anti-action $\fra g \to \roman{Vect}(X)$ of
$\fra g$ on $X$.
The exterior algebra $\Lambda[s \fra g]$ 
being canonically a graded subalgebra of  $\roman U[C\fra g]$
(not a differential graded subalgebra), the
$(\Lambda[s \fra g])$-invariants are then precisely
the ordinary {\it horizontal\/} elements, that is, the forms
$\alpha$ that are horizontal in the sense that
$\alpha(Y,Y_1,\ldots,Y_m)=0$ whenever 
$Y$ is a fundamental vector field on $X$, i.~e. a smooth vector field
on $X$ coming from $\fra g$ via the $G$-action.
For a general
$(G, C\fra g)$-module, we will therefore refer to a
$(\Lambda[s \fra g])$-invariant element as being {\it horizontal\/}.

\smallskip
\noindent
{\smc 1.5. The de Rham complex of a Lie group.\/} 
Let $G$ be a Lie group, and let $\fra g$ be its Lie algebra,
the Lie algebra of left invariant vector fields on $G$ as usual.
Let $V$ be a vector space.
The 
$V$-valued de Rham complex 
$
\Cal A(G,V)
$
of $G$ is well known to amount to the CCE complex
calculating the Lie algebra cohomology of
$\fra g$ with values in $\Cal A^0(G,V)$
relative to the $\fra g$-module structure
coming from left translation or, equivalently, relative to that coming from
right translation.
For later reference,
we will now spell out that CCE complex
relative to the {\it left\/} translation action of $G$ on itself. 
The associated fundamental vector field map is the {\it right\/}
trivialization of the tangent bundle of $G$. 
Beware: This is {\it not\/} the standard identification, which proceeds via
the {\it left\/} trivialization of the tangent bundle of $G$
(and will be explored in the next subsection).

Given $Y \in \fra g$, let $\overline Y$ be the 
associated {\it right invariant\/} vector field,
that is, the vector field 
on $G$ coming from {\it right\/} translation of 
the associated tangent vector $Y_e$ at the identity
element $e$ of $G$.
The fundamental vector field map under discussion
is the right trivialization
$$
\fra g \times G \longrightarrow \roman TG,\ 
(Y,q)\longmapsto \overline Y_q\ (Y \in \fra g,\, q \in G)
\tag1.5.1
$$
of the tangent bundle of $G$.
Relative to the Lie bracket, the resulting morphism 
$$
\fra g \longrightarrow \roman{Vect}(G),
\ Y \longmapsto \overline Y
$$
is anti-Lie and the induced
$\fra g$-action on $\Cal A^0(G,V)$ is from the right,
that is, $\Cal A(G,V)$ appears as a {\it right\/} $\fra g$-module.

To obtain an explicit expression for the identification,
in terms of the fundamental vector field
isomorphism (1.5.1),
of $\Cal A(G,V)$ with the appropriate CCE complex,
given the $p$-form $\alpha$ on $\roman TG$ and $p$ vectors
$Y_1,\ldots,Y_p$ in $\fra g$, let
$$
(\overline\Phi(\alpha))(Y_1,\ldots,Y_p) =\alpha(\overline Y_1,\ldots,\overline Y_p) 
\in \Cal A^0(G,V) .
$$

\proclaim{Proposition 1.5.2}
The morphism
$$
\overline \Phi\colon\Cal A(G,V)
\longrightarrow
\roman{Hom}^{\tau_{\fra g}}(\Lambda'_{\partial}[s\fra g],\Cal
A^0(G,V)) 
\tag1.5.3
$$
is
an isomorphism of chain complexes
between $\Cal A(G,V)$ and
the CCE complex calculating the Lie algebra cohomology of
$\fra g$ with values in the right $\fra g$-module
$\Cal A^0(G,V)$ (coming from left translation in $G$).
When $V$ is a chain complex, the isomorphism
{\rm (1.5.3)} is compatible with the operators 
on both sides of {\rm (1.5.3)}
that are induced by the differential of $V$.
\endproclaim

\demo{Proof} For the left-trivialization of the tangent bundle
of $G$, the corresponding statement is straightforward and classical;
see also Proposition 1.6.3 below.
The argument translates to the right-trivialization
by the standard trick which involves the antipode of the
Hopf algebras coming into play, that is,
the inversion mapping from $G$ to itself and multiplication
by $-1$ on $\fra g$. We leave the details to the reader. \qed
\enddemo

\smallskip\noindent
{\smc 1.6. The de Rham complex of a fiber bundle over a homogeneous space.\/} 
Let $H$ be a Lie group,
let $\fra h$ be its Lie algebra, let
$G$  be a closed subgroup of $H$, let $\fra g$ denote the Lie algebra
of $G$,
let $X$ be a left $G$-manifold, and
consider the de Rham complex $\Cal A(X)$ of $X$, with its induced
right $(G,C\fra g)$-module structure.
An obvious adjointness isomorphism
$$
\Cal A(H,\Cal A(X)) \longrightarrow \Cal A(H\times X)
\tag1.6.1
$$ 
identifies the de Rham complex
$\Cal A(H,\Cal A(X))$
of $H$ with values in the de Rham complex $\Cal A(X)$ of $X$
with the de Rham complex of the product 
$H \times X$ in a $G$-equivariant manner.
Our present aim is to describe the de Rham complex of $H \times_GX$
in terms of the induced  $G$- and $(C\fra g)$-module structures on
a suitable object naturally isomorphic to $\Cal A(H,\Cal A(X))$,
to be spelled out as the right-hand side of (1.6.3.1) below.

The construction of the quotient $H \times_GX$ involves the 
{\it right\/} translation
action of $G$ on $H$. 
The fundamental vector field map 
associated with the $H$-action on itself via right translation
is the {\it left translation\/} trivialization
$$
H \times \fra h_0 \longrightarrow \roman TH
\tag1.6.2
$$
of the tangent bundle $\roman TH \to H$ of $H$.
Thus, unlike the situation of (1.5) above, 
given $Y \in \fra h$, the associated fundamental vector field on $H$
is then simply just $Y$, viewed as a {\it left invariant vector field\/},
and the resulting injection of $\fra h$ into $\roman{Vect}(H)$
is a morphism of Lie algebras; in fact, this is simply the inclusion
of the ordinary Lie algebra $\fra h$ of {\it left-invariant\/} vector fields
into the Lie algebra  $\roman{Vect}(H)$ of all vector fields on $H$.
Thus, via right translation in $H$,
the chain complex
$\Cal A^0(H,\Cal A(X))$ of $\Cal A(X)$-valued functions
on $H$ acquires a left $\fra h$-chain complex
structure which does not involve  
$\Cal A(X)$, and the operator $\delta^{\tau_{\fra h}}$
determined by the universal Lie algebra twisting cochain
$\tau_{\fra h}\colon \Lambda'_{\partial}[s\fra h] \to \roman U[s\fra h]$,
cf. (1.2) above, is defined on
$\roman{Hom}(\Lambda'_{\partial}[s\fra h],\Cal
A^0(H,\Cal A(X)))$.
This operator is a perturbation of the obvious differential 
on $\roman{Hom}(\Lambda'_{\partial}[s\fra h],\Cal
A^0(H,\Cal A(X)))$ coming from $\partial$ and the differential on
$\Cal A(X)$.
For later reference, we spell out the following.

\proclaim{Proposition 1.6.3}
The fundamental vector field
isomorphism {\rm (1.6.2)} induces an isomorphism
$$
\Phi\colon\Cal A(H,\Cal A(X))
\longrightarrow
\roman{Hom}^{\tau_{\fra h}}(\Lambda'_{\partial}[s\fra h],\Cal
A^0(H,\Cal A(X)))
\tag1.6.3.1
$$
of chain complexes that is natural in terms of $H$ and $X$.
This isomorphism admits the following description:
Given the $\Cal A(X)$-valued $p$-form $\alpha$ on $\roman TH$ and $p$ vectors
$Y_1,\ldots,Y_p$ in $\fra h$,
viewed as fundamental vector fields on $H$,
$$
(\Phi(\alpha))(Y_1,\ldots,Y_p) =\alpha(Y_1,\ldots,Y_p) \in \Cal
A^0(H,\Cal A(X)) .
$$
\endproclaim

\demo{Proof} We leave the details to the reader. We only note that, 
for the special case where $X$ is a point, $\Phi$ amounts to the standard
isomorphism of the de Rham complex of $H$ onto the $\Cal A^0(H)$-valued
CCE complex of $\fra h$. \qed
\enddemo

To adjust the situation to the standard principal bundle formalism
where the action of the structure group is from the right,
view the projection from $H \times X$ to $H \times_GX$ as a principal
{\it right\/} $G$-bundle; thus the $G$-action on the product  $H \times X$
from the {\it right\/} is given by the association
$$
H \times X \times G \longrightarrow H \times X,
\quad
(q,x,y) \longmapsto (qy,y^{-1}x),\ q \in H, x \in X, y \in G.
\tag1.6.4
$$
In the standard way, this action induces a 
left $G$-action on $\Cal A(H,\Cal A(X))$. 

\proclaim{Proposition 1.6.5}
Rewritten as a right $G$-action
$$
\Cal A(H,\Cal A(X)) \times G 
\longrightarrow \Cal A(H,\Cal A(X)),
$$
the action of $G$ on $\Cal A(H,\Cal A(X))$ is given by
the assignment
to an alternating $(\Cal A(X))$-valued $p$-form
$\alpha \colon (\roman TH)^p \to \Cal A(X)$ on $\roman TH$ 
($p \geq 0$) and
$y \in G$ of $\alpha \cdot y$, the value $\alpha \cdot y$ 
on  a $p$-tuple $(Z_1,\ldots,Z_p)$ of $p$ vector fields
$Z_1,\ldots,Z_p$ on $H$
being given by
$$
(\alpha \cdot y)(Z_1,\ldots,Z_p) = (\alpha(Z_1 y^{-1},\ldots,Z_p y^{-1}))y .
\tag1.6.5.1
$$
\endproclaim

\demo{Proof} Let $q \in H$, let $a$ be a point of $X$,
let $Z_q \in \roman T_qH$, and let  $U_a \in \roman T_aX$.  The association 
$$
(Z_q ,U_a) \longmapsto
(Z_q \cdot y,y^{-1}\cdot \bold U_a)\ (y \in G)
$$
is the canonical extension
of the right $G$-action (1.6.4) on $H \times X$
to a right $G$-action on  $\roman T(H \times X)$.
Now, under the circumstances of Proposition 1.6.5,
let $\alpha \colon (\roman TH)^{\times p} \to \Cal A(X)$ be
an  $\Cal A^n(X)$-valued $p$-form on $H$,
let $\bold Z =(Z_1,\ldots,Z_p)$, let $y \in G$,
and let  $\bold U=(U_1,\ldots, U_n)$ be an $n$-tuple of vector fields on $X$.
Then
$$
((\alpha \cdot y)_q(\bold Z_q))( \bold U_a)
= (\alpha_{qy^{-1}}(\bold Z_q \cdot y^{-1}))_{ya}
(y\cdot \bold U_a). \qed
\tag1.6.5.2
$$
\enddemo

\proclaim{Corollary 1.6.6}
On the right-hand side
$\roman{Hom}^{\tau_{\fra h}}(\Lambda'_{\partial}[s\fra h],\Cal
A^0(H,\Cal A(X)))$
of {\rm (1.6.3.1)},
the right $G$-action is given by the formula
$$
(\alpha \cdot y)(Y_1,\ldots,Y_p) =
(\alpha(\roman{Ad}_{y} Y_1,\ldots, \roman{Ad}_{y}Y_p))\cdot  y, \ 
y \in G, Y_j \in \fra h;
\tag1.6.7
$$
here $\alpha$ ranges over $\Cal A^0(H,\Cal A(X))$-valued
alternating $p$-forms on $\fra h$, $p \geq 0$, 
and
the expression $(\ldots)\cdot y$ refers to the right $G$-action
on $\Cal A^0(H, \Cal A(X))$ induced by the right translation 
action of $G$ on $H$
and by the left $G$-action on $X$.
\endproclaim

\demo{Proof}
In the formula (1.6.5.2), when each 
vector field
$Z_j$ on $H$ is left-invariant, that is,
a member of $\fra h$, for each $Z_j$, given $y \in G$,
$$
(Z_j)_q\cdot y^{-1} =qZ_j y^{-1} =q y^{-1}\roman{Ad}_{y}Z_j 
= (\roman{Ad}_{y}Z_j)_{q y^{-1}}
$$
whence
$$
(\alpha \cdot y)_q((Z_1)_q,\ldots,(Z_p)_q)((U_1)_a,\ldots,((U_n)_a)
= (\alpha_{q y^{-1}}(\roman{Ad}_{y}\bold Z)_{q y^{-1}})_{ ya}
(y\cdot \bold U_a) . \qed
$$
\enddemo

Thus the $G$-invariant forms relative to the action (1.6.4)
are the $\Cal A(X)$-valued $G$-equivariant forms on $H$ 
and these, in turn, 
in terms of the right-hand side
$\roman{Hom}^{\tau_{\fra h}}(\Lambda'_{\partial}[s\fra h],\Cal
A^0(H,\Cal A(X)))$
of {\rm (1.6.3.1)},
amount to the
$\Cal A^0(H,\Cal A(X))$-valued $G$-equi\-variant alternating forms 
on $\fra h$. 
We will now accordingly characterize the forms that are horizontal
relative to the action (1.6.4)
in terms of the appropriate equivariance property
for $\Cal A^0(H,\Cal A(X))$-valued 
alternating forms on $\fra h$.
To this end,
we first complete the description of the right $(G,C\fra g)$-action on
the right-hand side of (1.6.3.1), as announed at the beginning of 
Subsection 1.6. 
We remind the reader that, as a graded 
Lie algebra, $C\fra g = s\fra g \rtimes \fra g$, the constituent 
$s\fra g$ being abelian.
We will denote the fundamental vector field on $X$ associated with 
$Y \in \fra h$ by $Y_X$ and, accordingly, 
the fundamental vector field on $H\times X$ associated with 
$Y \in \fra h$ by $Y_{H\times X}$.
The fundamental vector field map associated with (1.6.4) takes the form
$$
H \times X \times \fra g \longrightarrow \roman TH \times \roman TX,
\quad
(q,x,Y) 
\longmapsto (Y_q, -(Y_X)_x),\ q \in H, x \in X, Y \in \fra h.
$$
Thus the resulting injection
$\fra g \to \roman{Vect}(H \times X)$
is given by $Y \mapsto (Y,-Y_X)$ ($Y \in \fra g$)
and, the $G$-action on $H \times X$ being from the right,
the resulting infinitesimal $\fra g$-action 
on $\Cal A(H\times X) \cong \Cal A(H,\Cal A(X))$ 
via the operation of Lie derivative
is from the {\it left\/}, i.~e. a 
morphism of Lie algebras (not anti-Lie).

\proclaim{Corollary 1.6.8}
On the right-hand side
$\roman{Hom}^{\tau_{\fra h}}(\Lambda'_{\partial}[s\fra h],\Cal
A^0(H,\Cal A(X)))$
of {\rm (1.6.3.1)},
the induced right $(C\fra g)$-action 
admits the following description:

\noindent
{\rm (i)} The right $\fra g$-action
$$
[\,\cdot\, , \, \cdot \,]\colon
\roman{Hom}^{\tau_{\fra h}}(\Lambda'_{\partial}[s\fra h],\Cal
A^0(H,\Cal A(X))) \times \fra g
\longrightarrow
\roman{Hom}^{\tau_{\fra h}}(\Lambda'_{\partial}[s\fra h],\Cal
A^0(H,\Cal A(X)))
$$
is given by the formula
$$
[\alpha, Y](Y_1,\ldots,Y_p) = \sum
\alpha(Y_1,\ldots, [Y,Y_j],\ldots, Y_p) +
[\alpha(Y_1,\ldots, Y_p),Y],
\tag1.6.9
$$
where $Y \in \fra g$ and $Y_j \in \fra h$;
here $\alpha$ ranges over $\Cal A^0(H,\Cal A(X))$-valued alternating
$p$-forms on $\fra h$, and
the right-most expression $[\ldots, Y]$ 
in {\rm (1.6.9)}
refers to the right $\fra g$-action
on $\Cal A^0(H, \Cal A(X))$ induced by the right translation 
action of $G$ on $H$
and by the left $G$-action on $X$;
thus, given $q \in H$,
$$
[\alpha(Y_1,\ldots, Y_p),Y](q) = -Y_q(\alpha(Y_1,\ldots, Y_p))
+\lambda_{Y_X}((\alpha(Y_1,\ldots, Y_p))(q)).
$$
{\rm N.B. Given 
$Y_1,\ldots,Y_p\in \fra h$, the value
$\alpha(Y_1,\ldots, Y_p)$ is an $\Cal A(X)$-valued function on $H$.}

\noindent
{\rm (ii)}
Given $z = sY\in s \fra g$ where  $Y \in \fra g$ 
and, furthermore,  
the $n$-form  $\alpha=(\alpha_0,\ldots,\alpha_n)$ in
$$
\roman{Hom}^{\tau_{\fra h}}(\Lambda'_{\partial}[s\fra h],\Cal
A^0(H,\Cal A(X))),
$$
with components $
\alpha_j\in 
\roman{Hom}(\Lambda_{j}[s\fra h],\Cal
A^0(H,\Cal A^{n-j}(X)))$
{\rm ($0 \leq j \leq n$)}, the result
$\alpha \cdot z$ 
relative to the corresponding operation
$$
\,\cdot \ \colon
\roman{Hom}^{\tau_{\fra h}}(\Lambda'_{\partial}[s\fra h],\Cal
A^0(H,\Cal A(X))) \times s\fra g
\longrightarrow
\roman{Hom}^{\tau_{\fra h}}(\Lambda'_{\partial}[s\fra h],\Cal
A^0(H,\Cal A(X)))
$$
is given by
$$
\alpha \cdot z= i_{Y_X}\alpha_0 -i_{Y}\alpha_1 +
i_{Y_X}\alpha_1 - i_Y\alpha_2 + \ldots +
i_{Y_X}\alpha_{n-1} - i_Y\alpha_n ;
\tag1.6.10
$$
here, for $0 \leq j \leq n-1$, given $q\in H$ and
$Y_1,Y_2,\ldots,Y_j\in \fra h$,
$$
(i_{Y_X}\alpha_j)(Y_1,Y_2,\ldots,Y_j)(q)
=i_{Y_X}((\alpha_j(Y_1,Y_2,\ldots,Y_j))(q)),
$$
that is, given the vector fields $U_2,\ldots,U_{n-j}$ on
$X$,
$$
\align
(i_{Y_X}&\alpha_j)(Y_1,Y_2,\ldots,Y_j)(q)(U_2,\ldots,U_{n-j})
\\
&=(-1)^j
((\alpha_j(Y_1,Y_2,\ldots,Y_j))(q))(Y_X,U_2,\ldots,U_{n-j}).
\endalign
$$
{\rm N.B. Under these circumstances 
$(\alpha_j(Y_1,Y_2,\ldots,Y_j))(q)\in \Cal A^{n-j}(X)$
($0 \leq j \leq n$).}
\endproclaim

Indeed, on the constituent $\Cal A^p(H,\Cal A^{\ell}(X))$
 ($p,\ell \geq 0$),
the infinitesimal $\fra g$-action 
$$
\Cal A(H,\Cal A(X)) \times \fra g
\longrightarrow \Cal A(H,\Cal A(X)),
\ (\alpha,Y) \longmapsto [\alpha,Y],
$$
on
$\Cal A(H,\Cal A(X))$ from the {\it right\/} associated with (1.6.4)  
is given by
the formula
$$
[\alpha,Y]
=-\lambda_Y\alpha, \ \alpha \colon \roman T^{\times p}H \to \Cal A^{\ell}(X),
$$
suitably interpreted,
in particular, $\Cal A^{\ell}(X)$ is viewed as a right $\fra g$-module
via the left $G$-action on $X$ and the associated operation of Lie 
derivative.
Explicitly, given the vector fields $Z_1,\ldots, Z_p$ on $H$,
$$
[\alpha,Y](Z_1,\ldots,Z_p)= 
\sum \alpha(Z_1,\ldots, [Y,Z_j],\ldots,Z_p) 
+\lambda_{Y_X}(\alpha(Z_1,\ldots,Z_p))
\tag1.6.11
$$
(beware the parentheses: $\alpha(Z_1,\ldots,Z_p)$ is an $\ell$-form on $X$), 
that is,  given in addition the vector fields $U_1,\ldots, U_\ell$ on $X$,
$$
\aligned
[\alpha,Y](Z_1,\ldots,Z_p)(U_1,\ldots, U_\ell)  &= 
\sum \alpha(Z_1,\ldots, [Z_j,Y],\ldots,Z_p)(U_1,\ldots, U_\ell)  
\\
&\quad + Y_X(\alpha(Z_1,\ldots,Z_p)(U_1,\ldots, U_\ell))
\\
&\quad -\sum \alpha(Z_1,\ldots,Z_p)(U_1,\ldots,[Y_X,U_j],\ldots, U_\ell).
\endaligned
$$
Likewise,
on $\Cal A(H,\Cal A(X))$,
 the operation of contraction
with vectors in $\fra g$ induced by (1.6.4) 
can be described as follows: Given $Y \in \fra g$,
the operation of contraction 
$$
i_{Y_{H\times X}}\colon \Cal A(H,\Cal A(X)) \cong\Cal A(H \times X)
\longrightarrow \Cal A(H \times X) \cong \Cal A(H,\Cal A(X))
$$
with the fundamental vector field $i_{Y_{H\times X}}=(Y,-Y_X)$ 
associated with $Y$
is the operation
$$
i_{Y_{H\times X}}\alpha = i_{Y}\alpha - i_{Y_{X}}\alpha.
$$
An $n$-form $\alpha$ in $\Cal A(H,\Cal A(X))$ has $n+1$
components 
$
\alpha_j\in \Cal A^j(H, \Cal A^{n-j}(X)) 
$
($0 \leq j \leq n$) and, 
relative to the right $G$-action (1.6.4) on $H \times X$,
given $Y \in \fra g$, when $z=sY \in s\fra g$,
$$
\alpha \cdot z=-i_{Y_{H \times X}}\alpha = i_{Y_X}\alpha_0 -i_{Y}\alpha_1 +
i_{Y_X}\alpha_1 - i_Y\alpha_2 + \ldots +
i_{Y_X}\alpha_{n-1} - i_Y\alpha_n .
\tag1.6.12
$$
Summing up, we arrive at the following.

\proclaim{Proposition 1.6.13}
{\rm (i)} 
The $n$-form $\alpha=(\alpha_0,\ldots,\alpha_n)$ in $\Cal A(H,\Cal A(X))$ 
is horizontal if and only if $\alpha \cdot z$ is zero for
every $z \in  s\fra g$, that is,
if and only if
$$
\align
 i_{Y_X}\alpha_0 &=i_{Y}\alpha_1 \in  \Cal A^0(H, \Cal A^{n-1}(X))
\\
i_{Y_X}\alpha_1 &= i_Y\alpha_2 \in  \Cal A^1(H, \Cal A^{n-2}(X))
\\
\cdot &\cdot  \cdot
\\
i_{Y_X}\alpha_{n-1} &= i_Y\alpha_n \in  \Cal A^{n-1}(H, \Cal A^{0}(X))
\endalign
$$
for every $Y \in \fra g$.

\noindent
{\rm (ii)} The $n$-form  $\alpha$ in $\Cal A(H,\Cal A(X))$ 
is basic in the sense that it descends to an $n$-form on $H \times_GX$,
i.~e. is horizontal and $G$-invariant, if and only if
it is $(G,C\fra g)$-invariant. \qed
\endproclaim

\noindent
{\smc 1.7. The diagonal structure on the 
de Rham complex with values in a $(G,C\fra g)$-module.\/}
Abstracting from the material in Subsection 1.6, we now replace
the de Rham complex $\Cal A(X)$ with a general $(G,C\fra g)$-module
$\VV$.
Right translation in $H$ and the $G$-action on $\VV$ induce a
right $(G,C\fra g)$-module structure
$$
\Cal A(H,\VV) \times (G,C \fra g) 
\longrightarrow \Cal A(H,\VV)
\tag1.7.1
$$
on $\Cal A(H,\VV)$.
An explicit description thereof
is given by
the formulas (1.6.5.1),  (1.6.11) and (1.6.12),
with $\VV$ substituted for $\Cal A(X)$.

Similarly as before,
via right translation in $H$,
the chain complex
$\Cal A^0(H,\VV)$ of $\VV$-valued functions
on $H$ acquires a left $\fra h$-chain complex
structure that does not involve 
$\VV$, and the operator $\delta^{\tau_{\fra h}}$
determined by the universal Lie algebra twisting cochain
$\tau_{\fra h}\colon \Lambda'_{\partial}[s\fra h] \to \roman U[s\fra h]$,
cf. (1.2) above, is defined on
$\roman{Hom}(\Lambda'_{\partial}[s\fra h],\Cal
A^0(H,\VV))$.
This operator is a perturbation of the obvious differential 
on $\roman{Hom}(\Lambda'_{\partial}[s\fra h],\Cal
A^0(H,\VV))$ coming from $\partial$ and the differential on
$\VV$.
Now, the formulas (1.6.7), (1.6.9) and (1.6.10),
with $\VV$ instead of $\Cal A(X)$,
yield a right $(G,C\fra g)$-module structure
$$
\roman{Hom}^{\tau_{\fra h}}(\Lambda'_{\partial}[s\fra h],\Cal
A^0(H,\VV)) \times (G,C\fra g) 
\longrightarrow 
\roman{Hom}^{\tau_{\fra h}}(\Lambda'_{\partial}[s\fra h],\Cal
A^0(H,\VV))
\tag1.7.2
$$
on $\roman{Hom}^{\tau_{\fra h}}(\Lambda'_{\partial}[s\fra h],\Cal
A^0(H,\VV))$.

For completeness, we spell out the result of the action with
an element $z=sY$ of the constituent
$s \fra g$ of $C \fra g= s \fra g \rtimes \fra g$ where $Y \in \fra g$.
To this end,
for $0 \leq j \leq n-1$, given $q\in H$, vector fields
$Y_1,Y_2,\ldots,Y_j$ on $H$, and 
the $\VV^{n-j}$-valued $j$-form
$\alpha_j\in \Cal A^j(H,\VV^{n-j})$,
the value $(i_{Y_{\VV}}\alpha_j)(Y_1,Y_2,\ldots,Y_j)(q)$
is given by
$$
(i_{Y_{\VV}}\alpha_j)(Y_1,Y_2,\ldots,Y_j)(q)
=i_{Y_{\VV}}((\alpha_j(Y_1,Y_2,\ldots,Y_j))(q))=
((\alpha_j(Y_1,Y_2,\ldots,Y_j))(q))\cdot z.
$$
N.B. The value $(\alpha_j(Y_1,Y_2,\ldots,Y_j))(q)$ lies in $\VV^{n-j}$,
and $((\alpha_j(Y_1,Y_2,\ldots,Y_j))(q))\cdot z$ 
lies in $\VV^{n-j-1}$.
With this preparation out of the way,
given the $n$-form $\alpha=(\alpha_0,\ldots,\alpha_n)$ 
$\Cal A(H,\VV)$,
with components 
$\alpha_j\in \Cal A^j(H,\VV^{n-j})$ ($0 \leq j \leq n$),
the value $\alpha \cdot z$ 
relative to the corresponding operation
$$
\,\cdot \ \colon
\roman{Hom}^{\tau_{\fra h}}(\Lambda'_{\partial}[s\fra h],\Cal
A^0(H,\VV)) \times s\fra g
\longrightarrow
\roman{Hom}^{\tau_{\fra h}}(\Lambda'_{\partial}[s\fra h],\Cal
A^0(H,\VV))
$$
is given by 
$$
\alpha \cdot z= i_{Y_{\VV}}\alpha_0 -i_{Y}\alpha_1 +
i_{Y_{\VV}}\alpha_1 - i_Y\alpha_2 + \ldots +
i_{Y_{\VV}}\alpha_{n-1} - i_Y\alpha_n .
$$

\proclaim{Proposition 1.7.3}
The fundamental vector field
isomorphism {\rm (1.6.2)} induces an isomorphism
$$
\Phi\colon\Cal A(H,\VV)
\longrightarrow
\roman{Hom}^{\tau_{\fra h}}(\Lambda'_{\partial}[s\fra h],\Cal
A^0(H,\VV))
\tag1.7.3.1
$$
of right $(G,C\fra g)$-modules that is natural in terms of $H$ and
$\VV$.
This isomorphism admits the following explicit description:
Given the $\VV$-valued $p$-form $\alpha$ on $\roman TH$ and $p$ vectors
$Y_1,\ldots,Y_p$ in $\fra h$,
viewed as fundamental vector fields on $H$,
$$
(\Phi(\alpha))(Y_1,\ldots,Y_p) =\alpha(Y_1,\ldots,Y_p) \in \Cal
A^0(H,\VV) .
$$
\endproclaim

\demo{Proof} The reasoning is exactly the same as that for Corollary 1.6.6
and Corollary 1.6.8. \qed
\enddemo

\medskip\noindent{\bf 2. Relative differential homological algebra}
\smallskip\noindent
{\smc 2.1. Adjunctions and (co)monads.\/} 
Before going into details we note that we avoid the terminology \lq\lq
triple\rq\rq\ etc. and exclusively use the monad-comonad terminology.
An adjunction is well known to determine a monad and a comonad 
\cite\maclbotw:
Let $\Cat$ and $\Cal M$ be categories, let $\Cal G \colon \Cat \to \Cal M$
be a functor, suppose that the functor $\square\colon \Cal M \to \Cat$ 
is left-adjoint to $\Cal G$, and let 
$$
\Cal T= \Cal G \square\colon \Cal M \longrightarrow \Cal M.
$$
Let $\Cal I$ denote the identity functor,
let $\eta\colon \Cal I \to  \Cal T$ be the {\it unit\/},
$\varepsilon \colon   \square \Cal G \to \Cal I$ the {\it counit\/}
of the adjunction, and let $\mu$ be the natural transformation
$$
\mu = \Cal G \varepsilon \square\colon\Cal G \square\Cal G \square= 
\Cal T^2 \longrightarrow \Cal T =\Cal G \square.
$$
The data $(\Cal T,\eta, \mu)$ constitute a {\it monad\/}
over the category $\Cal M$.
The  {\it dual standard construction\/},
cf. \cite\duskinon, \cite\godebook\ (\lq\lq construction fondamentale\rq\rq\ 
on p.~271), \cite\maclbotw,
then yields the cosimplicial object
$$
\left(\Cal T^{n+1},  \varepsilon^j\colon \Cal T^{n+1} \to \Cal T^{n+2},
\eta^j\colon
\Cal T^{n+2} \to \Cal T^{n+1}
\right)_{n \in \Bbb N};
$$
here, for $n \geq 0$,
$$
\align
\varepsilon^j&= \Cal T^j \eta \Cal T^{n-j+1}\colon \Cal T^{n+1} \to \Cal T^{n+2},
\ j = 0, \ldots, n+1,
\\
\eta^j&= \Cal T^j \mu \Cal T^{n-j}\colon \Cal T^{n+2} \to \Cal T^{n+1},
\ j = 0, \ldots, n.
\endalign
$$
Thus, given an object $V$ of $\Cal M$,
$$
\bold T(V)=\left(\Cal T^{n+1}(V),\varepsilon^j,\eta^j\right)
$$
is a cosimplicial object in $\Cal M$;
here we do not distinguish in notation
between the natural transformations  $\eta^j$ and $\varepsilon^j$ 
and the morphisms they induce after evaluation of
the corresponding functors in an object.

Under suitable circumstances,
e.~g. when $\Cal M$ is a category of modules, the
associated chain complex $\left|\bold T(V)\right|$ 
is a relatively injective resolution
of $V$, more precisely,  a resolution of $V$ in the category $\Cal M$
that is injective relative to the category $\Cat$.
We will use this construction to introduce and exploit various relative
differential $\roman{Ext}_{(\Cal M,\Cat)}$-functors 
over suitable categories
$\Cal M$ and $\Cat$. Examples will
be given shortly.
For sheaves, this kind of construction goes back to \cite\godebook\ 
(pp. 270--279).

Likewise, let $\Cal F \colon \Cat \to \Cal M$
be a functor, suppose that the functor $\square\colon \Cal M \to \Cat$ 
is right-adjoint to $\Cal F$, and let 
$$
\Cal L= \Cal F \square\colon \Cal M \longrightarrow \Cal M.
$$
Let $\eta\colon \Cal I \to  \square\Cal F$ be the {\it unit\/},
$\varepsilon \colon  \Cal L \to \Cal I$ the {\it counit\/}
of the adjunction, and let $\delta$ be the natural transformation
$$
\delta = \Cal F \eta \square\colon \Cal L =\Cal F \square\longrightarrow 
\Cal F \square\Cal F \square= \Cal L^2.
$$
The data $(\Cal L,\epsilon, \delta)$ constitute a {\it comonad\/}
over the category $\Cal M$.
The {\it standard construction\/}
then yields the simplicial object
$$
\left(\Cal L^{n+1}, d_j\colon
\Cal L^{n+2} \to \Cal L^{n+1}, s_j\colon \Cal L^{n+1} \to \Cal L^{n+2}
\right)_{n \in \Bbb N};
$$
here, for $n \geq 0$,
$$
\align
d_j&= \Cal L^j \varepsilon \Cal L^{n-j+1}\colon \Cal L^{n+2} \to \Cal L^{n+1},
\ j = 0, \ldots, n+1,
\\
s_j&= \Cal L^j \delta \Cal L^{n-j}\colon \Cal L^{n+1} \to \Cal L^{n+2},
\ j = 0, \ldots, n.
\endalign
$$
Thus, given an object $W$ of $\Cal M$, 
$$
\bold L(W)=\left(\Cal L^{n+1}(W), d_j, s_j\right)_{n \in \Bbb N}
$$
is a simplicial object in $\Cal M$,
the {\it standard object associated with\/} $W$ {\it and the comonad\/};
again, we do not distinguish in notation
between the natural transformations  $d_j$ and $s_j$ 
and the morphisms they induce after evaluation of
the corresponding functors in an object.
Under suitable circumstances,
the associated chain complex $\left|\bold L(W)\right|$ is a relatively projective resolution
of $W$, more precisely,  a resolution of $W$ in the category $\Cal M$
that is projective relative to the category $\Cat$.
We will use this construction to introduce and exploit
certain relative
differential $\roman{Tor}^{(\Cal M,\Cat)}$-
and $\roman{Ext}_{(\Cal M,\Cat)}$-functors
 over certain categories
$\Cal M$ and $\Cat$.

\smallskip
\noindent
{\smc 2.2. Differentiable cohomology.\/} 
Let the ground ring be that of the reals, $\Bbb R$.
Recall that $\Cal A$ refers to the de Rham functor on smooth manifolds.
Thus  $\Cal A^0$ refers to ordinary smooth functions.
Let $G$ be a Lie group.
It is well known that, contrary to what is the case for projective
resolutions, the mechanism of {\it injective\/} resolutions can be adapted
to take account of additional structure, here that of
{\it differentiability\/} of a $G$-module.
Indeed, the appropriate way to resolve an object of the
category of differentiable $G$-modules 
is by means of a {\it differentiably injective resolution\/}
\cite\hochmost:
Let $\Cat=\roman{Vect}$, the category of real vector spaces,
$\Cal M=\roman{Mod}_{G}$ that of (differentiable) right $G$-modules, and let
$\Cal G_G \colon \roman{Vect} \to  \roman{Mod}_{G}$
be the functor
which assigns to the real vector space
$V$ the 
$G$-representation
$$
\Cal G_GV= \Cal A^0(G, V),
$$
endowed with the right $G$-module structure coming from 
left translation on $G$.
For our purposes it would be more appropriate to endow
$\Cal A^0(G, V)$ with the right module structure coming from right 
translation in $G$ 
combined with the inversion mapping of $G$, but to arrive
at formulas consistent with 
what is in the literature we proceed with the
left translation action of $G$ on itself.
We use the font $\Cal G$ merely for convenience since this is reminiscent
of the notation $G$ in \cite\maclbotw\ 
for this kind of functor; this usage of the font $\Cal G$
has nothing to do with our usage of the notation $G$ 
for the group variable.
The functor $\Cal G_G$
is right adjoint to the forgetful functor $\square 
\colon \roman{Mod}_{G} \to  \roman{Vect}$
and hence defines a monad  $(\Cal T,\eta,\mu)$  over the category 
$\roman{Mod}_{G}$.
Given a $G$-module $V$, 
the chain complex arising from the
{\it dual standard construction\/}
$\bold T(V)$ associated with $V$
is the
standard differentiably injective resolution of $V$
in the category of $G$-modules
which defines the differentiable cohomology $\roman H_{\roman{cont}}^*(G,V)$.
Here we write $\roman H^*_{\roman{cont}}$ since the differentiable 
cohomology with coefficients in a differentiable module coincides
with the continuous cohomology with coefficients in that module.
See \cite\hochmost\ for details. An explicit description of this resolution
will be given in (2.5) below.
The same kind of construction works for continuous cohomology
but, in this paper, we shall exclusively exploit the differentiable version.

\smallskip
\noindent
{\smc 2.3. The relative differential\/} $\roman{Ext}_{(C \fra g, \fra g)}$.
Let $R$ be an arbitrary commutative ring with 1 and $\fra g$ 
an $R$-Lie algebra, projective as an $R$-module as 
throughout the paper.
Let $\Cat=\chain_{\fra g}$, the category of right $\fra g$-chain complexes,
let $\Cal M=\roman{Mod}_{C\fra g}$,
and let $\Cal G^{\fra g}_{C\fra g} \colon\chain_{\fra g}\to \roman{Mod}_{C\fra g}$
be the functor given by
$$
\Cal G^{\fra g}_{C\fra g}(V)= \roman{Hom}_{\fra g}(\roman U[C\fra g],V)
\cong
\roman{Hom}^{\tau_{\fra g}}(\Lambda'_{\partial}[s\fra g],V) \cong
(\roman{Alt}(\fra g, V),d),
\tag2.3.1
$$
the total object arising from the bicomplex having 
$\roman{Alt}^*(\fra g, V_*)$ as underlying bigraded $R$-module;
here $V$ ranges over  right $\fra g$-chain complexes,
$(\roman{Alt}(\fra g, V),d)$ is
endowed with the obvious right $(C\fra g)$-module structure coming from
the obvious left $(\roman U[C\fra g])$-module structure on itself or,
equivalently, that given by the operations of 
contraction and Lie derivative on the CCE complex
$(\roman{Alt}(\fra g, V),d)$, cf. (1.3) above.
The functor $\Cal G^{\fra g}_{C\fra g}$
is right adjoint to the forgetful functor $\square 
\colon \roman{Mod}_{C\fra g}\to \chain_{\fra g}$
and hence defines a monad  $(\Cal T,\eta,\mu)$  over the category 
$\roman{Mod}_{C\fra g}$.
Given a right $(C\fra g)$-module $\VV$, 
the chain complex $\left|\bold T(\VV)\right|$ arising from the
{\it dual standard construction\/}
$\bold T(\VV)$ associated with $\VV$
is a resolution of $\VV$
in the category of $(C\fra g)$-modules that is injective relative to
the category $\chain_{\fra g}$ of right $\fra g$-chain complexes.
Given a right  $(C\fra g)$-module $\bold W$, the {\it relative differential 
$\roman{Ext}_{(C \fra g, \fra g)}(\bold W,\VV)$\/} is the homology of
the chain complex
$$
\roman{Hom}_{C\fra g}\left(\bold W,\left|\bold T(\VV)\right|\right).
$$
In particular, for $\bold W=R$, 
the relative differential graded 
$\roman{Ext}_{(C \fra g, \fra g)}(R,\VV)$ is the homology of
the chain complex
$
\left|\bold T(\VV)\right|^{C\fra g}.
$
\smallskip
\noindent
{\smc 2.4. The relative differential\/} $\roman{Ext}_{((G,C \fra g);\chain)}$.
As before, view the group  $G$ as a left $G$-manifold via left translation.
Given the chain complex $V$,
let $\Cal A(G,V)$, the $V$-{\it valued (totalized)
de Rham complex $\Cal A(G,V)$ of $G$\/} be the 
chain complex arising from the operation of totalization
applied to the bicomplex $(\Cal A^*(G,V_*),\delta, d)$,
where $\delta$ refers to the de Rham complex operator and $d$ to 
the differential induced by the differential of $V$,
and endow $\Cal A(G,V)$
with the
right $(G,C\fra g)$-module structure explained above.
Let $\chain$ be the category of real chain complexes.
Consider the {\it pair of categories\/}
$
\left(\Cal M,\Cat\right)=
\left(\roman{Mod}_{(G,C \fra g)},\chain\right)
$
and let
$$
\Cal G_{(G,C \fra g)} \colon \chain \longrightarrow  \roman{Mod}_{(G,C\fra g)}
\tag2.4.1
$$
be the functor
which assigns to the chain complex 
$V$ the right $(G,C\fra g)$-module
$$
\Cal G_{(G,C \fra g)}V= \Cal A(G,V).
\tag2.4.2
$$

\proclaim{Proposition 2.4.3}
The functor $\Cal G_{(G,C \fra g)}$
is right adjoint to the forgetful functor $\square 
\colon \roman{Mod}_{(G,C\fra g)}\to  \chain$
and hence defines a monad  $(\Cal T,\eta,\mu)$  over the category 
$\roman{Mod}_{(G,C\fra g)}$.
\endproclaim

\demo{Proof}
Let $W$ be a vector space, $\VV$ a right $(G,C\fra g)$-module, and denote the
graded vector space which underlies $\VV$ by $\VV^{\flat}$
The obvious linear map
$$
\roman{Hom}(\VV^{\flat},W) \longrightarrow \roman{Hom}_{(G,C\fra g)}(\VV,\Cal A(G,W))
\tag2.4.4
$$
sends the homogeneous linear map
$\varphi \colon \VV^{\flat} \to W$ to the 
$(G,C\fra g)$-linear morphism
$$
\Phi \colon \VV \longrightarrow \Cal A(G,W) 
\cong \roman{Hom}_{\fra g}(\roman U[C\fra g],\Cal A^0(G,W))
$$
determined by the requirement that, 
for a homogeneous member $\bold v$ of $\VV$ of degree $-k \leq 0$, 
the value $\Phi(\bold v)$ be the $W$-valued $k$-form on $G$ such that,
given $Y_1,\ldots,Y_k \in \fra g$,
$$
\Phi(\bold v) (sY_1 \ldots sY_k)=\Phi(\bold v sY_1 \ldots sY_k ) 
$$
is
the smooth $W$-valued
function 
on $G$ given by
$$
\Phi(\bold v sY_1 \ldots sY_k ) (x) = 
 \varphi(\bold v sY_1 \ldots sY_k x), \ x \in G.
$$
Here the juxtaposition
$$
(\bold v sY_1 \ldots sY_k, x) \longmapsto \bold v sY_1 \ldots sY_k x
$$
refers to the $G$-action on $\VV$.
The linear map (2.4.4) is an isomorphism of vector spaces.
Replacing $W$ with a chain complex and taking 
the totalized object $\Cal A(G,W)$,
we arrive at the desired adjunction
$$
\roman{Hom}_{\Cal C}(\square \VV,W) \longrightarrow \roman{Hom}_{(G,C\fra g)}(\VV,\Cal A(G,W)). \qed
$$
\enddemo

In view of the general theory reproduced in (2.1) above,
the adjunction spelled out in Proposition 2.4.3
yields a monad $(\Cal T, \eta,\mu)$ over the category 
$\roman{Mod}_{(G,C\fra g)}$.

Let $\chain_G$ be the category of right $G$-chain complexes.
On the category $\chain_G$, the obvious variant of the functor
$\Cal G^{\fra g}_{C\fra g}$, cf. (2.3.1), takes the form
of the functor
$$
 \Cal G^G_{(G,C\fra g)}\colon \chain_G \longrightarrow
\roman{Mod}_{(G,C\fra g)},
\  \Cal G^G_{(G,C\fra g)}(V) 
=\roman{Hom}^{\tau_{\fra g}}(\Lambda'_{\partial}[s\fra g], V)
\cong
(\roman{Alt}(\fra g, V),d)
\tag2.4.5
$$
where $V$ ranges over right $G$-chain complexes.

\proclaim{Proposition 2.4.6}
The functor
$\Cal G_{(G,C\fra g)}$ admits the decomposition
$$
\Cal G_{(G,C\fra g)}= \Cal G^G_{(G,C\fra g)}\circ \Cal G_{G}
\colon \chain \longrightarrow \chain_G \longrightarrow
\roman{Mod}_{(G,C\fra g)}.
$$
\endproclaim

\demo{Proof} This is an immediate consequence of 
Proposition 1.5.2. \qed
\enddemo

Let $\VV$ be  a right $(G,C\fra g)$-module. 
The chain complex $\left|\bold T(\VV)\right|$ arising from the
{\it dual standard construction\/}
$\bold T(\VV)$ associated with 
the monad $(\Cal T, \eta,\mu)$ and 
the right $(G,C\fra g)$-module
$\VV$
is a resolution of $\VV$
in the category of right
$(G,C\fra g)$-modules
that is {\it injective relative to the category of chain complexes\/}.
Given a right  
$(G,C\fra g)$-module $\bold W$, the {\it differential graded\/}
$\roman{Ext}_{((G,C \fra g);\chain)}(\bold W,\VV)$ is the homology of
the chain complex
$$
\roman{Hom}_{{(G,C\fra g)}}\left(\bold W,\left|\bold T(\VV)\right|\right).
$$
In particular, relative to the obvious 
trivial $(G,C \fra g)$-module structure
on $\Bobb R$, the  differential graded
$
\roman{Ext}_{((G,C \fra g);\chain)}(\Bbb R,\VV)
$ 
is the homology of
the chain complex
$
\left|\bold T(\VV)\right|^{(G,C\fra g)}
$
of $(G,C\fra g)$-invariants in $\left|\bold T(\VV)\right|$.

\smallskip
\noindent
{\smc 2.5. The Borel construction\/}.
As before, $G$ denotes a Lie group, neither necessarily compact
nor necessarily finite dimensional. Let $X$ be a left
$G$-manifold. The  {\it simplicial Borel construction\/} 
takes the form of either a {\it nonhomogeneous\/} construction
or  of a {\it homogeneous construction\/}. 

Recall that any object $Y$ of a 
symmetric monoidal 
category endowed with 
a cocommutative diagonal---we will take 
the categories of spaces, of
smooth manifolds, of groups, of vector spaces, 
of Lie algebras,
etc.,---defines two simplicial 
objects in the category, the {\it trivial\/}
object which, with an abuse of notation,
we still write as $Y$, and the {\it total object\/} $EY$
(\lq\lq total object\rq\rq\ 
not being standard terminology in this generality); the
trivial object $Y$ has a copy of $Y$ in each 
degree and 
all simplicial
operations are the identity while, for $p \geq 0$, 
the degree $p$  constituent
$EY_p$ of the total object $EY$ is a product of
$p+1$ copies of $Y$ with the familiar 
face operations given by omission and 
degeneracy operations given by insertion.
See e.~g. \cite\bottone\ and \cite\duaone\
(1.1).
When $Y$ is an ordinary 
$R$-module,
the simplicial $R$-module
associated with
$Y$ is in fact the result of application of the 
{\it Dold-Kan\/} functor $DK$ 
from chain complexes
to simplicial $R$-modules, 
cf. e.~g. \cite\doldpupp\ (3.2 on p.~219).

When $G$ is substituted for $Y$, the resulting simplicial
object is a simplicial group  $EG$, and the diagonal injection
$G \to EG$ turns $EG$ into a
simplicial principal right (or left) $G$-space.
The simplicial manifold  $N(G,X)=EG\times_G X$
is what we refer to as the {\it homogeneous\/} Borel-construction.
The term  \lq\lq homogeneous\rq\rq\ is intended to hint at the fact that the 
formulation uses the group  structure only for 
the $G$-action on $EG$ and {\it not\/} for the simplicial structure
on $EG$; cf. \cite\maclaboo\ (IV.5 p.~119) where this distinction is 
discussed relative to the bar resolution.

Let $\VV$ be a right $(G,C\fra g)$-module.
The simplicial structure of $EG$ induces 
a cosimplicial structure on the $\VV$-valued de Rham complex
$$
\Cal A(EG,\VV),
\tag2.5.1
$$
and the degreewise right diagonal $(G,C\fra g)$-module structures
(1.7.1),
relative to the right translation $G$-action on $EG$ where $G$ is viewed
as a subgroup of $EG$ and relative to the $(G,C\fra g)$-module structure
on $\VV$, turn $\Cal A(EG,\VV)$ into a 
{\it cosimplicial object in the category of right\/} 
$(G,C\fra g)$-modules.
Our aim is to prove the following.

\proclaim{Theorem 2.5.2}
Given the $(G,C\fra g)$-module $\VV$, the dual standard construction
$\bold T(\VV)$ associated,
over the category 
$\roman{Mod}_{(G,C\fra g)}$,
with 
the monad $(\Cal T,\eta,\mu)$  
spelled out in {\rm (2.4)} above and the $(G,C\fra g)$-module
$\VV$, is naturally isomorphic,
as a cosimplicial object in the category of right $(G,C\fra g)$-modules,
to $\Cal A(EG,\VV)$.
\endproclaim

We begin with the preparations for the proof.
At the risk of 
making a mountain out of a molehill we first explain briefly
the right and left nonhomogeneous versions 
$(EG)^{\roman{right}}$ and $(EG)^{\roman{left}}$
of the Borel construction;
we need them {\it both\/} to arrive at consistent formulas at a 
later stage. Indeed,
the dual standard construction
$\bold T(\VV)$ associated, over the category 
$\roman{Mod}_{(G,C\fra g)}$, with 
the monad $(\Cal T,\eta,\mu)$ 
spelled out in {\rm (2.4)} above and the $(G,C\fra g)$-module
$\VV$, is more naturally identified with
$\Cal A((EG)^{\roman{left}}, \VV)$, cf. Proposition 2.5.7 below,
whereas, given the left $G$-manifold $X$, within the framework
of comonads and standard construction,
the ordinary Borel construction leads 
to the simplicial manifold $(EG)^{\roman{right}}\times_GX$.

The nonhomogeneous {\it right} ({\it left\/}) Borel-construction
arises as follows: Let $\roman{Smooth}$ be the category of smooth manifolds,
$\roman{Smooth}_G$ (${}_G\roman{Smooth}$) that of smooth right 
(left) $G$-manifolds, and let
$\Cal F \colon \roman{Smooth} \to \roman{Smooth}_G $ 
($\Cal F \colon \roman{Smooth} \to {}_G\roman{Smooth}$)
be the functor
which assigns to the
smooth manifold $Z$ the smooth right $G$-manifold $Z\times G$
(left $G$-manifold $G\times Z$), 
endowed with the obvious right 
(left) $G$-action induced by right (left) translation in $G$.
This functor is left adjoint to the forgetful functor 
$\square \colon \roman{Smooth}_G \to \roman{Smooth}$
($\square \colon {}_G\roman{Smooth} \to \roman{Smooth}$),
and the standard construction applied to the 
resulting comonad and the right (left) $G$-manifold $Z$ yields 
a simplicial manifold  $\Cal E(Z,G)$ ($\Cal E(G,Z)$)
endowed with a free right (left)
$G$-action. For $Z$ a point $o$, we will write
$(EG)^{\roman{right}} = \Cal E(o,G)$
($(EG)^{\roman{left}} = \Cal E(G,o)$).
This is the
{\it nonhomogeneous version\/} 
of the total simplicial $G$-object for $G$ in the category of right (left)
$G$-manifolds.

The various constructions are isomorphic as simplicial free right (left) 
$G$-manifolds.
Indeed, 
$(EG)^{\roman{right}}$ is the simplicial group having
the iterated semi-direct product
$$
(EG)^{\roman{right}}_n
= G \rtimes G \rtimes \ldots \rtimes G\ 
(n+1\ \text{copies of}
\  G) 
$$
as degree $n$ constituent,
with $G$-action from the right via right translation on the rightmost 
copy of $G$.
The {\it nonhomogeneous\/} face operators 
$\partial_j$ are given by
the expressions
$$
\aligned\partial_0(x_0,x_1,\dots, x_n) &= (x_1,\dots, x_n),
\\
\partial_j(x_0,x_1,\dots, x_n) &= 
(x_0,\dots, x_{j-2},x_{j-1}x_j,
x_{j+1},\dots, x_n)\ 
(1 \leq j \leq n)
\endaligned
\tag2.5.3.r
$$
and, likewise,
the {\it nonhomogeneous\/} degeneracy operators 
$s_j$ are given by
$$
s_j(x_0,x_1,\dots, x_n) = 
(x_0,\dots, x_{j-1},e,x_j,\dots, x_n)
\ (0 \leq j \leq n).
\tag2.5.4
$$
The notation $X$ being maintained for the left $G$-manifold at the 
beginning of the present subsection, let
$$
\Cal N(G,X) =(EG)^{\roman{right}}\times_G X.
$$
The associations
$$
(x_0,x_1,\dots, x_n)\longmapsto
(x_0x_1x_2\ldots x_n,
\ldots,x_{n-1}x_n, x_n),\ x_j \in G,
$$
as $n$ ranges over the natural numbers,
induce an isomorphism of simplicial principal right $G$-manifolds
$$
(EG)^{\roman{right}} \longrightarrow EG,
\tag2.5.5.r
$$ 
in fact, an isomorphism of 
simplicial groups, and thence an isomorphism
of simplicial manifolds from
$\Cal N(G,X)$ onto $N(G,X)=EG \times_GX$. 
In degree $n$, the inverse of (2.5.5.r) is 
plainly given by the association 
$$
(y_1,\ldots,y_n,x) \longmapsto 
(x_0,x_1,\dots, x_n)=
(y_1y^{-1}_2,y_2y_3^{-1},\ldots,y_{n-1}y_n^{-1}, y_nx^{-1},x).
$$
This association explains, in a somewhat more down to earth manner than the elegant categorical 
explanation on p. 107 of \cite\gsegatwo,
the formulas in \cite\duponone\ and \cite\gsegatwo,
cf. also p. 573 of \cite\kan,
for the projection
from  the homogeneous total simplicial
object $EG$, written there as $N\overline G$ 
(the nerve of a suitably defined category $\overline G$ associated with $G$) 
to the base 
$(EG)^{\roman{right}}\big/ G$ 
of the universal simplicial $G$-bundle, written out in
nonhomogeneous form so that the simplicial structure on the base
becomes more perspicuous.

Likewise,
$(EG)^{\roman{left}}$ is the simplicial group having
the iterated semi-direct product
$$
(EG)^{\roman{left}}_n
= G \ltimes G \ltimes \ldots \ltimes G\ 
(n+1\ \text{copies of}
\  G) 
$$
as degree $n$ constituent,
with $G$-action from the left via left translation on the leftmost 
copy of $G$.
The {\it nonhomogeneous\/} face operators 
$\partial_j$ are given by
the familiar expressions
$$
\aligned
\partial_j(x_0,x_1,\dots, x_n) &= 
(x_0,\dots, x_{j-2},x_{j-1}x_j,
x_{j+1},\dots, x_n)\ 
(0 \leq j < n)\\
\partial_n(x_0,x_1,\dots, x_n) &= (x_0,\dots, x_{n-1})
\endaligned
\tag2.5.3.l
$$
and 
the {\it nonhomogeneous\/} degeneracy operators 
$s_j$ are still given by (2.5.4).
The associations
$$
(x_0,x_1,\dots, x_n)\longmapsto
(x_0,x_0x_1,\ldots,x_0x_1x_2\ldots x_{n-1},x_0x_1x_2\ldots x_n),\ x_j \in G,
$$
as $n$ ranges over the natural numbers,
induce an isomorphism of simplicial principal left $G$-manifolds
$$
(EG)^{\roman{left}} \longrightarrow EG, 
\tag2.5.5.l
$$
in fact, an isomorphism of 
simplicial groups.
In degree $n$, the inverse of (2.5.5.l) is 
plainly given by the association 
$$
(x,y_1,\ldots,y_n) \longmapsto 
(x_0,x_1,\dots, x_n)=
(x,x^{-1}y_1, y^{-1}_1y_2, \ldots,y^{-1}_{n-1}y_n).
$$

Let $V$ be a right $G$-module. 
We will now consider $\Cal A^0((EG)^{\roman{left}}, V)$ as a cosimplicial
right $G$-module, the right $G$-module structure being the diagonal structure
relative to left $G$-translation on $(EG)^{\roman{left}}$ and the cosimplicial
structure being induced from the simplicial structure on 
$(EG)^{\roman{left}}$; we recall that the diagonal structure is given by the
association
$$
\Cal A^0((EG)^{\roman{left}}, V) \times G
\longrightarrow
\Cal A^0((EG)^{\roman{left}}, V),
\ (\alpha,x) \longmapsto \alpha\cdot x,
\tag2.5.6
$$
where $(\alpha\cdot x)y = (\alpha(xy))x$, $x \in G$, 
$y \in (EG)^{\roman{left}}$;
here $\alpha$ ranges over smooth functions from $(EG)^{\roman{left}}$ to $V$.
We remind the reader that we write various forgetful functors 
as $\square$.
Thus, given the $G$-representation $V$, the notation  $\square V$
refers to the vector space which underlies $V$, possibly endowed with 
trivial $G$-action.

\proclaim{Proposition 2.5.7} 
Let  
$\bold T^0(V)$ be the dual standard construction associated,
 over the category 
$\roman{Mod}_{G}$, with 
the monad $(\Cal T,\eta,\mu)$ 
spelled out in {\rm (2.2)} above and the $G$-module
$V$. Relative to the diagonal $G$-action on
$\Cal A^0((EG)^{\roman{left}}, V)$,
the morphism
$$
\Phi^0=(\varphi_0, \ldots ) 
\colon \Cal A^0((EG)^{\roman{left}}, V) \longrightarrow \bold T^0(V)
\tag2.5.8
$$
of graded $R$-modules which, in degree $n$, is given by the association
$$
\aligned
\varphi_n &\colon \Cal A^0(G^{\times(n+1)}, V)
\longrightarrow \Cal A^0(G^{\times(n+1)}, \square V),
\\
\varphi_n(\alpha)(x_0,\ldots,x_n) 
&=(\alpha(x_0,\ldots,x_n))\cdot x_0\cdot \ldots \cdot x_n,
\ x_0,\ldots,x_n\in G,
\endaligned
$$
where $\alpha$ ranges over smooth maps from $G^{\times(n+1)}$ to $V$,
is an isomorphism of cosimplicial right $G$-modules.
\endproclaim

\demo{Proof} This is certainly folk lore; a direct argument
comes down to a tedious but straightforward verification.
We leave the details to the reader. \qed
\enddemo

N.B. In view of $G$-equivariance, 
the isomorphism $\varphi_0$ is determined by the requirement 
$$
(\varphi_0(\alpha))(e)= \alpha(e),\ \alpha \colon G \to V.
$$

Let $\VV$ be a right $(G,C\fra g)$-module.
We will now consider $\Cal A((EG)^{\roman{left}}, \VV)$ as a cosimplicial
right $(G,C\fra g)$-module, 
the cosimplicial structure being induced from the simplicial structure on 
$(EG)^{\roman{left}}$ and
the right $(G,C\fra g)$-module structure being the {\it diagonal\/} structure
relative to {\it left\/} $G$-translation on $(EG)^{\roman{left}}$ and
the  right $(G,C\fra g)$-module structure on $\VV$ (beware:
in (2.5.1) above
that kind of structure was considered relative to {\it right\/} translation
on $EG$ via the explicit description (1.7.1)).

For intelligibility, we spell out this diagonal structure explicitly:
The right $G$-module structure
on  $\Cal A((EG)^{\roman{left}}, \VV)$ is given by the  extension 
of the association (2.5.6) above to the present situation. That is,
in a cosimplicial degree $n$, 
given  the
$\VV$-valued $p$-form $\alpha$ on $(EG)_n^{\roman{left}}$, the 
$p$-tuple $Z_1,\ldots,Z_p$ of vector fields on
$(EG)_n^{\roman{left}}$,
 and $x\in G$, let
$$
(\alpha\cdot x)(Z_1,\ldots,Z_p)
= (\alpha(xZ_1,\ldots,xZ_p))x;  
$$
here the notation $xZ_j$ ($1 \leq j\leq p$) refers to the induced
left $G$-action on the vector space of smooth vector fields
$\roman{Vect}((EG)_n^{\roman{left}})$ on $(EG)_n^{\roman{left}}$.
The right $G$-module structure
on  $\Cal A((EG)_n^{\roman{left}}, \VV)$ is given by the
pairing 
$$
\Cal A((EG)_n^{\roman{left}},\VV) \times G
\longrightarrow
\Cal A((EG)_n^{\roman{left}}, \VV),
\ (\alpha,x) \longmapsto \alpha\cdot x, \ x \in G,
$$
where $\alpha$ ranges over $\VV$-valued $p$-forms on 
$(EG)_n^{\roman{left}}$, for $p \geq 0$.
Furthermore, let $Y\in\fra g$;
given the  $\VV$-valued $p$-form $\alpha$,
the $p$-form $[\alpha,Y]$,
evaluated at the $p$-tuple $Z_1,\ldots,Z_p$ of vector fields on
$(EG)_n^{\roman{left}}$, is given by
$$
[\alpha,Y](Z_1,\ldots,Z_p) = \sum\alpha(Z_1,\ldots,[Y,Z_j],\ldots, Z_p)
+[\alpha(Z_1,\ldots,Z_p),Y].
$$
Likewise, let $z=sY$ where $Y\in \fra g$.
Let $\alpha=(\alpha_0,\ldots,\alpha_m)$ 
be an $m$-form in 
$\Cal A((EG)_n^{\roman{left}},\VV)$,
with components 
$\alpha_j\in \Cal A^j((EG)_n^{\roman{left}},\VV^{m-j})$ ($0 \leq j \leq m$). 
Now, for $0 \leq j \leq m-1$, given $q\in H$ and the vector fields
$Y_1,Y_2,\ldots,Y_j$ on $H$,
the value $(i_{Y_{\VV}}\alpha_j)(Y_1,Y_2,\ldots,Y_j)(q)$ is given by
$$
(i_{Y_{\VV}}\alpha_j)(Y_1,Y_2,\ldots,Y_j)(q)
=i_{Y_{\VV}}((\alpha_j(Y_1,Y_2,\ldots,Y_j))(q))=
((\alpha_j(Y_1,Y_2,\ldots,Y_j))(q))\cdot z.
$$
N.B. The value $(\alpha_j(Y_1,Y_2,\ldots,Y_j))(q)$ lies in $\VV^{m-j}$,
and $((\alpha_j(Y_1,Y_2,\ldots,Y_j))(q))\cdot z$ 
lies in $\VV^{m-j-1}$.
With this preparation out of the way,
$\alpha \cdot z$ is given by
$$
\alpha \cdot z= i_{Y_{\VV}}\alpha_0 -i_{Y}\alpha_1 +
i_{Y_{\VV}}\alpha_1 - i_Y\alpha_2 + \ldots +
i_{Y_{\VV}}\alpha_{m-1} - i_Y\alpha_m .
$$

\proclaim{Proposition 2.5.10} 
Let
$\bold T(\VV)$ be the dual standard construction associated,
 over the category 
$\roman{Mod}_{(G,C\fra g)}$,
 with 
the monad $(\Cal T,\eta,\mu)$ 
spelled out in {\rm (2.4)} above and the $(G,C\fra g)$-module
$\VV$.
Relative to the diagonal $(G,C\fra g)$-action on
$\Cal A((EG)^{\roman{left}},\VV)$,
the unique extension
$$
\Phi^{\roman{left}} \colon \Cal A((EG)^{\roman{left}}, \VV) \longrightarrow \bold T(\VV)
\tag2.5.11
$$
of the morphism $\Phi^0$ given as {\rm (2.5.8)}
above
is an isomorphism of cosimplicial right $(G,C\fra g)$-modules.
This extension is characterized as follows:
In a simplicial degree $n$, given the $\VV$-valued
$j$-form on $(EG)_n^{\roman{left}}$,
$Y_1,\ldots,Y_j \in \fra g$, and $z_j=sY_j \in s \fra g$,
$$
(\Phi\alpha) (Y_1,\ldots,Y_j)(e)
=(\alpha \cdot z_1\cdot \ldots \cdot z_j)(e) \in \VV.
$$
\endproclaim

We will now consider the simplicial left $G$-manifold
$(EG)^{\roman{left}}$ as a simplicial right $G$-manifold
via the pairing
$$
(EG)^{\roman{left}} \times G \longrightarrow (EG)^{\roman{left}},
\ (y,x) \longmapsto x^{-1}y,\ 
y \in (EG)^{\roman{left}}, x \in G.
\tag2.5.12
$$

\proclaim{Proposition 2.5.13}
The diffeomorphisms
$$
G^{\times (n+1)} \longrightarrow G^{\times (n+1)},\ 
(x_0,\ldots,x_n)\longmapsto 
(x_1,\ldots,x_n,x^{-1}_nx^{-1}_{n-1}\ldots x^{-1}_0)
\ (n \geq 0)
$$
induce an isomorphism
$$
(EG)^{\roman{left}} \longrightarrow (EG)^{\roman{right}}
\tag2.5.14
$$
of simplicial right $G$-manifolds.
\endproclaim

We can now complete the proof of Theorem 2.5.2:
The isomorphism (2.5.14), together with (2.5.11),
induces an isomorphism
$$
\Phi^{\roman{right}} \colon \Cal A((EG)^{\roman{right}}, \VV)
 \longrightarrow \bold T(\VV)
\tag2.5.15
$$
of cosimplicial right $(G,C\fra g)$-modules.
This isomorphism, in turn, combined with  the isomorphism 
induced by the isomorphism (2.5.5.r),
yields the desired isomorphism between (2.5.1)
and the dual standard construction $\bold T(\VV)$
associated with 
the monad $(\Cal T,\eta,\mu)$ 
under discussion  over the category $\roman{Mod}_{(G,C\fra g)}$
and the $(G,C\fra g)$-module $\VV$.
This proves Theorem 2.5.2.

Henceforth we shall no longer distinguish in notation between
the homogeneous version $EG$ and the nonhomogeneous versions
$(EG)^{\roman{left}}$ and $(EG)^{\roman{right}}$.

\smallskip
\noindent {\smc 2.6. Extension of Bott's decomposition lemma.\/} 
Let $\VV$ be a right $(G,C\fra g)$-module.
Application of the functor $\Cal A^0$ to $EG$ 
yields the cosimplicial algebra $\Cal A^0(EG)$;
likewise application of that functor to $EG$ and 
$\VV$ yields the cosimplicial
chain complex $\Cal A^0(EG,\VV)$,
and the latter inherits  a cosimplicial
differential graded $(\Cal A^0(EG))$-module structure. 
Just as for (2.5.1) above,
the degreewise right diagonal $(G,C\fra g)$-module structures
(1.7.2),
relative to the right translation $G$-action on $EG$ where $G$ is viewed
as a subgroup of $EG$ and relative to the $(G,C\fra g)$-module structure
on $\VV$, turn 
$\roman{Hom}^{\tau_{ E\fra g}}(\Lambda'_{\partial}[s E\fra g],\Cal A^0(EG, \VV))$
into a 
{\it cosimplicial object in the category of right\/} 
$(G,C\fra g)$-modules.

The decomposition of the functor $\Cal G_{(G,C\fra g)}$ spelled out in
Proposition 2.4.6 translates to a {\it decomposition\/} for the 
corresponding standard constructions.
The explicit description thereof leads to the
following:

\proclaim{Theorem 2.6.1}{\rm (Extended decomposition lemma)} The 
degreewise left trivialization of the tangent bundle of the
simplicial group
$EG$, that is, the morphism {\rm (1.6.3.1)}, 
evaluated degreewise,
yields an  
isomorphism
$$
\Phi\colon\Cal A(EG,\VV)
\longrightarrow
\roman{Hom}^{\tau_{E\fra g}}(\Lambda'_{\partial}[s E\fra g],\Cal
A^0(EG,\VV)) 
\tag2.6.2
$$
of cosimplicial right $(G,C\fra g)$-modules
from the cosimplicial chain
complex $\Cal A(EG,\VV)$
onto the differential graded cosimplicial diagonal object on the
right-hand side of {\rm (2.6.2)}.
\endproclaim

\demo{Proof} In a degree $p \geq 0$, the cosimplicial diagonal
object on the right-hand side of (2.6.2) comes down to
$$
\roman{Hom}^{\tau_{E\fra g}} (\Lambda'_{\partial}[s E\fra g],\Cal
A^0(EG,\VV))_p = \roman{Hom}^{\tau_{(E\fra g)_p}} (
\Lambda'_{\partial}[s (E\fra g)_p],\Cal A^0(EG,\VV)_p );
$$
since $(E\fra g)_p =\fra g^{p+1}$ and $(EG)_p=G^{p+1}$, 
in view of Proposition 1.7.3, the
isomorphism (1.6.3.1),
with $H=G^{p+1}$, $\fra h=\fra g^{p+1}$, and $V =
\VV$, identifies this cosimplicial diagonal object with the
$\VV$-valued de Rham complex of $G^{p+1}$. \qed
\enddemo

\smallskip
\noindent {\smc Remark 2.6.3.} 
For the special case where $\VV$ is 
the ground field $\Bbb R$ 
and $G$ finite dimensional, a version of the
isomorphism (2.6.2) (in a language different from ours) is given
in the Decomposition Lemma in \cite\bottone. 

\smallskip
\noindent
{\smc 2.7. Equivariant de Rham 
theory as a  differential Ext\/}.
As before, $G$ denotes a Lie group
(neither necessarily compact
nor necessarily finite dimensional) and $X$  a left
$G$-manifold. Application of the de Rham functor to
the {\it simplicial Borel construction\/} relative to $G$ and $X$
yields a cosimplicial differential graded algebra
whose total object defines the $G$-equivariant de Rham
cohomology of $X$. 

Exploiting the monad  $(\Cal T,\eta,\mu)$  over the category
$\roman{Mod}_{(G,C\fra g)}$ introduced in Subsection 2.4 above
we can now spell out
the de Rham theory replacement for Theorem 3.1 in
\cite\duaone\  (which, in turn, refers to ordinary cohomology).
For our purposes, this 
de Rham theory replacement
reduces equivariant
de Rham theory to ordinary  homological algebra
and thereby provides a high amount of flexibility.
In a sense we will explore this flexibility
in the rest of the paper.

\proclaim{Theorem 2.7.1}  The cosimplicial chain complex $\Cal
A(\Cal N(G,X))$ associated with the nonhomogeneous 
simplicial Borel 
construction $\Cal N(G,X)$
is canonically isomorphic to the
cosimplicial chain complex 
$\left|\bold T(\Cal A(X))\right|^{(G,C\fra g)}$,
the $(G,C\fra g)$-invariants of the
chain complex associated with the
dual standard construction
relative to the monad  $(\Cal T,\eta,\mu)$  over the category 
$\roman{Mod}_{(G,C\fra g)}$ and the
$(G,C\fra g)$-module $\Cal A(X)$.
Consequently
the $G$-equivariant de Rham cohomology of
$X$ is canonically isomorphic to the 
 differential
$\roman{Ext}_{((G,C \fra g);\chain)}(\Bobb R,\Cal A(X))$. 
\endproclaim

\demo{Proof} This is a consequence of Theorem 2.5.2,
combined with the observation, cf. Proposition 1.6.13, that
the projection $EG \times X \to N(G,X)$ identifies the cosimplicial
differential graded algebra 
$\Cal A(N(G,X))\cong\Cal A(NG,\Cal A(X))$ with the
cosimplicial differential graded subalgebra $\Cal A(EG,\Cal
A(X))^{(G,C\fra g)}$ of $G$-invariant $\Cal A(X)$-valued
horizontal forms on $EG$. \qed
\enddemo

Thus the category of $(G,C\fra g)$-modules serves, in the smooth category,
as a replacement for the
non-existent category of comodules relative to the de Rham
complex of $G$. Since the  differential
$\roman{Ext}_{((G,C \fra g);\chain)}(\Bobb R,\, \cdot\,)$
yields $G$-equivariant de Rham cohomology, we occasionally refer to
the associated infinitesimal object, that is, to the
relative differential
$\roman{Ext}_{(C\fra g,\fra g)}(\Bobb R,\, \cdot\,)$ isolated in
(2.3) above, as
{\it infinitesimal\/} $\fra g$-{\it equivariant cohomology\/}.

The present
approach to equivariant de Rham theory is completely formal; suitably rephrased
it works perfectly well in similar situations and yields 
e.~g. the
{\sl equivariant sheaf cohomology of a complex manifold relative to a 
holomorphic action
of a complex Lie group\/}.

\smallskip
\noindent
{\smc 2.8. The relative differential\/} $\roman{Tor}^{(C \fra g, \fra g)}$.
Here again the ground ring $R$ is a general commutative ring with 1.
Any pair of rings $(\Cal R,\Cal S)$ with $\Cal R \supset \Cal S$
gives rise to a  resolvent pair of 
categories, cf.
\cite \maclaboo\ (IX.6), that is to say, the functor
$\Cal F \colon \roman{Mod}_{\Cal S} \to  \roman{Mod}_{\Cal R}$ 
which assigns to the
right $\Cal S$-module $N$ the induced  $\Cal R$-module 
$\Cal F(N)=N\otimes_{\Cal R}\Cal S$
is left adjoint to the forgetful functor $\square 
 \colon \roman{Mod}_{\Cal R} \to  \roman{Mod}_{\Cal S}$.
Relative (co)homology is then defined and
calculated in terms of a relatively projective 
resolution in the sense of  \cite\hochsone.
Given a right $\Cal R$-module $N$, the
{\it standard construction\/} $\bold L(N)$
arising from $N$ and the {\it comonad\/}
$(\Cal L,\varepsilon,\delta)$
associated with the adjunction is a simplicial object
whose associated chain complex
$\big|\bold L(N)\big|$
coincides with the
standard relatively projective resolution of $N$ 
in the sense of  \cite\hochsone.

In our case where $(\Cal R,\Cal S)=(\roman U[C \fra g], \roman U[\fra g])$,
the functor 
$$
\Cal F \colon \chain_{\fra g} \longrightarrow \roman{Mod}_{C\fra g}
$$ 
assigns to the right $\fra g$-chain complex $N$ the
totalized CCE complex 
$$
\Cal FN = N \otimes_{\fra g} \roman U[C \fra g] \cong
N \otimes_{\tau_{\fra g}}\Lambda'_{\partial}[s\fra g]
$$
calculating the Lie algebra homology of $\fra g$ with coefficients
in $N$ (suitably interpreted relative to the chain complex structure on $N$), 
and the differential
$\roman{Tor}^{(C \fra g, \fra g)}$-
and
$\roman{Ext}_{(C \fra g, \fra g)}$-functors are defined
on the category of $(C\fra g)$-modules.
The functor 
$\roman{Ext}_{(C \fra g, \fra g)}$ is the same as that introduced
in Subsection 2.3 above.

\smallskip
\noindent
{\smc 2.9. Lie-Rinehart algebras and Lie algebroids.\/} 
Even though this is not relevant later in the paper,
we spell out briefly a generalization of the situation of (2.3) to illustrate the flexibility
of the present formal approach.

Let $(A,L)$ be a Lie-Rinehart algebra and let $\roman U(A,L)$ be the universal algebra
associated with $A$ and $L$; see e.~g. \cite\extensta\ or \cite\lradq\ for details.
Let 
$\Cal F \colon {}_A\roman{Mod} \to  {}_{\roman U(A,L)}\roman{Mod}$ 
be the functor which assigns to the
$A$-module $M$ the induced  $\roman U(A,L)$-module 
$$
\Cal F(M)=\roman U(A,L)\otimes_A M
$$ 
and, as before, denote the forgetful functor by $\square 
\colon {}_{\roman U(A,L)}\roman{Mod}\to  {}_A\roman{Mod}$.
The resulting relative $\roman{Ext}_{(\roman U(A,L),A)}$ is precisely the
cohomology
theory introduced in \cite\rinehone\ by means of a generalized
CCE 
complex adapted to Lie-Rinehart algebras.
In particular, when $L$ is projective as an $A$-module,
the relative $\roman{Ext}_{(\roman U(A,L),A)}$ is an absolute
$\roman{Ext}_{\roman U(A,L)}$.
In particular, when $L$ is the $(R,A)$-Lie algebra
$D_{\{\,\cdot\, , \,\cdot \,\}}$ 
associated with a Poisson structure
$\{\,\cdot\, , \,\cdot \,\}$ on $A$ \cite\poiscoho,
the relative $\roman{Ext}_{(\roman U(A,L),A)}(A,A)$ coincides with 
the {\it Poisson cohomology\/} of $A$, cf. \cite\poiscoho.

The situation in (2.3) does not extend directly to Lie-Rinehart algebras
since the {\it cone on a Lie-Rinehart algebra is ill-defined\/}.
Indeed, given the $(A,L)$-module $M$, for $a \in A$ and $\alpha \in L$,
the operations $i$ of contraction and
$\lambda$ of Lie derivative satisfy the familiar 
identity
$$
\lambda_{a\alpha}(\omega)= a\lambda_{\alpha}(\omega) + 
da\cup i_{\alpha}(\omega)\ (a \in A, \ \alpha \in L,\ 
\omega \in \roman{Alt}_A(L,M))
\tag2.9.1
$$
involving the term $da\cup i_{\alpha}(\omega)$ which does not arise
for an ordinary Lie algebra.
Thus the values of the
functor $\Cal G$ on the category 
${}_{(A,L)}\roman{Mod}$ of $(A,L)$-modules
which assigns the {\it Rinehart complex\/} (generalized de Rham complex)
$$
\Cal G(V)=  (\roman{Alt}_A(L, V),d)
$$
to the $(A,L)$-module $V$ lie in a certain category $\Cal M$ of $(A,L)$-modules
which are also endowed with an action of the {\it ordinary cone\/}
$CL$ on $L$ in the category of Lie algebras 
(beware: not Lie-Rinehart algebras),
subject to certain identities including (2.9.1).
The resulting adjunction defines a monad and the corresponding
dual standard construction yields a relative differential graded Ext.
We believe that this is a formally correct approach to phrase
developments like the BRST-complex, the variational bicomplex, 
and the Noether identities.
The corresponding constructions for Lie algebroids can then presumably 
be globalized via Lie groupoids by a suitable
comonadic construction.
We hope to come back to these issues at another occasion.

We conclude our discussion with two more examples which illustrate the universalness
of the present relative  homological algebra approach to equivariant cohomology.
While these examples are not strictly needed for the rest of the paper,
they will make it clear that the dual standard construction
in (2.2) above which defines differentiable cohomology 
as well as that
which defines the differential $\roman{Ext}_{((G,C \fra g);\chain)}$
in (2.4) above are both versions of completed {\it cobar\/} constructions.

\smallskip
\noindent
{\smc 2.10. Rational cohomology of algebraic groups.\/} 
Let $k$ be a field,  let ${}_k\roman{Vect}$ be the category of
$k$-vector spaces,
let $G$ be an algebraic group defined over $k$, and let
$k[G]$ be the coordinate ring of $G$. The group structure turns $k[G]$ into a Hopf algebra.
A {\it rational\/} $G$-representation is, by definition, a
$k[G]$-comodule. Let ${}_{k[G]}\roman{Comod}$ be the category of
$k[G]$-comodules or, equivalently, rational $G$-representations, and let
$\Cal G\colon {}_k\roman{Vect} \to {}_{k[G]}\roman{Comod}$ be the functor
which assigns to the $k$-vector space $V$ the induced comodule 
$k[G]\otimes V$. This functor is right adjoint to the forgetful
functor  $\square \colon {}_{k[G]}\roman{Comod}\to {}_k\roman{Vect}$
whence the two functors define a monad,
and the resulting dual standard construction yields the appropriate
{\it cobar construction\/} which defines the
$\roman{Cotor}_{k[G]}(\,\cdot \, , \, \cdot \,)$ and in particular
the {\it rational cohomology\/} of $G$ with coefficients in a rational $G$-module.

\smallskip
\noindent
{\smc 2.11. Equivariant de Rham cohomology for algebraic varieties.\/} 
Let $k$ be a field and let $G$ be an algebraic group defined over $k$.
The algebraic de Rham algebra $\Cal A[G]$ of $G$ acquires a differential graded Hopf algebra structure.
Let ${}_{\Cal A[G]}\roman{Comod}$ be the category of
$\Cal A[G]$-comodules, ${}_k\Cal C$ that of $k$-chain complexes,
and let $\Cal G\colon {}_k\Cal C \to {}_{\Cal A[G]}\roman{Comod}$ be the functor
which assigns to the $k$-chain complex $V$ the induced comodule 
$\Cal A[G]\otimes V$, appropriately totalized. This functor is right 
adjoint to the forgetful
functor  $\square \colon {}_{\Cal A[G]}\roman{Comod}\to {}_k\Cal C$
whence the two functors define a monad,
and the resulting dual standard construction yields the appropriate
{\it cobar construction\/} which defines the functor
$\roman{Cotor}_{\Cal A[G]}(\,\cdot \, , \, \cdot \,)$.
Given a non-singular $G$-variety $X$,
the algebraic de Rham algebra $\Cal A[X]$ acquires an
$\Cal A[G]$-comodule structure, and the algebraic $G$-equivariant
de Rham cohomology of $X$ is given by
$\roman{Cotor}_{\Cal A[G]}(k,\Cal A[X])$.

\beginsection 3. Infinitesimal equivariant (co)homology

In this section we will explore the 
{\it infinitesimal equivariant cohomology\/}
functor $\roman{Ext}_{(C \fra g, \fra g)}(R,\,\cdot \,)$
by means of
homological algebra techniques.

The ordinary Weil algebra of a Lie algebra 
was introduced as an object which arises from abstraction of the
operations of contraction and Lie derivative
and serves as the principal tool for the description of 
the Chern-Weil map and of the Weil and Cartan models for equivariant cohomology
relative to a compact Lie group

In \cite\cartanei\ (Ex. XIII.14),
the CCE resolution is denoted by 
$V(\fra h)$ and the
universal differential graded algebra
$\roman U[C\fra h]$ which, as recalled above,
reproduces the CCE resolution,
is written as $W(\fra h)$. 
Is the usage of the letter $W$ just a 
notational incidence
or, at the time, was the 
notation $W$ intended to hint at 
the relationship
with the Weil algebra 
we are about to explain?

\smallskip
\noindent
{\smc 3.1. The Weil coalgebra.\/}
We will show that a classical construction in \cite\hochsone, 
adapted to our situation,
leads to what we refer to as the {\it ordinary Weil coalgebra\/}.
We have introduced that Weil coalgebra in
in \cite\extensta.

As before, let $R$ be a commutative ring
and $\fra \hgg$  an $R$-Lie algebra,
which we suppose throughout to be projective as an $R$-module.
As noted earlier,
the CCE coalgebra 
$\SSS_{\partial}[sC\fra
g]$ is defined for the differential graded Lie algebra
$C\fra \hgg$ (the cone on $\fra g$).
We will write
this  differential
graded coalgebra as
 $W'[\fra \hgg] = \Sigm_{\partial}'[sC\fra \hgg]$.
Thus $d+\partial$ turns $W'[\fra \hgg]$ into a
differential graded coalgebra which we refer to as the {\it Weil coalgebra\/}
of $\fra \hgg$. The dual $\roman{Hom}(W'[\fra \hgg],R)$ is 
the ordinary {\it Weil
algebra\/} of $\fra \hgg$.
For later reference,
we will now spell out some additional structure on the Weil coalgebra.

\smallskip\noindent
{(3.1.1)} {\sl Since, as a graded $R$-module,
$sC\fra \hgg = s^2 \fra \hgg \oplus s \fra \hgg$,
as a graded coalgebra,
the Weil coalgebra decomposes canonically as}
$$
W'[\fra \hgg] \cong \SSS[s^2 \fra \hgg] \otimes \Lambda' [s \fra \hgg].
$$
{(3.1.2)}
{\sl The underlying graded Lie algebra of $C\fra \hgg$
acts on $\SSS[s^2\hgg]$ in a canonical way through the
projection from $C\hgg$ to $\hgg$
(not a differential graded projection).
The graded $(C\hgg)$-module structures on
$\Lambda' [s\hgg]$ and $\SSS[s^2\hgg]$ combine to a  graded right
$(C\hgg)$-module structure
$$
W'[\hgg] \otimes C\hgg
@>>>
W'[\hgg]
\tag3.1.2.1
$$
which is, in fact, a differential
graded right $(C\hgg)$-module structure since the construction
of the differential graded CCE coalgebra is functorial
in the differential graded Lie algebra variable.}

\noindent {(3.1.3)}
{\sl Relative to the $\fra \hgg$-action
on the right of $\SSS[s^2\fra \hgg]$ coming from
the adjoint action of $\fra \hgg$ on itself,
$$
(W'[\fra \hgg],\partial)=
\SSS[s^2\fra \hgg] \otimes_{\tau_{\fra \hgg}}
\Lambda_{\partial}'[s \fra \hgg]\cong
\SSS[s^2\fra \hgg] \otimes_{\fra \hgg}\roman U[C\fra \hgg],
\tag3.1.3.1
$$
that is, $(W'[\fra \hgg],\partial)$ 
is precisely the standard complex computing
the Lie algebra homology of $\fra \hgg$ with values
in $\SSS[s^2\fra \hgg]$, viewed as a right 
$\fra \hgg$-module.}

\noindent
{(3.1.4)} {\sl Let $
\tau^{\SSS}\colon
\SSS[s^2\fra \hgg] \longrightarrow 
\Lambda[s \fra \hgg]
$
be the standard universal twisting cochain.
Relative to
the decomposition
$\SSS[s^2\fra \hgg] \otimes \Lambda[s \fra \hgg]$ of $W'[\fra \hgg]$,
the differential $d$ is the operator 
$\partial^{\tau^{\SSS}}= - ( \tau^{\SSS} \cap\,\cdot \, ) 
\colon \SSS[s^2\fra \hgg] \otimes \Lambda[s \fra \hgg]
\to \SSS[s^2\fra \hgg] \otimes \Lambda[s \fra \hgg],
$
so
$\Lambda[s \fra \hgg]$ is viewed as 
{\it fiber\/} and $\SSS[s^2\fra \hgg]$ as {\it base\/}
of the corresponding twisted tensor product
whence, as chain complexes,}
$$
(W'[\fra \hgg],d)=\left(
\SSS[s^2\fra \hgg] \otimes \Lambda[s \fra \hgg]
,d\right)
=\SSS[s^2\fra \hgg] \otimes_{\tau^{\SSS}}
\Lambda[s \fra \hgg].
\tag3.1.4.1
$$ 

\noindent
{(3.1.5)}
{\sl The graded right $(C\hgg)$-module structures on
$\Lambda' [s\hgg]$ and $\SSS[s^2\hgg]$ combine to a 
differential graded right
$(C\hgg)$-module structure
on
$(W'[\fra \hgg],\partial)$ as well
(on $W'[\fra \hgg]$ endowed merely with the operator $\partial$), 
that is,
the pairing {\rm (3.1.2.1)}, written out in the form
$$
\left(\SSS[s^2\fra \hgg] \otimes_{\tau_{\fra \hgg}}
\Lambda_{\partial}'[s \fra \hgg]\right)
 \otimes C\hgg
@>>>
\SSS[s^2\fra \hgg] \otimes_{\tau_{\fra \hgg}}
\Lambda_{\partial}'[s \fra \hgg]
\tag3.1.6
$$
is compatible with the differentials\/}. This reflects the
familiar fact that the effect of the adjoint action 
on Lie algebra homology is trivial, cf. (1.3) above.

Henceforth we will denote by  $\SSS_{2k}[s^2\fra g]$ 
the
$k$-th
homogeneous constituent of the graded symmetric 
coalgebra $\SSS[s^2\fra g]$.
In \cite\hochsone, 
for a pair $(\fra a, \fra b)$ of ordinary Lie algebras,
Hochschild has introduced
an acyclic relatively projective CCE complex
which yields the {\it relative Lie algebra cohomology of the pair
$(\fra a, \fra b)$ in the sense of Chevalley-Eilenberg\/} \cite\cheveile.
This CCE complex arises by abstraction from the situation of
the invariant de Rham complex of a homogeneous space of compact 
connected Lie groups.
Hochschild has furthermore shown that, when $\fra b$ is reductive in $\fra a$,
that CCE complex is a relatively projective 
resolution of the ground ring, that is, that CCE complex
admits a $\fra b$-equivariant
contracting homotopy.
The literal translation of that CCE construction,
to the pair $(C\fra g,\fra g)$ of differential graded Lie algebras,
yields the following:

\proclaim{Proposition 3.1.7} The Weil coalgebra
$W'[\fra g]$, written out in the form
$$
\ldots @>{d}>>
\SSS_{2p}[s^2 \fra g]
\otimes_{\tau_{\fra g}} \Lambda'_{\partial}[s \fra g]
@>{d}>> \ldots  
@>{d}>>
\SSS_2[s^2 \fra g]
\otimes_{\tau_{\fra g}} 
\Lambda'_{\partial}[s \fra g]
@>{d}>>\Lambda'_{\partial}[s \fra g],
\tag3.1.8
$$
is a relatively
projective 
complex
in the category of right
$(\roman U[C\fra g])$-modules which, augmented
by the obvious augmentation map
$\varepsilon \colon \Lambda'_{\partial}[s \fra g] \to R$
(the counit of $\Lambda'_{\partial}[s \fra g]$),
is an acyclic complex.
Here, for each $p \geq 1$, the right
$(\roman U[C\fra g])$-module structure is
the obvious one on
$$
\SSS_{2p}[s^2 \fra g] \otimes _{\tau_{\fra g}} \Lambda'[s \fra g] 
\cong 
\SSS_{2p}[s^2 \fra g] \otimes _{\roman U[\fra g]} \roman U[C\fra g].
$$
\endproclaim

\demo{Proof} In view of the preparatory steps
(3.1.1)--(3.1.5) this is straightforward and left to the reader. \qed
\enddemo

We do {\it not\/} claim that (3.1.8) 
has a $\fra g$-linear contracting homotopy.
Thus we do {\it not\/} assert that (3.1.8) 
is a relatively
projective 
resolution of $R$
in the category of right
$(\roman U[C\fra g])$-modules.

To spell out the appropriate structure that 
the Weil coalgebra
$W'[\fra g]$ acquires,
we return to the situation of (2.8) above:
Thus, consider
a pair of algebras $(\Cal R,\Cal S)$ 
with $\Cal R \supset \Cal S$ and 
suppose that $\Cal R$ is actually an augmented {\it differential\/} graded
algebra and that $\Cal S$ is an augmented differential graded
subalgebra with zero differential.
Given a chain complex $V$ we denote by
$V^{\sharp}$  the graded $R$-module underlying $V$.
Let $\overline M$ be a differential graded right $\Cal S$-module.
We will refer to an augmented differential graded right $\Cal R$-module
$M$ whose underlying graded $\Cal R$-module $M^{\sharp}$ is an induced
graded module of the kind
$M^{\sharp}=\overline M^{\sharp} \otimes_{\Cal S}\Cal R^{\sharp}$
as a  {\it construction\/} for $\Cal R$ {\it relative\/} to $\Cal S$
provided the induced isomorphism
$$
\overline M \longrightarrow  M \otimes_{\Cal R}R
$$
of graded $R$-modules is an isomorphism of
chain complexes.
When $\Cal S$ is the ground ring this notion of construction
comes down to the usual notion of 
construction in the sense of H. Cartan, cf. \cite\mooretwo.

Proposition 3.1.7 plainly says that the Weil coalgebra $W'[\fra g]$ is 
an $R$-acyclic, even $R$-contractible construction for 
$\roman U[C \fra g]$  relative to $\roman U[\fra g]$.

\proclaim{Theorem 3.1.9} Suppose that the ground ring $R$ is a field
of characteristic zero
and that $\fra g$ is reductive.
For any right $(C \fra g)$-module $N$,
the relative  
differential graded
$\roman{Ext}_{(C \fra g, \fra g)}(R,N)$
is the homology of the chain complex
$$
\roman{Hom}(W'[\fra g],N)^{C \fra g}.
\tag3.1.10
$$
Likewise for any left $(C \fra g)$-module $M$,
the relative  differential graded 
$\roman{Tor}^{(C \fra g, \fra g)}(R,M)$
is the homology of the chain complex
$$
W'[\fra g]\otimes_{C \fra g}M . 
\tag3.1.11
$$
\endproclaim

We will prove Theorem 3.1.9 in (3.3.6) and in (3.4) below.
Under the circumstances of Theorem 3.1.9, we will refer to 
(3.1.10) as the {\it Weil\/} model for
the differential graded
$\roman{Ext}_{(C \fra g, \fra g)}(R,N)$. In the language of ordinary 
differential geometry, 
(3.1.10) consists of the {\it basic\/} elements of
$\roman{Hom}(W'[\fra g],N)$: indeed, 
$C \fra g= (s\fra g)\rtimes \fra g$,
the invariants relative to the constituent $\fra g$
are the {\it invariant\/} elements 
(in the usual sense) and
the invariants relative to the constituent $s\fra g$
are the {\it horizontal\/} elements whence
the $(C \fra g)$-invariants are the elements which are
horizontal and invariant.
In Section 6 below we shall show
that a variant of (3.1.10) yields the 
familiar {\it Weil\/} model for 
the equivariant cohomology relative to a finite-dimensional
compact connected Lie group.

Recall that, for purely formal reasons, 
$\roman{Ext}_{(C \fra g, \fra g)}(R,R)$
acquires a graded commutative algebra structure.

\proclaim{Corollary 3.1.12}
When the ground ring $R$ is a field of characteristic zero and
when $\fra g$ is reductive,
as a graded commutative algebra, 
$\roman{Ext}_{(C \fra g, \fra g)}(R,R)$ is 
canonically isomorphic to 
the algebra $\roman{Hom}(\SSS[s^2\fra \hgg],R)^{\fra g}$  of
$\fra g$-invariants of 
the algebra $\roman{Hom}(\SSS[s^2\fra \hgg],R)$,
that is,
to the algebra of
$\fra g$-invariants of 
the symmetric algebra $\Sigm[s^{-2} \fra g^*]$
on the double desuspension $s^{-2} \fra g^*$ of the dual
$\fra g^*$ of $\fra g$. 
Likewise,
$\roman{Tor}^{(C \fra g, \fra g)}(R,R)$ then
acquires a graded coalgebra structure
and, as a graded coalgebra,
is canonically isomorphic to 
the coalgebra $\SSS[s^2\fra \hgg]\otimes_{\fra g}R$  of
$\fra g$-coinvariants of 
the graded coalgebra $\SSS[s^2\fra \hgg]$.
\qed
\endproclaim

\noindent
{\smc Remark 3.1.13.\/}
{\sl 
When $\fra g$ is abelian, the Weil coalgebra $W'[\fra g]$,
written out as a chain complex as in {\rm (3.1.8)}, plainly comes down to the ordinary
Koszul resolution of the ground ring 
in the category of 
$(\Lambda[s\fra g])$-modules.\/}

\smallskip
\noindent
{\smc 3.2. The relative  bar resolution.\/}
Given a right
$(C\fra g)$-module, 
this resolution is the {\it standard resolution\/} of
that module
associated with the comonad 
$(\Cal L,\varepsilon,\delta)$ mentioned in  
Subsection 2.8 above; cf. \cite\duskinon\ 
and \cite\maclbotw\ 
for this notion of standard resolution.

Let $\beta_{\roman U[\fra g]} (\roman U[C \fra g])$
denote the two-sided simplicial bar resolution 
of $\roman U[C \fra g]$  
{\it relative\/} to the category of 
$(\roman U[\fra g])$-modules;
this is a simplicial $(\roman U[C \fra g])$-bimodule.
Normalization yields the associated
two-sided {\it normalized relative bar resolution\/}
$ \roman B_{\roman U[\fra g]} (\roman U[C \fra g])$
of $\roman U[C \fra g]$ in the 
category
of $(\roman U[C \fra g])$-bimodules,

Let $M$ be a right $(\roman U[\fra g])$-module.
The simplicial 
{\it right\/} $(\roman U[C \fra g])$-module 
$$
\beta_{\roman U[\fra g]} 
(M,\roman U[C \fra g],\roman U[C \fra g])
=
M \otimes_{\roman U[C \fra g]} 
\beta_{\roman U[\fra g]} (\roman U[C \fra g])
\tag3.2.1
$$
is the
simplicial  bar resolution of $M$ 
in the category of {\it right\/}
$(\roman U[C \fra g])$-modules
{\it relative\/} to the category of 
$(\roman U[\fra g])$-modules. 
This is precisely the standard resolution of
$M$ relative to the comonad 
$(\Cal L,\varepsilon,\delta)$.
The normalized chain complex
of the simplicial object
$M \otimes_{\roman U[C \fra g]} 
\beta_{\roman U[\fra g]} (\roman U[C \fra g])$ 
 yields
the 
{\it normalized relative bar resolution\/}
$ \roman
B_{\roman U[\fra g]} (M,\roman U[C \fra g],\roman U[C \fra g])$
of $M$ in the 
category
of right $(\roman U[C \fra g])$-modules,
plainly  a  
resolution
of $M$ in the category of
right $(\roman U[C \fra g])$-modules that is
projective relative to the category of 
$(\roman U[\fra g])$-modules.
In particular, when $M=R$, viewed as a trivial
right $(\roman U[C \fra g])$-module,
$R \otimes_{\roman U[C \fra g]} \roman
B_{\roman U[\fra g]} (\roman U[C \fra g])$
leads to the 
{\it normalized relative bar resolution\/}
$ \roman
B_{\roman U[\fra g]} (R,\roman U[C \fra g],\roman U[C \fra g])$
of $R$ in the 
category
of right $(\roman U[C \fra g])$-modules.
Our notation for the functor $\roman B_{\roman U[\fra g]}$ etc.
is that in
\cite\gugenmay, with $\roman U[\fra g]$ substituted for the ground ring.
With the notation $\roman I\Lambda$ for the 
augmentation ideal of the exterior algebra 
$\Lambda$,
the resolution 
$\roman B_{\roman U[\fra g]} (R,\roman U[C \fra g],\roman U[C \fra g])$
has the form
$$
\ldots
@>{d}>>
\left((\roman I\Lambda[s\fra g])^{\otimes n}
\otimes \Lambda[s\fra g],\partial\right)
@>d>> \ldots @>d>>
\left(\roman I\Lambda[s\fra g]
\otimes \Lambda[s\fra g]
,\partial\right)
@>d>>(\Lambda[s\fra g],\partial)
\tag3.2.2
$$
where $d$ is the differential arising from
the operation of normalization.
Here the (differential graded) right
$(\roman U[C\fra g])$-module structure 
is induced by the
obvious pairing map
$$
(\Lambda[s\fra g])^{\otimes n}
\otimes \Lambda[s\fra g]
\otimes (\roman U[\fra g]\odot \Lambda[s\fra g])
\longrightarrow
(\Lambda[s\fra g])^{\otimes n}
\otimes \Lambda[s\fra g]
$$
coming from the obvious diagonal action
of $\roman U[\fra g]$ on the right of
$(\Lambda[s\fra g])^{\otimes n}
\otimes \Lambda[s\fra g]$ 
and from the right multiplication action of
$\Lambda[s\fra g]$ on itself.

To spell out the operator $\partial$ explicitly, 
we note that the augmentation map
$\varepsilon$ of $\Lambda[s\fra g]$
yields the counit
$$
\varepsilon\colon
\Lambda'_{\partial}[s\fra g]
\longrightarrow
R
$$
of 
the CCE coalgebra
$\Lambda'_{\partial}[s\fra g]$ of $\fra g$ as well;
this kernel is plainly a differential graded
cocommutative coalgebra
without counit and without coaugmentation
and, whenever this coalgebra structure 
is under discussion,
we use the notation $\roman I\Lambda'_{\partial}[s\fra g]$
for the kernel of $\varepsilon$
and refer to it
as the {\it augmentation coalgebra\/}
of $\Lambda'_{\partial}[s\fra g]$.
For $n \geq 1$,
the iterated twisted tensor product
$$
(\roman I\Lambda'_{\partial}[s\fra g])
\otimes_{\tau_{\fra g}}
(\roman I\Lambda'_{\partial}[s\fra g])
\otimes_{\tau_{\fra g}}
\ldots \otimes_{\tau_{\fra g}}
(\roman I\Lambda'_{\partial}[s\fra g]) 
\tag3.2.3
$$
of $n$ copies of
$\roman I\Lambda'_{\partial}[s\fra g]$
is still defined, e.~g. as a subcomplex of
the twisted tensor  product 
$$
\Lambda'_{\partial}[s\fra g]\otimes_{\tau_{\fra g}}
\Lambda'_{\partial}[s\fra g]\otimes_{\tau_{\fra g}}
\ldots \otimes_{\tau_{\fra g}}
\Lambda'_{\partial}[s\fra g] 
$$
of $n$ copies of $\Lambda'_{\partial}[s\fra g]$,
and 
the chain complex
$\left((\roman I\Lambda[s\fra g])^{\otimes n}
\otimes \Lambda[s\fra g],\partial\right)
$
in the resolution (3.2.2)
amounts to 
the  iterated twisted tensor product
$$
\left((\roman I\Lambda'_{\partial}[s\fra g])
\otimes_{\tau_{\fra g}}
(\roman I\Lambda'_{\partial}[s\fra g])
\otimes_{\tau_{\fra g}}
\ldots \otimes_{\tau_{\fra g}}
(\roman I\Lambda'_{\partial}[s\fra g])\right)
\otimes_{\tau_{\fra g}}
\Lambda'_{\partial}[s\fra g]
\tag3.2.4
$$
of $n$ copies of
$\roman I\Lambda'_{\partial}[s\fra g]$
with a single copy 
of 
$\Lambda'_{\partial}[s\fra g]$.

Following Mac Lane \cite\maclaboo,
we will use the notation $\left|\, \cdot \,\right|^{\bullet}$
for the {\it condensation\/} functor.
Thus
{\it condensation\/} transforms the resolution
(3.2.2) 
into the corresponding {\it construction\/}
$$
\left|R \otimes_{\roman U[C \fra g]} 
\beta_{\roman U[\fra g]} (\roman U[C \fra g])
\right|^{\bullet} =\left(
\roman T'[s\roman I\Lambda[s\fra g]]
 \otimes\Lambda[s \fra g],
d \otimes \roman{Id}_{\Lambda[s \fra g]}
-\tau^{\rbar}\cap\,\cdot \, + \partial\right)
$$
in the relative sense explained above, and we will use the notation
$$
\roman B\Lambda_{\partial}[s \fra g]
=\left|R \otimes_{\roman U[C \fra g]} 
\beta_{\roman U[\fra g]} (\roman U[C \fra g])
\right|^{\bullet} 
\tag3.2.5
$$
Thus $\roman B\Lambda_{\partial}[s \fra g]$ is
a  construction for $\roman U[C \fra g]$  relative to $\roman U[\fra g]$.
The ordinary bar construction contracting homotopy
is a $\fra g$-linear contracting homotopy for this construction.

Pushing a bit further, we observe that, 
for each $n \geq 1$,
via the obvious morphism
$$
s^{\otimes n}
\colon
\left(\roman I\Lambda[s\fra g]\right)^{\otimes n}
\longrightarrow
\left(s\roman I\Lambda[s\fra g]\right)^{\otimes n},
$$
the operator $\partial$
given by the chain complex (3.2.3)
induces an operator 
$$
\partial
\colon
\left(s\roman I\Lambda[s\fra g]\right)^{\otimes n}
\longrightarrow
\left(s\roman I\Lambda[s\fra g]\right)^{\otimes n}
$$
(where the notation $\partial$ is abused again).
These operators assemble to an operator 
$\partial$ on the graded
tensor coalgebra
$\roman T'[s\roman I\Lambda[s\fra g]]$, and we 
will write
$$
\rbar \Lambda_{\partial}[s\fra g]=
\left(\roman T'[s\roman I\Lambda[s\fra g]],
\partial\right).
$$
The object
$
\rbar \Lambda_{\partial}[s\fra g]
$
carries an obvious right $\fra g$-module structure, 
and
$$
\left(
\roman T'[s\roman I\Lambda[s\fra g]]
 \otimes\Lambda[s \fra g],
\partial\right) =
\rbar \Lambda_{\partial}[s\fra g]
\otimes_{\tau_{\fra g}}\Lambda'_{\partial}[s\fra g].
$$ 
Thus, relative to the operator $\partial$,
the object under discussion
appears as
a twisted tensor product
relative to
$\Lambda'_{\partial}[s\fra g]$
as
{\it base\/} and $\rbar \Lambda_{\partial}[s\fra g]$ 
as {\it fiber\/};
cf. (3.14) above where the corresponding twisted tensor
product decomposition is spelled out relative to the differential
$d$.

It is obvious that, when $\fra g$ is abelian,
the normalized relative
resolution comes down to the ordinary
bar resolution of $R$ in the category
of right $(\Lambda[s\fra g])$-modules.

\proclaim{Theorem 3.2.6} 
Relative to the tensor product coalgebra structure
on \linebreak
$\roman T'[s\roman I\Lambda[s\fra g]]
 \otimes\Lambda'[s \fra g]$,
the condensed object $\roman B\Lambda_{\partial}[s \fra g]$
(cf. {\rm (3.2.5)})
is a differential graded coalgebra.
\endproclaim

\demo{Proof} This comes down to a tedious verification.
The statement is also a consequence of Theorem 4.5 below, see Remark 4.6.
\enddemo

\noindent{\smc 3.3. Comparison between the Weil coalgebra and the 
relative bar construction.\/}
The bracket on $\fra g$ being momentarily ignored,
let $\tau^{\SSS}\colon\SSS[s^2 \fra g] \to 
\Lambda[s\fra g]$ be
the obvious acyclic twisting cochain; its adjoint
is the canonical injection 
$$
\overline{\tau^{\SSS}}\colon
\SSS[s^2 \fra g] @>>> \rbar\Lambda[s \fra g]
$$
of differential graded coalgebras.
Let
$$
\iota= \overline{\tau^{\SSS}} \otimes\roman{Id}
\colon \SSS[s^2 \fra g]\otimes_{\tau^{\SSS}}
\Lambda[s \fra g] 
\longrightarrow
\rbar\Lambda[s \fra g] 
\otimes_{\tau^{\rbar}}\Lambda[s \fra g].
$$
This is a morphism
$$
\iota\colon
\SSS[s^2 \fra g]\otimes_{\tau^{\SSS}}
\Lambda[s \fra g]  \longrightarrow \roman B\Lambda[s \fra g]
\tag3.3.1
$$
of differential graded coalgebras.
The canonical comparison 
between 
$W'[\fra g]$ and 
$\roman B\Lambda_{\partial}[s \fra g]$
is achieved
by the following.

\proclaim{Theorem 3.3.2} The 
perturbation $\partial$
determined by the Lie bracket on $\fra g$ being 
taken into account, $\iota$ is a morphism
$$
\iota\colon
W'[\fra g]  \longrightarrow 
\roman B\Lambda_{\partial}[s \fra g]
\tag3.3.3
$$
of differential graded coalgebras
which is, furthermore, compatible with the right
$(C\fra g)$-module structures.
\endproclaim

This comparison is formally exactly of the 
same kind as the classical
comparison, spelled out in detail
in \cite\cartanei\ (chap. XIII),
between the CCE complex for an ordinary Lie algebra $\fra g$
and the bar complex for $\roman U[\fra g]$.

The proof will rely on Lemma 3.3.4 below.
To prepare for it, we note that,
as a graded coalgebra,
$$
\bbar \Lambda_{\partial}[s \fra g]=
\bbar \Lambda[s \fra g]=
\rbar \Lambda[s \fra g] \otimes \Lambda'[s \fra g]
=
\roman T'[s\roman I\Lambda[s \fra g]] \otimes 
\Lambda'[s \fra g]
$$ 
and that, as a graded algebra,
$\roman U[C\fra g]= 
\Lambda[s \fra g] \odot\roman U[\fra g]$, the crossed
product algebra.
Abusing the notation $\tau^{\rbar}$ and
$\tau_{\fra g}$,
slightly,
write
$$
\tau^{\rbar}=\tau^{\rbar}\otimes \eta \varepsilon
\colon \bbar \Lambda[s \fra g]=
\roman T's[\roman I\Lambda[s \fra g]] \otimes 
\Lambda'[s \fra g]
\longrightarrow
\Lambda[s \fra g] \odot\roman U[\fra g]=\roman U[C\fra g]
$$
and, likewise,
write
$$
\tau_{\fra g}=
\eta \varepsilon \otimes \tau_{\fra g}
\colon \bbar \Lambda[s \fra g]=
\roman T's[\roman I\Lambda[s \fra g]] \otimes 
\Lambda'[s \fra g]
\longrightarrow
\Lambda[s \fra g] \odot\roman U[\fra g]=\roman U[C\fra g].
$$

\proclaim{Lemma 3.3.4} The sum
$\tau^{\rbar} + \tau_{\fra g}$
is a twisting cochain
$$
\tau^{\rbar} + \tau_{\fra g}
\colon
\bbar \Lambda_{\partial}[s \fra g] \longrightarrow
\roman U[C\fra g].
$$
\endproclaim

\demo{Proof}
We must prove that
$$
D(\tau^{\rbar} + \tau_{\fra g})
 =
(\tau^{\rbar} + \tau_{\fra g})\cup
(\tau^{\rbar} + \tau_{\fra g})
$$
that is,
$$
D\tau^{\rbar} + D\tau_{\fra g} =
\tau^{\rbar}\cup \tau^{\rbar} +
\tau^{\rbar}\cup\tau_{\fra g} +
\tau_{\fra g}\cup \tau^{\rbar}+
\tau_{\fra g}\cup \tau_{\fra g}.
$$
We note first that
$$
\align
 D\tau_{\fra g} &= \tau_{\fra g} (\partial -
\tau^{\rbar} \cap)
\\
 D\tau^{\rbar} &= d\tau^{\rbar}
+\tau^{\rbar}(d+\partial -
\tau^{\rbar} \cap) .
\endalign
$$
By construction,
$\Lambda'_{\partial}[s \fra g]$
is a differential graded subcoalgebra of
$\bbar \Lambda_{\partial}[s \fra g]$,
the algebra $\roman U[\fra g]$ is a differential graded
subalgebra of $\roman U[C\fra g]$,
and the restriction of $\tau^{\rbar} + \tau_{\fra g}$ to 
$\Lambda'_{\partial}[s \fra g]$ amounts to the composite of
$\tau_{\fra g}\colon \Lambda'_{\partial}[s \fra g] \to 
\roman U[\fra g]$ with the injection of
$\roman U[\fra g]$ into 
$\roman U[C\fra g]$.
Consequently
$$ 
\tau_{\fra g}\partial =
\tau_{\fra g}\cup \tau_{\fra g}.
$$
Furthermore,
$\tau^{\rbar}$ is the bar construction twisting cochain
for $\Lambda[s \fra g]$ whence
$$
 \tau^{\rbar}d =  \tau^{\rbar} \cup \tau^{\rbar}
$$
and, likewise,
$$
d \tau^{\rbar}- \tau_{\fra g}(\tau^{\rbar} \cap)=0.
$$
The identities established so far in particular 
show that $\tau$ is a twisting
cochain in the special case where the bracket is zero.

It remains to show that
$$
\tau^{\rbar}\cup\tau_{\fra g} +
\tau_{\fra g}\cup \tau^{\rbar}
=\tau^{\rbar}\partial.
$$
Given the elements $x$ and $y$ of $\fra g$,
$$
\partial(s^2y \otimes sx)=-s\partial(sy \otimes sx)
=s(s[y,x])=-s^2[y,x]
$$
whereas, since $\tau_{\fra g}(sx)=x$,
$$
(\tau^{\rbar}\cup\tau_{\fra g} +
\tau_{\fra g}\cup \tau^{\rbar})(s^2y \otimes sx)=
-s[x,y]
$$
whence, indeed,
$$
\left(\tau^{\rbar}\cup\tau_{\fra g} +
\tau_{\fra g}\cup \tau^{\rbar}\right)
(s^2y \otimes sx)=
-s[x,y]
=\tau^{\rbar}\partial(s^2y \otimes sx) . \qed
$$
\enddemo

\demo{Proof of Theorem {\rm 3.3.2}}
The adjoint
$$
\overline{\tau^{\rbar} + \tau_{\fra g}}
\colon
\bbar \Lambda_{\partial}[s \fra g] 
\longrightarrow
\rbar \roman U[C\fra g]
$$
is a morphism of differential graded coalgebras,
necessarily injective.
The ordinary Weil coalgebra $W'[\fra g]$
has been defined as the CCE coalgebra
$ \Sigm_{\partial}'[sC\fra g]$ for the differential graded
Lie algebra $C\fra g$.
The composite 
$$
(\overline{\tau^{\rbar} + \tau_{\fra g}})
\circ \iota
\colon
W'[\fra g]
\longrightarrow
\rbar \roman U[C\fra g]
$$
plainly coincides with the adjoint
$$
\overline{\tau_{C\fra g}}
\colon W'[\fra g] \longrightarrow 
\rbar \roman U[C\fra g]
$$
of the universal twisting cochain
$\tau_{C\fra g}\colon W'[\fra g] \longrightarrow 
\roman U[C\fra g]$. Consequently $\iota$
is compatible with the structure as asserted. 
\qed \enddemo

\smallskip\noindent
{\smc Remark 3.3.5.\/}
The classical fact that,
relative to the zero bracket,
the canonical injection 
$\SSS[s^2 \fra g] @>>> \rbar\Lambda[s \fra g]$
is an isomorphism on 
homology corresponds to, indeed, is 
equivalent to the statement 
\lq\lq $\roman H(\Lambda^q \Sigma \fra g^*)= \Sigm^q \fra g^*$ 
(in dim $q$)\rq\rq\ in Lemma 3.1 of \cite\bottone.
This statement, in turn, is established there by means of the 
observation that
the appropriate {\it Dold-Puppe\/} 
derived functor \cite\doldpupp\ of the 
$q$-th exterior power functor $\Lambda^q$
is the $q$-th symmetric power functor $\Sigm^q$, 
suitably shifted.

\smallskip\noindent
{\smc 3.3.6. A spectral sequence proof
of Theorem 3.1.9.\/}
Let $N$ be a right $(C\fra g)$-module. 
The comparison $\iota$ plainly induces a morphism
$$
\iota^*\colon\roman{Hom}(W'[\fra g],N)^{C \fra g}
\longrightarrow
\roman{Hom}(\roman B\Lambda_{\partial}[s \fra g],N)^{C \fra g}
$$
of chain complexes.
As a morphism of the underlying graded objects, $\iota^*$ can be written
as
$$
\roman{Hom}(\SSS[s^2\fra g],N)^{\fra g}
\longrightarrow
\roman{Hom}(\rbar\Lambda_{\partial}[s \fra g],N)^{\fra g}.
$$
The coaugmentation filtrations of
$\SSS[s^2\fra g]$ and $\rbar \Lambda_{\partial}[s \fra g]$
induce Serre filtrations on both sides of $\iota^*$
and $\iota^*$ is compatible with the filtrations whence it induces
a morphism between the associated spectral sequences.
At the $E_0$-level, the comparison comes down to the 
standard comparison,
between the complexes induces by the bar
and Koszul resolutions for the exterior algebra $\Lambda[s\fra g]$,
but
restricted to the $\fra g$-invariants,
and thence this comparison has the form
$$
\roman{Hom}(\SSS[s^2\fra g],N)^{\fra g}
\longrightarrow
\roman{Hom}(\rbar\Lambda_{\partial}[s \fra g],N)^{\fra g}.
$$
Since $\fra g$ is reductive,
from the $E_1$-level on, $\iota^{*}$ induces an isomorphism of 
spectral sequences.
Consequently $\iota^*$ is an isomorphism on homology
whence the
relative  
differential graded
$\roman{Ext}_{(C \fra g, \fra g)}(R,N)$
is the homology of the chain complex (3.1.10) as asserted.
The same kind of reasoning shows that,
for any left $(C \fra g)$-module $M$,
the relative  differential graded 
$\roman{Tor}^{(C \fra g, \fra g)}(R,M)$
is the homology of the chain complex (3.1.11).
This proves Theorem 3.1.9.

\smallskip
\noindent{\smc 3.4. The Weil coalgebra $W'[\fra g]$ as a relative
$\roman U[\fra g]$-contractible construction in the
 reductive case.\/}
Suppose that the ground ring $R$ is a field of characteristic zero
and let $\fra g$ be a {\it reductive\/} Lie algebra.
The diagonal map of $\fra g$ induces a graded commutative 
algebra structure on $\roman H^*(\fra g)$ and,
furthermore, a graded cocommutative coalgebra
structure on $\roman H_*(\fra g)$.
The projection
$\pi\colon \Lambda'_{\partial}[s\fra g] \to\Lambda'_{\partial}[s\fra g]\otimes_{\fra g}R$,
restricted to the invariants
$\Lambda'[s\fra g]^{\fra g}$,
is an isomorphism whence
the differential on $\Lambda'_{\partial}[s\fra g]\otimes_{\fra g}R$
is zero, and we will write $\Lambda'[s\fra g]\otimes_{\fra g}R$
rather than $\Lambda'_{\partial}[s\fra g]\otimes_{\fra g}R$. 
This quotient is naturally isomorphic to the homology 
$\roman H_*(\fra g)$; further,
as a chain complex, $\Lambda'_{\partial}[s\fra g]$ 
decomposes as
$$
\Lambda'_{\partial}[s\fra g]= 
\roman{ker}(\pi) \oplus \Lambda'_{\partial}[s\fra g]^{\fra g}
\cong \roman{ker}(\pi) \oplus
\Lambda'[s\fra g]\otimes_{\fra g}R,
$$
and $\roman{ker} (\pi)$ is a contractible chain complex.
The quotient $\Lambda'[s\fra g]\otimes_{\fra g}R\cong \roman H_*(\fra g)$ 
acquires
a graded cocommutative coalgebra structure
for purely formal reasons; the resulting coalgebra structure on 
$\roman H_*(\fra g)$ is the one induced by the diagonal of $\fra g$.

The kernel of the composite 
of the projection $\roman U[C\hgg]\to \Lambda'_{\partial}[s\fra g]$ with $\pi$
is the (differential graded)
two-sided ideal $\langle\hgg \rangle$ in $\roman U[C\hgg]$
generated by $\fra g$ 
whence  $\roman H_*(\hgg)$ acquires a graded commutative algebra
structure which combines with the coalgebra structure 
to a Hopf algebra structure. Dually the cohomology
$\roman H^*(\fra g)$ acquires a coalgebra structure
which combines with its  algebra structure
to a Hopf algebra structure.
As an algebra, $\roman H^*(\fra g)$  
is the exterior algebra $\Lambda[\roman{Prim}(\fra g)]$ 
generated by the {\it primitives\/}
$\roman{Prim}(\fra g)\subseteq \roman H^*(\fra g)$ relative to the 
coalgebra structure.
The dual $I(\fra g)$ of $\roman{Prim}(\fra g)$  
is the module of {\it indecomposables\/}
relative to the algebra structure on $\roman H_*(\fra g)$,
the injection of $\roman{Prim}(\fra g)$ into $\roman H^*(\fra g)$
dualizes to the canonical projection 
$\roman H_*(\fra g) \to I(\fra g)$ (defining the indecomposables)
and, as a graded coalgebra, $\roman H_*(\fra g)$ 
is the exterior coalgebra $\Lambda'[I(\fra g)]$ cogenerated by $I(\fra g)$.

As a graded coalgebra,
$\roman{Tor}^{(C \fra g, \fra g)}(R,R)$
is the cofree graded cocomutative
coalgebra $\SSS[sI(\fra g)]$ cogenerated by the suspension
$sI(\fra g)$ of $I(\fra g)$. 
Pick a section $j \colon I(\fra g) \to \Lambda[I(\fra g)]$
for the projection $\Lambda[I(\fra g)] \to I(\fra g)$.
Then the composite
$$
\tau\colon
\SSS[s I(\fra g)] @>{\roman{pr}}>> s I(\fra g) @>{s^{-1}}>> 
I(\fra g) @>{j}>>
 \Lambda[I(\fra g)]
$$
is a transgression twisting cochain.
By {\it transgression twisting cochain\/} we mean a twisting cochain
which induces the transgression in the corresponding spectral
sequence; cf. \cite\maclaboo\ for the notion of transgression in a 
spectral sequence.
We will 
write $\Lambda=\Lambda[I(\fra g)]\ (=\roman H_*(\fra g))$ 
and $\SSS = \SSS[sI(\fra g)]$.

\proclaim{Theorem 3.4.1} The Lie algebra $\fra g$ being reductive,
the Weil coalgebra
$W'[\fra g]$
admits a $\fra g$-equivariant contracting homotopy.
\endproclaim

\demo{Proof} Since $\fra g$ is reductive, Hodge theory yields a contraction
$$
\Nsddata {\SSS[s^2\fra \hgg] \otimes_{\tau_{\fra \hgg}}
\Lambda_{\partial}'[s \fra \hgg]}
{\nabla_1}{\pi_1}{\SSS \otimes \Lambda'}{h_1}
\tag3.4.1.1
$$
in the category of $\fra g$-modules, the $\fra g$-actions on 
$\SSS$  and $\Lambda'$ being trivial.
See e.~g. \cite\kostaeig\ for details.
Recall that $W'[\fra g]=( \SSS[s^2\fra \hgg] \otimes_{\tau_{\fra \hgg}}
\Lambda_{\partial}'[s \fra \hgg],d)$.
Relative to the Serre filtrations, 
the  contraction (3.4.1.1) is a filtered contraction, and
an application of the perturbation lemma
transforms the  contraction (3.4.1.1) into the (filtered) contraction
$$
\Nsddata {W'[\fra g]}
{\nabla_2}{\pi_2}{\SSS \otimes_{\tau} \Lambda'}{h_2}
\tag3.4.1.2
$$
in the category of $\fra g$-modules. 
The twisted tensor product $\SSS \otimes_{\tau} \Lambda'$ is contractible;
thus let 
$$
\Nsddata {\SSS \otimes_{\tau} \Lambda'}
{\eta}{\varepsilon}{R}{h_3}
\tag3.4.1.3
$$
be a contraction of
$\SSS \otimes_{\tau} \Lambda'$
onto the ground ring  $R$,
necessarily $\fra g$-equivariant, the $\fra g$-actions being trivial;
beware: the contracting homotopy $h_3$ is {\it not\/} unique.
Further, this kind of contraction is {\it not\/} a filtered one
relative to the Serre filtration of
$\SSS \otimes_{\tau} \Lambda'$, though.
Combining the two contractions, we obtain the contraction
$$
\Nsddata {W'[\fra g]}
{\eta}{\varepsilon}{R}{h}
\tag3.4.1.4
$$
in the category of $\fra g$-modules. \qed
\enddemo

Theorem 3.1.9 is a consequence of Theorem 3.4.1.
Indeed, the Lie algebra being reductive, pick
a $\fra g$-equivariant contracting homotopy
of $W'[\fra g]$ of the kind constructed in Theorem 3.4.1. Then
the canonical comparison,
cf. \cite\maclaboo\ (Theorem IX.6.2 on p.~267), \cite\mooretwo,
yields a 
$(\roman U[C\fra g])$-linear 
morphism 
$\alpha_{\partial}\colon \roman B \Lambda_{\partial}[s \fra g] \to W'[\fra g]$
and homogeneous  $(\roman U[C\fra g])$-linear operators
$h_{\partial}$ on  $\roman B \Lambda_{\partial}[s \fra g]$
and $h_W$ on $W'[\fra g]$ of degree 1
such that
$$
Dh_{\partial} = \roman{Id}- \iota \alpha_{\partial},
\quad
Dh_W = \roman{Id}- \alpha_{\partial} \iota .
\tag3.4.2
$$
Thus the data
$$
\Nsddata {\roman B \Lambda_{\partial}[s \fra g]}
{\iota_{\phantom {\partial}}}{\alpha_{\partial}}{h_W, W'[\fra g]}{h_{\partial}}
\tag3.4.3
$$
constitute a filtered chain equivalence 
which is, furthermore,  $(\roman U[C\fra g])$-linear.
The notion of filtered chain equivalence was introduced
in \cite\huebkade;
in the present paper we shall exclusively use the defining property
(3.4.2), though, and no reference to \cite\huebkade\ will be made.
The standard reasoning then immediately establishes Theorem 3.1.9.

\smallskip
\noindent
{\smc 3.5. The Cartan model\/}. Return momentarily to a general ground
ring $R$. Let $\fra g$ be an $R$-Lie algebra which, 
as an $R$-module, is projective. Given the right $(C\fra g)$-module $N$,
the chain complex
$\roman{Hom}(W'[\fra g],N)^{C \fra g}$, cf. (3.1.10) above, is still defined.
We will now rewrite this chain complex as a twisted object.

When the Lie bracket on $\fra g$ is ignored,
as a differential graded right $(\Lambda[s \fra g])$-module, 
$W'[\fra g]$
has the form 
$
\SSS[s^2 \fra g] \otimes_{\tau^{\SSS}} \Lambda[s \fra g],
$
the differential being the operator
$$
\partial ^{\tau^{\SSS}}=-( \tau^{\SSS}\cap \,\cdot\, )
\colon \SSS[s^2\fra \hgg]
 \otimes \Lambda[s \fra \hgg]
\longrightarrow
\SSS[s^2\fra \hgg] \otimes \Lambda[s \fra \hgg]
$$
relative to the universal twisting cochain
 $
\tau^{\SSS}\colon
\SSS[s^2\fra \hgg] \longrightarrow 
\Lambda[s \fra \hgg].
$
With reference to the graded right
$(\Lambda[s \fra g])$-module structure on $N$,
when the differential on $N$ is ignored,
the twisted Hom-object
$\roman{Hom}^{\tau^{\SSS}}(\SSS[s^2\fra g], N)$
is defined; 
we remind the reader that the operator 
$\delta^{\tau^{\SSS}}$ on
$\roman{Hom}(\SSS[s^2\fra g], N)$, cf. (1.2) above and
\cite\duaone\ (2.4.1),
is defined by
$\delta^{\tau^{\SSS}}(f) = (-1)^{|f|}f \cup \tau^{\SSS}$, the argument $f$
being a homogeneous morphism.
With the differential on $N$
and the Lie bracket on $\fra g$ incorporated,
on the $\fra g$-invariants, the operator 
$\delta^{\tau^{\SSS}}$ on
$\roman{Hom}(\SSS[s^2\fra g], N)$
is still a
perturbation of the differential on $\roman{Hom}(\SSS[s^2\fra g], N)$
(coming from that on $N$), and we write the resulting
twisted object as
$$
\roman{Hom}^{\tau^{\SSS}}(\SSS[s^2\fra g], N)^{\fra g}.
\tag3.5.1
$$

The exterior algebra $\Lambda=\Lambda[s \fra \hgg]$ 
is a Hopf algebra.
Recall that, in terms of the notation
$\Lambda= \eta\varepsilon + \iota$, the antipode $S$ of $\Lambda$ can be written as
$$
S=\eta\varepsilon - \iota+
\iota\cup\iota - \iota^{\cup 3} +
\ldots
= \sum (-1)^j\iota^{\cup j} 
\colon \Lambda \to \Lambda.
$$
At the risk of notational confusion
with our notation for a symmetric algebra, 
we use here the standard notation
$S$ for the antipode. This notation for the antipode is not 
used elsewhere in the paper.

Let $\Cal A$ be a general graded Hopf algebra and let $M$ and $N$ be
ordinary right $\Cal A$-modules.
Then $\Cal A$ acts on the right of
$\roman{Hom}(M,N)$ in various ways:

The $\Cal A$-actions on $M$ and $N$ induce the $\Cal A$-actions 
$$
\align\mu_M&\colon \roman{Hom}(M,N) \otimes 
\Cal A \longrightarrow \roman{Hom}(M,N)
\\
\mu_N&\colon \roman{Hom}(M,N) \otimes 
\Cal A \longrightarrow \roman{Hom}(M,N)
\endalign
$$
on the right of
$\roman{Hom}(M,N)$, and
the two pairings
$\mu_M$ and $\mu_N$ combine to an $\Cal A$-action
$$
\mu_{M,N}\colon \roman{Hom}(M,N) \otimes 
\Cal A \longrightarrow \roman{Hom}(M,N)
$$
on  $\roman{Hom}(M,N)$
given as the composite of the following two morphisms:
$$
\align
\roman{Hom}(M,N) \otimes 
\Cal A 
&@>{\roman{Hom}(M,N) \otimes\Delta}>>
\roman{Hom}(M,N) \otimes \Cal A \otimes
\Cal A 
\\
\roman{Hom}(M,N) \otimes \Cal A \otimes
\Cal A 
&@>{\mu_M \otimes\Cal A}>>
\roman{Hom}(M,N) \otimes \Cal A 
@>{\mu_N}>>
\roman{Hom}(M,N).
\endalign
$$

We will now take $M$ to be $\Cal A$ itself, viewed as a right
$\Cal A$-module via right multiplication. The following is well known and classical.

\proclaim{Lemma 3.5.2} The association $\alpha \longmapsto \alpha \cup \Cal A$
induces an isomorphism
$$
\psi\colon\left(\roman{Hom}(\Cal A,N),\mu_{\Cal A}\right) 
\longrightarrow
\left(\roman{Hom}(\Cal A,N),\mu_{\Cal A,N}\right) 
$$
of right $\Cal A$-modules.
The inverse isomorphism
$$
\phi\colon\left(\roman{Hom}(\Cal A,N),\mu_{\Cal A,N}\right) 
\longrightarrow
\left(\roman{Hom}(\Cal A,N),\mu_{\Cal A}\right) 
$$
is given 
by the association $\beta \longmapsto \beta \cup S$. \qed
\endproclaim

\proclaim{Theorem 3.5.3} {\rm [Cartan]}
The assignment to a homogeneous
$\alpha \in \roman{Hom}(\SSS[s^2 \fra g],N)^{\fra g}$ of
$$
\Phi_{\alpha} \colon
\SSS[s^2 \fra g] \otimes\Lambda[s \fra g]
\longrightarrow N,
\quad
\Phi_{\alpha} (w\otimes a)= \alpha(w) a,
\ w \in \SSS[s^2 \fra g],\ a \in \Lambda[s \fra g],
\tag3.5.4
$$
yields an injective morphism
$$
\roman{Hom}(\SSS[s^2\fra g], N)^{\fra g}
@>>>
\roman{Hom}
\left(\SSS[s^2\fra g],\roman{Hom}(\Lambda'[s\fra g],N)\right)
\tag3.5.5
$$
of bigraded $R$-modules, and the composite of 
{\rm (3.5.5)} with the  morphism
$$
\phi_*\colon\roman{Hom}
\left(\SSS[s^2\fra g],\roman{Hom}(\Lambda'[s\fra g],N)\right)
@>>>
\left(\SSS[s^2\fra g],\roman{Hom}(\Lambda'[s\fra g],N)\right)
\tag3.5.6
$$
of bigraded $R$-modules induced by
$
\phi\colon\roman{Hom}(\Lambda'[s\fra g],N)
@>>>
\roman{Hom}(\Lambda'[s\fra g],N)
$,
combined with the adjointness isomorphism
$$
\left(\SSS[s^2\fra g],\roman{Hom}(\Lambda'[s\fra g],N)\right)
\cong \roman{Hom}(W'[\fra g],N),
$$
yields
 an injective chain map
$$
\roman{Hom}^{\tau^{\SSS}}(\SSS[s^2\fra g], N)^{\fra g}
@>>>
\roman{Hom}(W'[\fra g],N)
\tag3.5.7
$$
which identifies the source {\rm (3.5.1)}
of {\rm (3.5.7)} with 
$\roman{Hom}(W'[\fra g],N)^{C\fra g}$
\endproclaim

\demo{Proof} In the special case where
the ground ring is that of the reals and where $\fra g$ is compact,
this goes back to {\it Cartan\/} \cite\cartanon,
and the reasoning in the general case is formally the same. \qed \enddemo

We will now suppose that $R$ is a field of characteristic zero and that
$\fra g$ is reductive.
We then refer to the twisted object (3.5.1) as the {\it Cartan\/} model for
the differential graded
$\roman{Ext}_{(C \fra g, \fra g)}(R,N)$
and to a morphism of the kind (3.5.6)
as a {\it Cartan twist\/}.
The following is immediate.

\proclaim{Corollary 3.5.8}
The chain map {\rm (3.5.7)} identifies
the source {\rm (3.5.1)}
of {\rm (3.5.7)} with {\rm (3.1.10)}. 
Thus the Cartan model indeed calculates the 
differential graded
$\roman{Ext}_{(C \fra g, \fra g)}(R,N)$.
Consequently
$\roman{Ext}_{(C \fra g, \fra g)}(R,R)$ is 
canonically isomorphic to 
the algebra of
$\fra g$-invariants of 
the algebra $\roman{Hom}(\SSS[s^2\fra \hgg],R)$,
that is, to
the algebra of
$\fra g$-invariants of 
the symmetric algebra $\Sigm[s^{-2} \fra g^*]$.
\qed
\endproclaim

\proclaim{Corollary 3.5.9}
The differential graded
$\roman{Ext}_{(C \fra g, \fra g)}(R,N)$ acquires
the structure
$$
\roman{Hom}(\SSS[s^2\fra \hgg],R)^{\fra g}
\otimes 
\roman{Ext}_{(C \fra g, \fra g)}(R,N)
\longrightarrow
\roman{Ext}_{(C \fra g, \fra g)}(R,N)
\tag3.5.10
$$
of a $\roman{Hom}(\SSS[s^2\fra \hgg],R)^{\fra g}$-module
via the induced
$(\SSS[s^2\fra g])$-comodule structure on $W'[\fra g]$. \qed
\endproclaim

\noindent{\smc 3.6. 
Cutting the Cartan model to size in the reductive case.\/}
Suppose that $R$ is a field of characteristic zero and
let $\hgg$ be a reductive Lie algebra.
Let $\VVV$ be a differential graded right
$(C\hgg)$-module. The space $\VVV^{\hgg}$ of invariants is
manifestly a $(C\hgg)$-submodule and the induced $(C\hgg)$-action
on $\VVV^{\hgg}$,
restricted to $\hgg$, is plainly trivial. Consequently the
$(\roman U[C\hgg])$-action on $\VVV^{\hgg}$,
restricted to the two-sided
differential graded ideal $\langle\hgg \rangle$
in $\roman U[C\hgg]$ generated by 
$\hgg$, is trivial whence the action passes through
an action of the quotient algebra
$\roman U[C\hgg]\big /\langle\hgg \rangle$ on $\VVV^{\hgg}$,
and this action is compatible with the differentials.
Hence the action of $\roman U[C\hgg]$ on $\VVV$
then passes to an action
$$
\roman H_*(\hgg ) \otimes \VVV^\hgg @>>> \VVV^\hgg
\tag3.6.1
$$
of $\roman H_*(\hgg)$ on $\VVV^\hgg$ which is compatible with the 
differential
on $\VVV^\hgg$.

Through the projection $\roman U[C \fra g] \to \roman H_*(\fra g)=\Lambda$ of 
differential graded algebras, the twisted object
$\SSS\otimes_{\tau}\Lambda$ acquires a canonical right 
$(C\fra g)$-module structure.
Exploiting the Hopf algebra structure of $\roman U[C\fra g]$,
we endow $\SSS\otimes_{\tau}\Lambda\otimes_{\tau} W'[\fra g]$ with the diagonal
right $(C\fra g)$-module structure and, likewise
exploiting the Hopf algebra structure of $\SSS$, we endow
$\SSS\otimes_{\tau}\Lambda\otimes_{\tau} W'[\fra g]$ with the diagonal
left $\SSS$-comodule structure. To spell out the latter,
we write the multiplication map of $\SSS$ as
$\mu \colon \SSS\otimes \SSS \to \SSS$, the comodule structure
maps
as $\Delta\colon \SSS\otimes\Lambda \to 
\SSS \otimes \SSS\otimes \Lambda$ and
$\Delta\colon W'[\fra g] \to 
\SSS \otimes W'[\fra g]$, and the twist map as
$T \colon \Lambda \otimes \SSS \to \SSS \otimes \Lambda$.
With these preparations out of the way,
the
structure map of the diagonal left $\SSS$-comodule structure
is the composite of the following three morphisms
$$
\align
\SSS\otimes\Lambda\otimes W'[\fra g]
&@>{\Delta \otimes \Delta}>>
\SSS\otimes\SSS\otimes\Lambda\otimes \SSS\otimes W'[\fra g]
\\
\SSS\otimes\SSS\otimes\Lambda\otimes \SSS\otimes W'[\fra g]
&@>{\SSS\otimes\SSS\otimes T \otimes W'[\fra g]}>>
\SSS\otimes\SSS\otimes \SSS\otimes\Lambda\otimes W'[\fra g]
\\
\SSS\otimes\SSS\otimes \SSS\otimes\Lambda\otimes W'[\fra g]
&@>{\mu\otimes \SSS\otimes\Lambda\otimes W'[\fra g]}>>
\SSS\otimes \SSS\otimes\Lambda\otimes W'[\fra g] .
\endalign
$$
This  left $\SSS$-comodule structure is compatible with 
the differentials.

The $(C\fra g)$-linear projections
$$
\align
\varepsilon \otimes \varepsilon \otimes W'[\fra g]
&\colon\SSS\otimes_{\tau}\Lambda\otimes_{\tau} W'[\fra g]
\longrightarrow W'[\fra g]
\tag3.6.2
\\
\SSS\otimes\Lambda\otimes\varepsilon
&\colon\SSS\otimes_{\tau}\Lambda\otimes_{\tau}W'[\fra g]
\longrightarrow\SSS\otimes_{\tau}\Lambda 
\tag3.6.3
\endalign
$$
are morphisms of right $(C\fra g)$-modules and left
$\SSS$-comodules; furthermore, (3.6.2) and (3.6.3)
are chain equivalences since the objects
involved are contractible.

\proclaim{Proposition 3.6.4}
The projections {\rm (3.6.2)} and {\rm (3.6.3)}
induce chain equivalences
$$
\align
\roman{Hom}(W'[\fra g], \VV)^{C\fra g}
&@>>>\roman{Hom}
(\SSS\otimes_{\tau}\Lambda\otimes_{\tau} W'[\fra g],
\VV)^{C\fra g}
\tag3.6.5
\\
\roman{Hom}(\SSS\otimes_{\tau}\Lambda, \VV^{\fra g})^{C\fra g}
&@>>>\roman{Hom}
(\SSS\otimes_{\tau}\Lambda\otimes_{\tau}W'[\fra g],
\VV)^{C\fra g}
\tag3.6.6
\endalign
$$
that are compatible with the induced differential graded
$\roman{Hom}(\SSS,R)$-module structures.
As a graded algebra,
$\roman{Hom}(\SSS, R)$
is therefore canonically isomorphic to 
$\roman{Ext}_{(C \fra g, \fra g)}(R,R)$ and, for general  $\VV$,
the twisted
Hom-object
$\roman{Hom}^{\tau}(\SSS, \VV^{\fra g})$
is a small model for the differential graded
$\roman{Ext}_{(C \fra g, \fra g)}(R,\VV)$
that is compatible with the bundle structures
in the sense that
the resulting pairing
$$
\roman{Hom}(\SSS, R)
\otimes \roman{Hom}^{\tau}(\SSS, \VV^{\fra g})
\longrightarrow
\roman{Hom}^{\tau}(\SSS, \VV^{\fra g})
\tag3.6.7
$$
induces the $\roman{Hom}(\SSS, R)$-module structure
{\rm (3.5.10)} on $\roman{Ext}_{(C \fra g, \fra g)}(R,\VV)$.
\endproclaim

\demo{Proof}
A spectral sequence comparison argument shows that
(3.6.5) and (3.6.6) are isomorphisms on homology and hence
chain equivalences.

Adjointness yields the isomorphism
$$
\roman{Hom}^{\tau}(\SSS, \VV^{\fra g})
@>>>
\roman{Hom}(\SSS\otimes_{\tau}\Lambda, \VV^{\fra g})^{\Lambda}
=\roman{Hom}(\SSS\otimes_{\tau}\Lambda, \VV)^{C\fra g}
$$
which combines with (3.6.6) to the chain equivalence
$$
\roman{Hom}^{\tau}(\SSS, \VV^{\fra g})
@>>>\roman{Hom}
(\SSS\otimes_{\tau}\Lambda\otimes_{\tau}W'[\fra g],\VV)^{C\fra g}.
\tag3.6.8
$$
Since (3.6.5) is a chain equivalence as well,
the left-hand side of (3.6.8) is a small model for the differential graded
$\roman{Ext}_{(C \fra g, \fra g)}(R,\VV)$ as asserted. 
The verification of the compatibility with the 
bundle structures is left to the reader. \qed
\enddemo

We will refer to
$\roman{Hom}^{\tau}(\SSS, \VV^{\fra g})$ as the {\it small Cartan model\/}
 for the differential graded
$\roman{Ext}_{(C \fra g, \fra g)}(R,\VV)$.
The following is an immediate consequence of Proposition 3.6.4.

\proclaim{Theorem 3.6.9}
The canonical map
$$
\roman{Ext}_{\Lambda}(R,\VV^{\fra g})\longrightarrow
\roman{Ext}_{(C \fra g, \fra g)}(R,\VV)
$$
induced by the projection $\roman U[C \fra g] \to \Lambda=\roman H_*(\fra g)$
is an isomorphism. \qed
\endproclaim

\smallskip\noindent
{\smc Remark 3.6.10.\/} The construction of an explicit section
for (3.6.2) is not entirely obvious whence
the construction of an explicit map between
$\roman{Hom}^{\tau}(\SSS, \VV^{\fra g})$
and the Weil model requires some care.
Incomplete reasoning and faulty constructions
aiming at comparing the Weil model
(or Cartan model) with the small Cartan model
led to a certain activity in the literature
\cite\allpuone, \cite\almei, \cite\gorkomac, \cite\maszwebe,
see in particular the introduction of
\cite\almei.

\medskip\noindent {\bf 4. The simplicial Weil coalgebra\/}
\smallskip\noindent
As before, let $R$ be a commutativ ring and
$\fra g$  an ordinary
$R$-Lie algebra  which is projective as an $R$-module.
In the previous section,
we explored the {\it nonhomogeneous\/} form of the relative bar resolution
for the pair $(\roman U[C\fra g], \roman U[\fra g])$.
The present aim is to show that
totalization carries a suitably defined simplicial Weil coalgebra
associated with $\fra g$
to the {\it homogeneous\/} form of the relative bar resolution under discussion.
This homogeneous bar resolution will enable us to 
introduce small models for the corresponding relative differential
Ext-functors.

A  cosimplicial version of the
Weil algebra associated with a Lie algebra has been introduced 
in \cite\kamtonfo\ and \cite\kamtontw\ (p.~59), see
also \cite\kamtonfi. 
The dual of our simplicial Weil coalgebra does {\it not\/} coincide
with the cosimplicial version of the
Weil algebra  explored in
\cite\kamtonfo--\cite\kamtonfi, though.
Yet we prefer to stick to our terminology  since Corollary 4.7
below will establish a canonical comparison
between the ordinary Weil coalgebra and the total object associated
with the simplicial Weil coalgebra in our sense.

As before, let $\fra g$ be an $R$-Lie algebra which
is projective as an $R$-module. 
We remind the reader that the notation $DK$ refers to the Dold-Kan functor,
cf. (2.5) above.
The simplicial $R$-module
$DK\fra g$ acquires a simplicial Lie algebra structure.
We will often discard the symbol $DK$
in notation and thus
view $\fra g$ as a simplicial Lie algebra
whenever necessary,
each structure map being the identity.
The total object $E\fra g$ associated 
with $\fra g$ 
in the category of Lie algebras 
relative to the obvious monoidal structure
is a simplicial Lie algebra as well.
We refer to the resulting
simplicial differential graded CCE coalgebra
$\Lambda'_{\partial}[s E\fra g]$
as the 
{\it simplicial Weil coalgebra\/}
associated with $\fra g$.

In the same vein, consider the symmetric monoidal category of 
coaugmented 
differential graded cocommutative coalgebras,
with the differential graded coalgebra tensor product $\otimes$
as monoidal structure,
and with the coalgebra diagonal as diagonal
structure for the category---it is here where the
requirement that the coalgebras be graded cocommutative
is needed. 
In this category,
the total object
$E^{\otimes}\Lambda'_{\partial}[s \fra g]$
is a simplicial
differential graded coalgebra.
By functoriality, 
the differential
graded $(C \fra g)$-action 
on $\Lambda'_{\partial}[s \fra g]$
given in (1.3) above
induces  a {\it differential
graded $(C \fra g)$-action 
on $E^{\otimes}\Lambda'_{\partial}[s \fra g]$
that is compatible with the coalgebra structure\/}.

\proclaim{Lemma 4.1}
As a simplicial differential graded coalgebra,
the CCE coalgebra $\Lambda'_{\partial}[s E\fra g]$, carried out
for the simplicial Lie algebra $E\fra g$, that is, for the total
object associated with $\fra g$,
is canonically isomorphic to the total object
$E^{\otimes}\Lambda'_{\partial}[s \fra g]$
associated with the CCE coalgebra $\Lambda'_{\partial}[s \fra g]$
for $\fra g$
in the category of differential graded coalgebras.
\endproclaim

\demo{Proof}
In a simplicial degree $p \geq 0$,
$$
(E\fra g)_p = \fra g \times \dots \times\fra g \quad (p+1\ \text{factors})
$$
whence, as $p$ ranges over the natural numbers,
$$
(\Lambda'_{\partial}[s E\fra g])_p =
(\Lambda'_{\partial}[(s \fra g)^{\oplus (p+1)}])
\cong
(\Lambda'_{\partial}[s \fra g])^{\otimes (p+1)}
\cong
(E^{\otimes}\Lambda'_{\partial}[s \fra g])_p.
$$
These isomorphisms are compatible with the simplicial operations.\qed
\enddemo

We will denote by
 $|\Lambda'_{\partial}[s E\fra g]|$
the complex
$$
@>d>>
|\Lambda'_{\partial}[s E\fra g]|_n
@>d>>
\ldots @>d>> |\Lambda'_{\partial}[s E\fra g]|_1
@>d>>
\Lambda'_{\partial}[s \fra g]
\tag4.2
$$
of right $(C\fra g)$-modules,
necessarily contractible,
arising from 
the simplicial differential graded coalgebra
$\Lambda'_{\partial}[s E\fra g]$
by normalization.
Our next aim is to
identify
the complex (4.2)
with the relative bar resolution (3.2.2).
To this end, we recall the (right)
{\it nonhomogeneous\/} version of
$E\fra g$. The situation is formally the same as that in (2.5) above:
View momentarily $E\fra g$ merely as a 
graded $R$-module and 
consider the familiar automorphism
$$
\Phi \colon E\fra g \longrightarrow
E\fra g
$$ 
of graded $R$-modules which, 
in degree $n$,
that is, on $(E\fra g)_n=\fra g^{\times(n+1)}$,
is given by the formula
$$
\Phi(x_0,x_1,\dots, x_n)=
(x_0+x_1+x_2+\ldots + x_n,
\ldots,x_{n-1}+x_n, x_n).
$$
In degree $n$, the inverse of $\Phi$ is 
plainly given by the association 
$$
(y_1,\ldots,y_n,x) \longmapsto 
(x_0,x_1,\dots, x_n)=
(y_1-y_2,y_2-y_3,\ldots,y_{n-1}-y_n, y_n-x,x).
$$
The {\it nonhomogeneous\/} face operators 
$\partial_j$ are given by
the familiar expressions
$$
\aligned\partial_0(x_0,x_1,\dots, x_n) &= (x_1,\dots, x_n),
\\
\partial_j(x_0,x_1,\dots, x_n) &= 
(x_0,\dots, x_{j-2},x_{j-1}+x_j,
x_{j+1},\dots, x_n)\ 
(1 \leq j \leq n)
\endaligned
\tag4.3
$$
and, likewise,
the {\it nonhomogeneous\/} degeneracy operators 
$s_j$ are given by
$$
s_j(x_0,x_1,\dots, x_n) = 
(x_0,\dots, x_{j-1},0,x_j,\dots, x_n)
\ (0 \leq j \leq n);
\tag4.4
$$
we will denote the resulting
simplicial $R$-module by 
$(E\fra g)^{\roman{right}}$.
The automorphism $\Phi$ of graded $R$-modules
is well known to be
 an isomorphism
$$
\Phi\colon (E\fra g)^{\roman{right}}
\longrightarrow E\fra g
$$
of simplicial $R$-modules.
Here the chosen
nonhomogeneous
formulas reflect the 
$(C\fra g)$-operations being 
from the {\it right\/}, and our 
constructions are written as
tensor product of \lq\lq {\it base\/}\rq\rq\  
times \lq\lq {\it fiber\/}\rq\rq.
The formulas in the literature for objects with
operators from the left differ from the 
above ones.
See e.~g. p. 75 of \cite\gugenmay.

\proclaim{Theorem 4.5} The isomorphism
$\Phi$ of simplicial $R$-modules
induces an isomorphism
$$
\beta_{\roman U[\fra g]} (R,\roman U[C \fra g],\roman U[C \fra g])
\longrightarrow
\Lambda'_{\partial}[s E \fra g]
\tag4.5.1
$$
of simplicial right $(\roman U[C\fra g])$-complexes.
Consequently the complex 
$|\Lambda'_{\partial}[s E \fra g]|$
of right $(\roman U[C\fra g])$-modules
comes down to the homogeneous form of the normalized relative 
bar resolution of $R$ in the category
of right $(\roman U[C\fra g])$-modules.
\endproclaim

\demo{Proof}
In terms of the nonhomogeneous description, in a fixed
degree $n\geq 0$,
the Lie algebra 
$(E\fra g)_n^{\roman{right}}$ simply 
comes down to the iterated
semi-direct product
$$
\fra g \rtimes \fra g \rtimes \ldots \rtimes \fra g\ 
(n+1\ \text{copies of}
\  \fra g) ,
$$
the formulas {\rm (4.3)} and {\rm (4.4)} 
for the 
simplicial structure
still being
valid.
A special case thereof
is the observation that the linear map
$$
\fra g \rtimes \fra g \longrightarrow \fra g \times \fra g,
\ (x_0,x_1)\longmapsto (x_0+x_1,x_1)
$$
is an isomorphism of Lie algebras with inverse mapping given by the
assignment to $(y,x)$ of $(y-x,x)$.
Consequently, for a fixed $n\geq 0$,
the degree $n$ differential graded coalgebra 
$\Lambda'_{\partial}[s (E\fra g)^{\roman{right}}]_n$
of the
simplicial differential graded coalgebra
$\Lambda'_{\partial}[s (E\fra g)^{\roman{right}}]$
takes the form of an iterated twisted 
tensor product
$$
\Lambda'_{\partial}[s\fra g]\otimes_{\tau_{\fra g}}
\Lambda'_{\partial}[s\fra g]\otimes_{\tau_{\fra g}}
\ldots \otimes_{\tau_{\fra g}}
\Lambda'_{\partial}[s\fra g] 
$$
of $n+1$ copies of the 
CCE coalgebra $\Lambda'_{\partial}[s\fra g]$
of $\fra g$ relative to the appropriate canonical actions
of the corresponding copy of $\fra g$ on the right of
that part of
$\Lambda'_{\partial}[s (E\fra g)^{\roman{right}}]_n$
left to $\fra g$ in the tensor product
decomposition. 
For fixed $n$,
the differential graded coalgebra
$\Lambda'_{\partial}[s (E\fra g)^{\roman{right}}]_n$
is precisely the corresponding constituent
(3.2.2) of the relative simplicial bar construction.
In particular, 
each face and degeneracy operator is manifestly a
morphism of differential graded coalgebras,
and these operators are exactly the same as those
in the simplicial bar construction. Hence
$\Phi$ induces an isomorphism from the complex
(3.2.2) onto (4.2). This observation establishes
Theorem 4.5. \qed \enddemo

\noindent
{\smc Remark 4.6.\/}
Since the condensed object 
$\roman B\Lambda_{\partial}[s \fra g]=
\left|R \otimes_{\roman U[C \fra g]} 
\beta_{\roman U[\fra g]} (\roman U[C \fra g])
\right|^{\bullet}
$
(cf. {\rm (3.2.5)})
in Theorem 3.2.6 arises from the 
simplicial differential graded coalgebra
$\Lambda'_{\partial}[s (E\fra g)]$ by normalization
and condensation, the condensed object $\roman B\Lambda_{\partial}[s \fra g]$
acquires a
differential graded coalgebra structure.
This establishes Theorem 3.2.6

\proclaim{Corollary 4.7} The comparison {\rm (3.3.3)},
combined with the induced morphism $|4.5.1|^{\bullet}$
between the condensed objects,
yields a morphism
$$
W'[\fra g] \longrightarrow \left|\Lambda'_{\partial}[s (E\fra g)]\right|^{\bullet}
$$
of differential graded coalgebras
between the ordinary Weil coalgebra and the total object associated with the simplicial
Weil coalgebra for $\fra g$. \qed
\endproclaim

\medskip\noindent
{\bf 5. Cutting the defining object for 
$\roman{Ext}_{((G,C\fra g);\chain)}$ to size}
\smallskip\noindent
As before, let $R$ be a commutative ring and $\fra g$ an $R$-Lie algebra
which is projective as an $R$-module.

\smallskip
\noindent {\smc 5.1. The model for 
$\roman{Ext}_{((G,C\fra g);\chain)}$ arising from the Weil coalgebra.\/}
The realization $|E\fra g|$ of the total simplicial
Lie algebra $E\fra g$  is a differential  graded Lie algebra, and the
simplicial twisting cochain
$$
\tau_{E\fra g}\colon \Lambda'_{\partial}[s E\fra g]
@>>>
\roman U[E\fra g]
$$
the constituents of which in each 
simplicial degree are given
in (1.2) above induces, via the 
twisted Eilenberg-Zilber theorem
\cite\gugenhtw,
an acyclic twisting cochain
$$
\tau_{|E\fra g|}\colon \LLL
=|\Lambda'_{\partial}[s E\fra g]|^{\bullet} 
@>>> \roman U[|E\fra g|].
\tag5.1.1
$$
More precisely, as chain complexes, $|s E\fra g| \cong s |E\fra
g|$, and the canonical projection from the differential graded
coalgebra $|\Lambda'_{\partial}[s E\fra g]|^{\bullet}$ to $|s E\fra g|$
determines, via the universal property of the differential graded
CCE coalgebra $\SSS_{\partial}[s |E\fra g|]$ for the
differential graded Lie algebra $|E\fra g|$, a unique morphism
$$
\LLL=|\Lambda'_{\partial}[s E\fra g]|^{\bullet} 
@>>> \SSS_{\partial}[s |E\fra
g|]
$$
of differential graded coalgebras which, combined with the
universal twisting cochain for the CCE construction of $|E\fra
g|$, yields the asserted twisting cochain $\tau_{|E\fra g|}$.

We now take the ground ring to be that of the reals, $\Bbb R$.
Until the end of the present section, 
we take $G$ to be a Lie group, $\fra g$ its Lie algebra,
an $\VV$ a right $(G,C \fra g)$-module.
In view of Proposition 2.5.7, the chain complex $|\Cal A^0(EG,\VV)|$ 
is the standard injective resolution of $\VV$
in the category of (differentiable) right $G$-modules
and, in particular,
carries a canonical right $G$-module structure. Furthermore, the
obvious componentwise actions of the constituents of the
simplicial Lie algebra $E\fra g$ on the constituents of $\Cal
A^0(EG,\VV)$ induce an action of the differential graded Lie algebra
$|E\fra g|$ on $|\Cal A^0(EG,\VV)|$; this action does {\it not\/}
involve $\VV$. Let 
$$
\Cal B^*_{(G,C \fra g)}(\Bobb R,G,\VV)
=\roman{Hom}^{\tau_{|E\fra g|}}(\LLL,
|\Cal A^0(EG,\VV)|),
\tag5.1.2
$$
the resulting twisted Hom-object. Since the $|E\fra g|$-action does not involve $\VV$,
the twisting cochain $\tau_{|E\fra g|}$ does not involve $\VV$.
The twisted Hom-object (5.1.2) inherits a canonical $G$-action.
Furthermore, the assignment to 
$G$ and $\VV$ of 
$\Cal B^*_{(G,C \fra g)}(\Bobb R,G,\VV)$
is plainly a 
functor in the group variable and, given the group $G$,
in the $(G,C \fra g)$-module
variable as well.
The functor $\Cal B^*_{(G,C \fra g)}(\,\cdot\, , \,\cdot\, , \,\cdot\, )$
is, in a somewhat generalized sense, a 
{\it dualized unreduced bar construction\/} for the category
of $(G,C \fra g)$-modules.
Indeed, the object (5.1.2) acquires a natural $(G,C \fra g)$-module structure.
For  $\LLL$ and $\VV$
inherit graded (not differential graded)
$(\Lambda[s \fra g])$-module structures
from their differential graded $(C\fra g)$-module structures, and these
$(\Lambda[s \fra g])$-module structures
induce a graded (not differential graded)
$(\Lambda[s \fra g])$-module structure on
the corresponding untwisted object
$$
\roman{Hom}(\LLL,|\Cal A^0(EG,\VV)|).
\tag5.1.3
$$
On the twisted object (5.1.2), the induced differential graded 
$\fra g$-module
structure (the diagonal structure coming from that on
$\LLL$ and the diagonal structure on
$|\Cal A^0(EG,\VV)|$) and the
graded $(\Lambda[s \fra g])$-module structure
combine to a differential graded right $(C \fra g)$-module structure.
Thus $\Cal B^*_{(G,C \fra g)}(\Bobb R,G,\VV)$ 
is a $(G,C \fra g)$-module functor,
that is, a functor having as range the category of $(G,C \fra g)$-modules.

For intelligibility we recall that, given the group $H$ (we will
then substitute $EG$ for $H$)
and the subgroup $G$, the twisted object
$$
\roman{Hom}^{\tau_{\fra h}}(\Lambda'_{\partial}[s\fra h],\Cal
A^0(H,\VV))
$$
is defined relative to the universal Lie algebra twisting cochain
$\tau_{\fra h}\colon \Lambda'_{\partial}[s\fra h] \to \roman U [\fra h]$.
Here $\Cal A^0(H,\VV)$ is a left $\fra h$-module
via right translation in $H$; this structure does not involve
$\VV$, and the operator 
$$
\delta^{\tau_{\fra h}}\colon \roman{Hom}(\Lambda'[s\fra h],\Cal
A^0(H,\VV))\longrightarrow
\roman{Hom}(\Lambda'[s\fra h],\Cal
A^0(H,\VV))
$$
determined by the universal Lie algebra twisting cochain
$\tau_{\fra h}\colon \Lambda'_{\partial}[s\fra h] \to \roman U[s\fra h]$,
cf. (1.2) above, is defined; 
cf. (1.6.3.1) above where this is explained for the special case
where $\VV$ is the de Rham complex $\Cal A(X)$ of a $G$-manifold $X$.
Thus the operator
$$
\delta^{\tau_{|E\fra g|}}\colon 
\roman{Hom}(\LLL,|\Cal A^0(EG,\VV)|)
\longrightarrow
\roman{Hom}(\LLL,|\Cal A^0(EG,\VV)|)
$$
determined by the universal Lie algebra twisting cochain
$\tau_{|E\fra g|}\colon \LLL \to \roman U[|E\fra g|]$
does not involve $\VV$.

Define the 
functors
$\overline {\Cal B}{^*_{(\,\cdot\, ,\,\cdot\,)}}(\Bobb R,\,\cdot \, ,\,\cdot\,)$,
$ \Cal B^*_{(\,\cdot\,,\,\cdot\,)}(\,\cdot\,)$
and 
$\overline{\Cal B}^*_{(\,\cdot\,,\,\cdot\,)}(\,\cdot\,)$
by
$$
\aligned
\overline {\Cal B}{^*_{(G,C \fra g)}}(\Bobb R,G,\VV)
&=\Cal B^*_{(G,C \fra g)}(\Bobb R,G,\VV)^{(G,C \fra g)}
\\
\Cal B^*_{(G,C \fra g)}(G)
&=\Cal B^*_{(G,C \fra g)}(\Bobb R,G,\Bobb R)
\\
\overline {\Cal B}{^*_{(G,C \fra g)}}(G)
&= \overline {\Cal B}{^*_{(G,C \fra g)}}(\Bobb R,G,\Bobb R).
\endaligned
\tag5.1.4
$$
We have chosen the notation  $\Cal B^*$ and
$\overline{\Cal B}^*$ since
$\Cal B^*$ and $\overline{\Cal B}^*$ are,
in a somewhat generalized sense, dualized respective 
unreduced and reduced  bar constructions, 
with reference to the category $\roman{Mod}_{(G,C \fra g)}$ and,
accordingly, we refer to
$\Cal B^*$ and
$\overline{\Cal B}^*$
as {\it unreduced\/} and {\it reduced\/}
constructions, respectively.
In particular,
$\Cal B^*_{(G,C \fra g)}(G)$
and
$\overline {\Cal B}{^*_{(G,C \fra g)}}(G)$
are differential graded algebras; further,
$\Cal B^*_{(G,C \fra g)}(\Bobb R,G,\VV)$
is a $\Cal B^*_{(G,C \fra g)}(G)$-module
and
$\overline {\Cal B}{^*_{(G,C \fra g)}}(\Bobb R,G,\VV)$
is a
$\overline{\Cal B}^*_{(G,C \fra g)}(G)$-module
in an obvious manner.

\proclaim{Theorem 5.1.5} 
The 
 differential graded $\roman{Ext}_{(G,C \fra g)}(\Bbb R, \VV)$
is canonically isomorphic to the homology
of
the twisted object
$\overline {\Cal B}^*_{(G,C \fra g)}(\Bobb R,G,\VV)$.
\endproclaim

Theorem 5.1.5 is an immediate consequence of the following lemma.

\proclaim{Lemma 5.1.6} The left trivialization of 
the tangent bundle
of $G$ induces a contraction of the totalized complex $\left|\Cal
A(EG,\VV)^{(G,C\fra g)}\right|$ onto 
the twisted object
$\overline {\Cal B}^*_{(G,C \fra g)}(\Bobb R,G,\VV)$.
\endproclaim

\demo{Proof} The de Rham theory
 Eilenberg-Zilber
theorem 
yields the contraction
$$
\Nsddata {|\roman{Hom}^{\tau_{E\fra g}} (\Lambda'_{\partial}[s
E\fra g],\Cal A^0(EG,\VV))|} {\alpha^\flat}{\nabla^\flat} {\Cal
B^*_{(G,C \fra g)}(\Bobb R,G,\VV)} {h^\flat} ,
$$
necessarily $G$- and $(C \fra g)$-equivariant.
Exploiting  Theorem 2.6.1 (the extended decomposition lemma),
we replace 
$ {|\roman{Hom}^{\tau_{E\fra g}} (\Lambda'_{\partial}[s
E\fra g],\Cal A^0(EG,\VV))|}$ with
$|\Cal A(EG,\VV)|$ and, thereafter, we take
$(G,C\fra g)$-invariants.
This yields the contraction
$$
\Nsddata {|\Cal A(EG,\VV)^{(G,C\fra g)}|} 
{\alpha^*}{\nabla^*} {\overline {\Cal
B}^*_{(G,C \fra g)}(\Bobb R,G,\VV)} {h^*} . \qed
\tag5.1.7
$$
\enddemo

Thus the homology of the twisted object
$
{\overline {\Cal B}^*_{(G,C \fra g)}(\Bobb R,G,\VV)}
$
coincides with the differential graded 
$\roman{Ext}_{((G,C \fra g);\chain)}(\Bobb R,\VV)$ whence Theorem 5.1.5.

Recall that $\iota\colon W'[\fra g]  \to
\roman B\Lambda_{\partial}[s \fra g]
$ refers to the canonical comparison (3.3.3), cf. also Corollary 4.7.
The composite of $\iota$ with the acyclic twisting cochain (5.1.1) is plainly
an acyclic twisting cochain
$$
\tau_{|E\fra g|}\circ \iota\colon W'[\fra g]
@>>> \roman U[|E\fra g|].
$$

\proclaim{Proposition 5.1.8} {\rm (i)} The comparison $\iota$
induces a homology isomorphism
$$
{\overline {\Cal B}^*_{(G,C \fra g)}(\Bobb R,G,\VV)}
\longrightarrow
{\roman{Hom}^{\tau_{|E\fra g|}\circ\iota}
(W'[\fra g],|\Cal A^0(EG,\VV)|)^{(G,C \fra g)} } 
\tag5.1.9
$$
between {\rm (5.1.2)} and
$$
\roman{Hom}^{\tau_{|E\fra g|}\circ\iota}(W'[\fra g],
|\Cal A^0(EG,\VV)|)^{(G,C\fra g)}.
\tag5.1.10
$$
{\rm(ii)} When $G$ is reductive,
the chain equivalence {\rm (3.4.3)} induces
a $(G,C \fra g)$-equivariant chain equivalence of the kind
$$
\Nsddata
{\Cal B^*_{(G,C \fra g)}(\Bobb R,G,\VV)}
{\alpha^\flat}{\nabla^\flat}
{h_W^\flat,\roman{Hom}^{\tau_{|E\fra g|}\circ\iota}(W'[\fra g],
|\Cal A^0(EG,\VV)|)}
{h^\flat}
\tag5.1.11
$$
and, taking $(G,C \fra g)$-invariants,
we obtain the chain equivalence
$$
\Nsddata
{\overline {\Cal B}^*_{(G,C \fra g)}(\Bobb R,G,\VV)}
{\alpha^\sharp}{\nabla^\sharp}
{h_W^\sharp,\roman{Hom}^{\tau_{|E\fra g|}\circ\iota}
(W'[\fra g],|\Cal A^0(EG,\VV)|)^{(G,C \fra g)}}
{h^\sharp}. 
\tag5.1.12
$$
\endproclaim

Taking  $(G,C \fra g)$-invariants means taking $(\Lambda[s \fra g])$-
and $G$-invariants;
the $(\Lambda[s \fra g])$-invariants are the horizontal elements
in a sense explained earlier.

\demo{Proof}
Filtering the twisted objects (5.1.2) and (5.1.10)
by the degree complementary to the $G$-resolution degree
we obtain spectral sequences
$(E_r(5.1.2),d_r)$ and $(E_r(5.1.10),d_r)$ 
converging to $\roman{Ext}_{((G,C \fra g);\chain)}(\Bobb R,\VV)$
and the total object of (5.1.10) respectively,
and the comparison $\iota$ induces a morphism
$$
(E_r(5.1.2),d_r) \longrightarrow (E_r(5.1.10),d_r)
$$
of spectral sequences.
The bigraded $R$-module
$$
\roman{Hom}(W'[\fra g],
|\Cal A^0(EG,\VV)|)^{(G,C\fra g)},
$$
endowed with the bar complex operator alone,
amounts to the chain complex
$$
|\Cal A^0(EG,\roman{Hom}(W'[\fra g],\VV)^{\Lambda[s\fra g]})|^G
\cong 
|\Cal A^0(EG,\roman{Hom}(\SSS[s^2\fra g],\VV))|^G
$$
whence
$$
E_1(5.1.10) \cong \roman H^*_{\roman{cont}}\left(G,
\roman{Hom}(W'[\fra g],\VV)^{\Lambda[s\fra g]}\right)
= \roman H^*_{\roman{cont}}\left(G,
\roman{Hom}(\SSS[s^2\fra g],\VV)\right),
$$
the $(\Lambda[s\fra g])$-action being given by contraction,
and the operator $d_1$ is the Koszul resolution operator,
 induced by the operator
$\partial^{\tau^{\SSS}}$ on 
$\roman{Hom}(\SSS[s^2\fra g],\VV)$.
Consequently
$$
E_2(5.1.10) \cong \roman H^*_{\roman{cont}}\left(G,
\roman{Ext}_{\Lambda[s\fra g]}(\Bobb R,\VV)\right).
$$

Likewise, the bigraded $R$-module
$$
\roman{Hom}(\roman B\Lambda_{\partial}[s \fra g],
|\Cal A^0(EG,\VV)|)^{(G,C\fra g)},
$$
endowed with the bar complex operator alone,
amounts to the chain complex
$$
|\Cal A^0(EG,\roman{Hom}(\roman B\Lambda[s \fra g],\VV)^{\Lambda[s\fra g]})|^G
\cong 
|\Cal A^0(EG,\roman{Hom}(\overline {\roman B}[\Lambda[s\fra g],\VV))|^G
$$
whence
$$
E_1(5.1.2) \cong \roman H^*_{\roman{cont}}\left(G,
\roman{Hom}(\roman B\Lambda[s \fra g],\VV)^{\Lambda[s\fra g]}\right)
= \roman H^*_{\roman{cont}}\left(G,
\roman{Hom}(\overline {\roman B}[\Lambda[s\fra g],\VV)\right),
$$
the $(\Lambda[s\fra g])$-action being given by contraction,
and the operator $d_1$ is the bar resolution operator.
Consequently
$$
E_2(5.1.2) \cong  \roman H^*_{\roman{cont}}\left(G,
\roman{Ext}_{\Lambda[s\fra g]}(\Bobb R,\VV)\right).
$$
Hence the ordinary spectral sequence comparison
establishes the assertion (i).
Assertion (ii) is immediate. \qed 
\enddemo

\proclaim{Theorem 5.1.13}
Via the comparison $\iota$,
the  differential graded $\roman{Ext}_{(G,C \fra g)}(\Bbb R, \VV)$
is canonically isomorphic to the homology of the twisted object {\rm (5.1.10)},
viz. of
$$
{\roman{Hom}^{\tau_{|E\fra g|}\circ\iota}
(W'[\fra g],|\Cal A^0(EG,\VV)|)^{(G,C \fra g)}}. 
$$
\endproclaim

\demo{Proof} This is an immediate consequence of Proposition 5.1.8. \qed
\enddemo

The twisted object (5.1.10) has somewhat the form of a {\it Weil model\/},
with $|\Cal A^0(EG,\VV)|$ instead of the module $\VV$ itself in the ordinary Weil model.
There is also a corresponding object which takes the form
of a {\it Cartan model\/}:
The composite of
the injective morphism
$$
\roman{Hom}(\SSS[s^2\fra g], |\Cal A^0(EG,
\VV)|)^{G}
@>>>
\roman{Hom}
\left(\SSS[s^2\fra g],\roman{Hom}(\Lambda'[s\fra g],|\Cal A^0(EG,
\VV)|)\right)
$$
of bigraded $R$-modules with the {\it Cartan twist\/}
$$
\phi_*\colon\roman{Hom}
\left(\SSS[s^2\fra g],\roman{Hom}(\Lambda'[s\fra g],|\Cal A^0(EG,
\VV)|)\right)
@>>>
\left(\SSS[s^2\fra g],\roman{Hom}(\Lambda'[s\fra g],|\Cal A^0(EG,
\VV)|)\right),
$$
cf. (3.5.6) above,
combined with the adjointness isomorphism
$$
\left(\SSS[s^2\fra g],\roman{Hom}(\Lambda'[s\fra g],|\Cal A^0(EG,
\VV)|)\right)
\cong \roman{Hom}(W'[\fra g],|\Cal A^0(EG,
\VV)|),
$$
yields an injective  morphism of graded vector spaces
$$
\roman{Hom}(\SSS[s^2\fra g], |\Cal A^0(EG,\VV)|)^{G}
@>>>
\roman{Hom}(W'[\fra g],|\Cal A^0(EG,
\VV)|)
$$
which induces an isomorphism from
$$
\roman{Hom}^{\tau^{\SSS},\tau_{|E\fra g|}\circ\iota}
(\SSS[s^2\fra g],|\Cal A^0(EG,
\VV)|)^G
\tag5.1.14
$$
onto the twisted object (5.1.10). 
The reasoning is essentially the same as that which establishes Theorem 3.5.3.
The total differential of (5.1.14) has the form $d + \partial$
of a perturbed differential:
The operator $d$ is the naive differential
on $\roman{Hom}(\SSS[s^2\fra g],|\Cal A^0(EG,\VV)|)^G$
coming from the bar complex operator $\delta$ and the differential on $\VV$.
Furthermore,
$$
\partial = \delta^{\tau^{\SSS}} +
\delta^{\tau_{|E\fra g|}\circ\iota},
\tag5.1.15
$$
where
$\delta^{\tau^{\SSS}}$
is the operator defined
with reference to the action of $\Lambda[s\fra g]$
on $\Cal A^0(EG,\VV)|$  coming from the
action of $\Lambda[s\fra g]$ on $\VV$,
and where
$\delta^{\tau_{|E\fra g|}\circ\iota}$
is the  operator
induced from the twisting cochain
$\tau_{|E\fra g|}\circ \iota$; see (1.2) above 
for details.
In view of Theorem 5.1.5, the twisted object (5.1.14)
calculates the   differential graded 
$\roman{Ext}_{((G,C \fra g);\chain)}(\Bbb R, \VV)$;
the twisted  object (5.1.14) is somewhat smaller than
the original object
defining $\roman{Ext}_{((G,C \fra g);\chain)}(\Bbb R, \VV)$.
By adjointness, we may rewrite the graded object
which underlies (5.1.14) as
$$
|\Cal A^0(EG,\roman{Hom}(\SSS[s^2\fra g],\VV))|^G.
\tag5.1.16
$$
Since $EG$ is contractible,
the chain complex $|\Cal A^0(EG,\roman{Hom}(\SSS[s^2\fra g],\VV))|$,
endowed with the operator $\delta$,
is a differentiably injective resolution of
$\roman{Hom}(\SSS[s^2\fra g] ,\VV)$,
and taking $G$-invariants we obtain precisely
the object (5.1.16) endowed merely with the bar complex operator
$\delta$, which therefore
calculates the differentiable cohomology of $G$ with coefficients in
$\roman{Hom}(\SSS[s^2\fra g],\VV)$.
For the special case 
where $\fra g$ is finite dimensional and $\VV$ the real
numbers with trivial action,
this is exactly Theorem 1 in \cite\bottone.

\smallskip\noindent
{\smc Remark 5.1.17.\/} (Relationship with the Bott spectral sequence)
Suppose that $G$ is a finite-dimensional Lie group and let $X$ be a left $G$-manifold.
Substituting $\Cal A(X)$ for $\VV$ in the spectral sequence
$(E_r(5.1.2),d_r)$,
we obtain a spectral sequence
$\left(E_r(G,X),d_r\right)$ having
$$
E_2 =  \roman H^*_{\roman{cont}}\left(G,
\roman{Ext}_{\Lambda[s\fra g]}(\Bobb R,\Cal A(X))\right).
$$
For $X$ a point, this is the  spectral sequence
explored by {\it Bott\/} in \cite\bottone.
In particular,
$\roman{Ext}_{\Lambda[s\fra g]}(\Bobb R,\Bobb R) = \roman{Hom}(\SSS[s^2\fra g],\Bobb R)$
and, in Theorem 1 in \cite\bottone, 
the object which corresponds
to our 
$\roman{Hom}(\SSS[s^2\fra g],\Bobb R)$ 
is written as $\Sigm \fra g^*$. 
Likewise, the spectral sequence $(E_r(5.1.10),d_r)$
has the form of a van Est spectral sequence. Indeed,
for a finite-dimensional connected Lie group $G$ and a $G$-representation
$V$, the complex $\Cal A(G,V)$ of $V$-valued forms
on $G$ can be written as $\roman{Alt}(\fra g,\Cal A^0(G,V))$
and the $G$-invariant subcomplex $\Cal A(G,V)^G$ 
amounts to the CCE complex $\roman{Alt}(\fra g,V)$
calculating $\roman H^*(\fra g,V)$ where
$V$ is
viewed as a $\fra g$-module.
However, since $G/K$ is contractible,
the cohomology of
$
\Cal A(G/K,\Cal A(G,V))^G
$
amounts to the cohomology of $\Cal A(G,V)^G$.
The spectral sequence of the 
form degree filtration 
of 
$
\Cal A(G/K,\Cal A(G,V))^G
$
relative to $G/K$ is the van Est spectral sequence \cite\vanestwo.
This spectral sequence has
$$
E_2 =  \roman H_{\roman{cont}}^*(G,\roman H_{\roman{top}}^*(G,V))
$$
and converges to $\roman H^*(\fra g,V)$.
See also Theorem 2.10 in \cite\kamtontw.

\smallskip
\noindent
{\smc 5.2. A small object 
for $\roman{Ext}_{((G,C \fra g);\chain)}$ in the strictly 
exterior case.\/} At the present stage,
the Lie group $G$ is not supposed to be reductive.
To simplify the exposition somewhat,
define the functor
$$
t^* \colon 
\roman{Mod}_{(G,C \fra g)} @>>> 
{}_{\overline {\Cal B}{^*_{(G,C \fra g)}}(G)}\roman{Mod}
\tag5.2.1
$$
by the assignment to a (right) $(G,C \fra g)$-module $\Nflat$
of the twisted object 
$$
t^*(\Nflat) = \overline {\Cal B}{^*_{(G,C \fra g)}}
(\Bobb R,G,\Nflat)
$$
which, in turn, {\it calculates\/} the differential graded 
$\roman{Ext}_{((G,C \fra g);\chain)}(\Bobb R,\Nflat)$.
This functor is one of two (Koszul) duality
functors; we shall come back to this duality in Section 7 below.

\proclaim{Theorem 5.2.2} Suppose that $G$ is of strictly
exterior type in such a way that $\roman H^*(G)$ is the exterior Hopf algebra
$\Lambda[y_1,\dots]$ in suitable universally transgressive
generators $y_1,y_2, \dots $ and,
for each $y_j$ in $\roman H^*G$, choose
a  cycle in $\overline {\Cal B}^*$  such that  $y_j$  transgresses to 
the class of this cycle. This choice determines a differential
$\Cal D$ such that
$$
(\roman{Hom}(\roman H_*(BG), \VV),\Cal D)
\tag5.2.3
$$
is a small model calculating $\roman{Ext}_{((G,C \fra g);\chain)}(\Bobb R,\VV)$.
\endproclaim

At this stage, the group $G$ is a {\it general\/} Lie group
of strictly exterior type, possibly
{\it infinite dimensional\/}.  The 
hypothesis of Theorem 5.2.2 
is of course automatically satisfied when
$G$ is finite dimensional and connected. In this particular case,
an explicit description of the differential $\Cal D$ will 
be given later.

We now begin with the preparations for the proof of Theorem 5.2.2.
For
simplicity, we will momentarily write ${\Cal B}^*_{(G,C \fra g)}(G)$
as ${\Cal B}^*$ and $\overline {\Cal B}^*_{((G,C \fra g),\fra
g)}(G)$ as $\overline {\Cal B}^*$.
We define the completed tensor
product ${\Cal B}^* \widehat \otimes {\Cal B}^*$ by
$$
{\Cal B}^* \widehat \otimes {\Cal B}^*
=\roman{Hom}^{\tau^{\otimes}}(\roman B\Lambda_{\partial}[s\fra g]
\otimes \roman B \Lambda_{\partial}[s\fra g], \Cal A^0(E(G \times
G)))
$$
where
$$
\tau^{\otimes}=\tau_{|E\fra g|} \otimes \eta \varepsilon + \eta
\varepsilon \otimes \tau_{|E\fra g|} \colon \roman
B\Lambda_{\partial} \otimes \roman B\Lambda_{\partial} \to \roman
U[|E\fra g|] \otimes \roman U[|E\fra g|] \cong \roman U[|E(\fra g
\times \fra g)|].
$$
Likewise, we define the reduced completed tensor product
$\overline {\Cal B}^* \widehat \otimes \overline {\Cal B}^*$ by
$$
\overline {\Cal B}^* \widehat \otimes \overline {\Cal B}^* =
({\Cal B}^* \widehat \otimes {\Cal B}^*)^{(G \times G, C \fra g
\times C\fra g)} .
$$
We will now exploit the graded coalgebra structure which underlies 
the graded Hopf algebra
$\roman H^*(G)$.

\proclaim{Lemma 5.2.4} For each $y_j$ in $\roman H^*G$, choose
a  cycle in $\overline {\Cal B}^*$  such that  $y_j$  transgresses to 
the class of this cycle.
This choice determines an acyclic twisting cochain
$$
\zeta_{\overline {\Cal B}^*}\colon \roman H^*G
@>>>
\overline {\Cal B}^*.
\tag5.2.5
$$
\endproclaim

\demo{Proof} Given two twisting cochains $\sigma \colon C \to
\overline {\Cal B}^*$ and $\sigma' \colon C' \to \overline {\Cal
B}^*$ defined on graded cocommutative coalgebras $C$ and $C'$, the
values of the twisting cochain $\sigma \otimes \eta \varepsilon +
\eta \varepsilon \otimes \sigma'$ lie in $ \overline {\Cal B}^*
\widehat \otimes \overline {\Cal B}^* $.
Moreover, the Eilenberg-Zilber theorem furnishes a contraction of
the kind
$$
\Nsddata
{\overline
{\Cal B}^*_{(G\times G,C (\fra g \times \fra g))}
(G\times G)}
{\alpha^*}{\nabla^*}
{
\overline {\Cal B}^*_{(G,C \fra g)}(G)
\widehat \otimes
\overline {\Cal B}^*_{(G,C \fra g)}(G)
}{h^*}.
\tag5.2.6
$$
The inductive construction of the twisting cochain (3.2.1${}_*$) in 
the proof
of Lemma 3.2${}_*$ in \cite\duaone\ dualizes, with
\cite\duaone\ (2.2.1$^*$) instead of \cite\duaone\ (2.2.1$_*$).
We leave the details to the reader. \qed
\enddemo

\demo{Proof of Theorem {\rm 5.2.2}} Suppose that a choice of cycle in
$\overline {\Cal B}^*$ for each $y_j$ has been made so that the
resulting twisting cochain $\zeta_{\overline {\Cal B}^*}$ is available. We
first construct a contraction from $(\roman H^*G)
\otimes_{\zeta_{\overline {\Cal B}^*}}t^*(\VV)$ onto $\VV$. 
To this end we note first that, in view of the naturality of 
the functor $\overline {\Cal B}$ in the group variable,
restriction of the construction to the trivial group
induces a canonical
projection $t^*(\VV) \to \VV$ 
which thus forgets the $(G,C \fra g)$-structure.  This projection
extends to a projection $\alpha$ from $(\roman H^*G)
\otimes_{\zeta_{\overline {\Cal B}^*}}t^*(\VV)$ to $\VV$, necessarily a chain
equivalence. A section for this projection which is compatible
with the differentials includes an sh-comodule structure over
$\roman H^*G$ on $\VV$. To construct such a section, we consider $\VV$
momentarily endowed with the trivial $(G,C \fra g)$-structure
which we refer to by the notation $\VV(0)$. Then the obvious
injection of $\VV$ into $t^*(\VV)$ which, for $v \in \VV$, assigns
$\psi_v \colon \rbar \Lambda_{\partial} \to \Cal A^0(EG,\VV)$ to $v
\in \VV$ given by
$$
\psi_v(w)(x) = \varepsilon(w) xv, \quad w \in
\rbar \Lambda_{\partial}, \ x \in EG,
$$
induces an injective chain map
$j(0) \colon \VV(0) \to (\roman H^*G) \otimes_{\zeta_{\overline {\Cal B}^*}}t^*(\VV(0))$,
and there is an obvious extension
$$
\Nsddata
{(\roman H^*G) \otimes_{\zeta_{\overline {\Cal B}^*}}t^*(\VV(0))}
{j(0)}
{\alpha}
{\VV(0)}
{h(0)}
$$
of the data to a contraction. Incorporating the non-trivial
$(G,C \fra g)$-structure on $\VV$ amounts to perturbing the
differential on the right-hand side via an operator $\partial$
which lowers the filtration
coming from the coaugmentation filtration
of $\rbar \Lambda_{\partial}[s \fra g]$ and
the simplicial degree filtration with reference to $EG$.
Application of the perturbation lemma, cf. \cite\duaone\ (2.4),
yields the  contraction
$$
\Nsddata
{(\roman H^*G) \otimes_{\zeta_{\overline {\Cal B}^*}}t^*(\VV)}
{j}
{\alpha}
{\VV}
{h}.
\tag5.2.7
$$
In this manner, the extended  $(\roman H^*G)$-comodule structure on
$(\roman H^*G) \otimes_{\zeta_{\overline {\Cal B}^*}}t^*(\VV)$
recovers an sh-comodule structure on $\VV$ over $\roman H^*G$.

To complete the construction of a small model,
we extend the contraction (5.2.7) to a contraction
$$
\Nsddata
{(\roman{Hom}(\roman H_*(BG),
(\roman H^*G) \otimes_{\zeta_{\overline {\Cal B}^*}}t^*(\VV))}
{j}
{\alpha}
{\roman{Hom}(\roman H_*(BG), \VV)}
{h}
\tag5.2.8
$$
in the obvious manner
where the notation $j$, $\alpha$, $h$ is abused somewhat.
Let $\tau \colon \roman H_*(BG) \to \roman H_*G$
be the transgression twisting cochain.
Application of the perturbation lemma, 
cf. \cite\duaone\ (2.3),
yields the contraction
$$
\Nsddata
{\roman{Hom}^{\tau}(\roman H_*(BG),
(\roman H^*G) 
\otimes_{\zeta_{\overline {\Cal B}^*}}t^*(\VV))}
{j_{\tau}}
{\alpha_{\tau}}
{(\roman{Hom}(\roman H_*(BG), \VV),\Cal D)}
{h_{\tau}}
\tag5.2.9
$$
The left-hand side of (5.2.9)
yields the desired small model 
calculating \linebreak
$\roman{Ext}_{((G,C \fra g);\chain)}(\Bobb R,\VV)$. \qed
\enddemo

\noindent
{\smc 5.3. The Weil and Cartan models for compact $G$\/.} 
Recall  the obvious acyclic twisting cochain
 $\tau^{\SSS}\colon\SSS[s^2 \fra g] \to 
\Lambda[s\fra g]$.
With reference to the graded right
$(\Lambda[s \fra g])$-module structure on $\VV$,
when the differential on $\VV$ is ignored,
the twisted Hom-object
$\roman{Hom}^{\tau^{\SSS}}(\SSS[s^2\fra g], \VV)$
is defined; 
we remind the reader that the operator 
$\delta^{\tau^{\SSS}}$ on
$\roman{Hom}(\SSS[s^2\fra g], \VV)$, cf. (1.2) above and
\cite\duaone\ (2.4.1),
is defined by
$\delta^{\tau^{\SSS}}(f) = (-1)^{|f|}f \cup \tau^{\SSS}$, the argument $f$
being a homogeneous morphism.
 Further we will 
write the differential on $\VV$ and the differential it induces on
$\roman{Hom}^{\tau^{\SSS}}(\SSS[s^2\fra g], \VV)$
as $d$, with an abuse of notation.
A slight extension of the reasoning for Theorem
3.5.3 establishes the following.

\proclaim{Proposition 5.3.1}
For a general Lie group $G$, 
the canonical morphism of graded objects from 
$\roman{Hom}(\SSS[s^2\fra g], \VV)$ to
$\roman{Hom}(W'[\fra g], \VV)$,
followed by the Cartan twist {\rm (3.5.6)} and the appropriate
adjointness isomorphism,
induces an isomorphism 
$$
\roman{Hom}^{\tau^{\SSS}}(\SSS[s^2\fra g], \VV)^G
\longrightarrow
\roman{Hom}(W'[\fra g], \VV)^{(G,C\fra g)}
$$
of graded $R$-modules
such that the differential on the right-hand side 
passes to the sum $\delta^{\tau^{\SSS}}+d$, restricted to the 
$G$-invariants. \qed
\endproclaim

Until the end of the present subsection we suppose $G$ to be
finite-dimensional and compact.

\proclaim{Theorem 5.3.2}
The chain complex 
$\roman{Hom}(W'[\fra g], \VV)^{(G,C\fra g)}$
and the twisted object
$\roman{Hom}^{\tau^{\SSS}}(\SSS[s^2\fra g], \VV)^G$
are small models  
for $\roman{Ext}_{((G,C \fra g);\chain)}(\Bobb R,\VV)$,
and these models are compatible with the bundle structures
in the sense that the obvious pairings
$$
\roman{Hom}(W'[\fra g], \Bbb R)^{(G,C\fra g)}
\otimes
\roman{Hom}(W'[\fra g], \VV)^{(G,C\fra g)}
\longrightarrow
\roman{Hom}(W'[\fra g], \VV)^{(G,C\fra g)}
$$
and
$$
\roman{Hom}(\SSS[s^2\fra g],\Bbb R)^G
\otimes
\roman{Hom}^{\tau^{\SSS}}(\SSS[s^2\fra g], \VV)^G
\longrightarrow
\roman{Hom}^{\tau^{\SSS}}(\SSS[s^2\fra g], \VV)^G
$$
induce the $(\roman H^*(BG))$-module structure on
$\roman{Ext}_{((G,C \fra g);\chain)}(\Bobb R,\VV)$.
\endproclaim

\demo{Proof}
For a general Lie group $G$,
the canonical injection $\overline j$ of right $G$-modules extends to a contraction
$$
\Nsddata {|\Cal A^0(EG,\VV)|}
{\overline j}{\overline \alpha}{\VV}{\overline h}
\tag5.3.3
$$
of chain complexes; this contraction  
encapsulates the fact that $|\Cal A^0(EG,\VV)|$ is a
differentiably injective $G$-resolution of $\VV$. The contraction (5.3.3), in turn,
induces a contraction  
$$
\Nsddata
{\roman{Hom}^{\tau_{|E\fra g|}\circ \iota}(W'[\fra g],
|\Cal A^0(EG,\VV)|)}
{j}
{\widetilde \alpha}
{\roman{Hom}(W'[\fra g], \VV)}
{\widetilde h}
\tag5.3.4
$$
of chain complexes where $j$ is a morphism of $(G,C\fra g)$-modules.

The group $G$ being compact, integration over $G$
transforms the contraction 
(5.3.4)
into the $(G,C \fra g)$-equivariant contraction
$$
\Nsddata
{\roman{Hom}^{\tau_{|E\fra g|}\circ\iota}(W'[\fra g],
|\Cal A^0(EG,\VV)|)}
{j}
{\alpha}
{\roman{Hom}(W'[\fra g], \VV)}
{h}.
\tag5.3.5
$$
Taking $(G,C \fra g)$-invariants on both sides,
we obtain the contraction
$$
\Nsddata
{\roman{Hom}^{\tau_{|E\fra g|}\circ\iota}
(W'[\fra g],|\Cal A^0(EG,
\VV)|)^{(G,C\fra g)}}
{j}
{\alpha}
{
\roman{Hom}(W'[\fra g], \VV)^{(G,C\fra g)}}
{h}
\tag5.3.6
$$
where the notation $\alpha$ and $h$ is abused somewhat;
notice that $j$ remains unchanged under integration, though,
since it was already $G$-equivariant.
Proposition 5.3.1 implies that
$
\roman{Hom}^{\tau^{\SSS}}(\SSS[s^2\fra g], \VV)^G
$
is a
{\it small model\/}
for 
$\roman{Ext}_{((G,C \fra g);\chain)}(\Bobb R,\VV)$ as well.  \qed 
\enddemo

We will refer to
$\roman{Hom}(W'[\fra g], \VV)^{(G,C\fra g)}$ as the {\it Weil model\/}
for $\roman{Ext}_{((G,C \fra g);\chain)}(\Bobb R,\VV)$
and to the twisted object 
$\roman{Hom}^{\tau^{\SSS}}(\SSS[s^2\fra g], \VV)^G$
as the {\it Cartan model\/}
for $\roman{Ext}_{((G,C \fra g);\chain)}(\Bobb R,\VV)$
associated with $G$ and $\VV$. 

\proclaim{Corollary 5.3.7}
Passing to $G$-invariants 
in the chain complex 
$\roman{Hom}^{\tau^{\SSS}}(\SSS[s^2\fra g], \VV)^{\fra g}$
calculating 
$\roman{Ext}_{(C \fra g, \fra g)}(\Bobb R,\VV)$ induces
a canonical isomorphism
$$
\roman{Ext}_{((G,C \fra g);\chain)}(\Bobb R,\VV)
\cong \roman{Ext}_{(C \fra g, \fra g)}(\Bobb R,\VV)^{\pi_{0}(G)} . \qed
$$
\endproclaim

Combining this corollary with
Theorem 3.6.4 we arrive at the following.

\proclaim{Corollary 5.3.8}
The twisted Hom-object
$\roman{Hom}^{\tau}(\SSS, \VV^G)$
is a small model for
\linebreak
$\roman{Ext}_{((G,C \fra g);\chain)}(\Bobb R,\VV)$
that is compatible with the bundle structures
in the sense that
the obvious pairing
$$
\roman{Hom}(\SSS,\Bbb R)
\otimes
\roman{Hom}^{\tau}(\SSS, \VV^G)
\longrightarrow
\roman{Hom}^{\tau}(\SSS, \VV^G)
$$
induces the $\roman H^*(BG)$-module structure on
$\roman{Ext}_{((G,C \fra g);\chain)}(\Bobb R,\VV)$. \qed
\endproclaim

Let $G_0$ be the connected component of the identity.
It is worthwhile noting that, under the circumstances of
Corollary 5.3.8, $\Lambda = \roman H_*(G_0)$
and $\SSS=\roman H_*(BG_0)$.

We will refer to the twisted object
$\roman{Hom}^{\tau}(\SSS, \VV^G)$
as the {\it small Cartan model\/}
for $\roman{Ext}_{((G,C \fra g);\chain)}(\Bobb R,\VV)$
associated with $G$ and $\VV$.

\medskip\noindent
{\bf 6. Small models in equivariant de Rham theory}
\smallskip\noindent
Let $G$ be a reductive Lie group and $X$ a left $G$-manifold.

Substitution of $\Cal A(X)$ for $\VV$ in (5.1.14) yields the model
$$
\roman{Hom}^{\tau^{\SSS},\tau_{|E\fra g|}\circ\iota}
(\SSS[s^2\fra g],|\Cal A^0(EG,\Cal A(X))|)^G. \tag6.1
$$
for the $G$-equivariant der Rham cohomology of $X$. An explicit
chain equivalence between $|\Cal A(N(G,X))|$ and (6.1) arises from
combination of the contraction (5.1.7) and the chain equivalence
(5.1.12) together with 
the Cartan twist, with $\Cal A(X)$ being substituted for $\VV$.

Suppose that $G$ is compact. Substitution
of $\Cal A(X)$ for $\VV$ in the Weil model
$\roman{Hom}(W'[\fra g], \VV)^{(G,C\fra g)}$
for  $\roman{Ext}_{((G,C \fra g);\chain)}(\Bobb R,\VV)$
associated with $G$ and $\VV$, cf. (5.3) above,
then yields the {\it Weil model\/} 
for the $G$-equivariant de Rham cohomology of $X$.
Likewise,
substitution
of $\Cal A(X)$ for $\VV$ in 
the Cartan model
$\roman{Hom}^{\tau^{\SSS}}(\SSS[s^2\fra g], \VV)^G$
for $\roman{Ext}_{((G,C \fra g);\chain)}(\Bobb R,\VV)$
associated with $G$ and $\VV$, cf. (5.3) above,
then yields the {\it Cartan
model\/}
$$
{\roman{Hom}^{\tau^{\SSS}}(\SSS[s^2\fra g], \Cal A(X))^G}
\tag6.2
$$
for the $G$-equivariant de Rham cohomology of $X$. An explicit
chain equivalence between $|\Cal A(N(G,X))|$ and the Cartan model arises
from combination of the above chain equivalence between $|\Cal A(N(G,X))|$
and (6.1) with (5.3.6) together with the Cartan twist, 
with $\Cal A(X)$ being substituted for
$\VV$. 

In the same vein,  substitution of $\Cal A(X)$ for $\VV$ in 
the small Cartan model
$\roman{Hom}^{\tau}(\SSS, \VV^G)$
for $\roman{Ext}_{((G,C \fra g);\chain)}(\Bobb R,\VV)$
associated with $G$ and $\VV$, cf. (5.3) above,
yields the {\it small model\/}
$$
\roman{Hom}^{\tau}(\roman H_*(BG),\Cal A(X)^G) \tag6.3
$$
{\it for the $G$-equivariant de Rham cohomology of\/} $X$.

\smallskip\noindent
{Remark 6.4\/.} {\sl In\/} \cite\gorkomac\ (Section 8 and thereafter) {\sl a 
small model for equivariant
de Rham theory of the kind\/} (6.3) {\sl is explored.}
Incomplete reasoning and faulty usage of this model
led to a certain activity in the literature
\cite\allpuone, \cite\almei, \cite\gorkomac, \cite\maszwebe,
see in particular the introduction of
\cite\almei;
cf. Remark 3.6.10 above.

\smallskip
\noindent {\smc 6.5. Homogeneous spaces\/.} For illustration,
suppose that $G$ is a closed subgroup of a compact connected Lie
group $K$. The $G$-equivariant cohomology of $K$ equals the
cohomology of the homogeneous space $K/G$. In the small model
(6.3), with $X=K$, we may replace the de $G$-invariant de Rham
algebra $\Cal A (K)^G$ with the cohomology $\roman H^*(K)$ which,
in fact, sits inside $\Cal A (K)$ as the graded subalgebra of
biinvariant forms. The resulting model 
for the cohomology of $G/K$ 
has the form
$\roman{Hom}^{\vartheta}(\roman H_*(BG),\roman H^*(K))$, the
twisting cochain $\vartheta$ being given as the composite of the
induced morphism from $\roman H_*(BG)$ to $\roman H_*(BK)$ with
the transgression twisting cochain from $\roman H_*(BK)$ to
$\roman H_*(K)$. This is the Cartan model for the de Rham
cohomology  of the homogeneous space $K/G$ \cite\cartantw.

\smallskip\noindent {\smc 6.6. Multiplicative cohomology generators\/.} The
group $G$  being supposed compact and connected, 
for the special case where $\VV=\Bobb R$, consider the
canonical injection
$$
\roman H^*(BG) @>>>
\overline {\Cal B}^*_{(G,C \fra g)}(G)
\tag6.6.1
$$
which is the  composite
of the injections 
in (5.1.12) and (5.3.6)
for $\VV=\Bobb R$;
the injection (6.6.1) is plainly 
is an isomorphism on cohomology. 
The
composite of (6.6.1) with the injection $\alpha^*$ in (5.1.7) for
the special case where $\VV=\Bbb R$ yields an injection
$$
\roman H^*(BG) @>>> |\Cal A(NG)|, \tag6.6.2
$$
manifestly  a cohomology isomorphism. This injection involves {\it no
choices\/} \ at all and is, in particular, natural in $G$; it is
certainly {\it not\/} multiplicative unless $\roman H^*(BG)$ is
trivial or a polynomial algebra in a single generator. This
injection includes the construction of representatives in $|\Cal
A(NG)|$ of the multiplicative cohomology generators of $\roman
H^*(BG)$, similar to that given in \cite\duponone\ and
\cite\shulmone\ via the simplicial Chern-Weil construction. It is
interesting to note that the present construction 
of the injection (6.6.2) does not involve
curvature arguments.

\medskip\noindent
{\bf 7. Duality}
\smallskip\noindent
Suppose that the group $G$ is of finite homological type
(that is, its homology is finite in each degree).
Recall that 
the functor
$$
t^* \colon 
\roman{Mod}_{(G,C \fra g)} @>>> 
{}_{\overline {\Cal B}{^*_{(G,C \fra g)}}(G)}\roman{Mod}
$$
has been defined above, cf. (5.2.1).
Define the functor
$$
h^* \colon  
{}_{\overline {\Cal B}{^*_{(G,C \fra g)}}(G)}\roman{Mod}
@>>> 
\roman{Mod}_{(G,C \fra g)}
\tag7.1 
$$
by the assignment to a  
$(\overline {\Cal B}^*_{(G,C \fra g)}(G))$-module $\Mflat$ 
of the twisted object 
$$
h^*(\Mflat) 
= \Cal B^*_{(G,C \fra g)}
(\Bobb R,G,\Mflat)^{\overline {\Cal B}^*_{(G,C \fra g)}(G)}.
\tag7.2
$$
The latter inherits a canonical $(G,C \fra g)$-module structure:
Indeed,
$$
\Cal B^*_{(G,C \fra g)}(\Bobb R,G,\Mflat)
=
\roman{Hom}^{\tau_{|E\fra g|}}
(\roman B\Lambda_{\partial}[s \fra g], \Cal A^0(EG, \Mflat))
$$
inherits a $(G,C \fra g)$-module structure from the 
$G$-actions on $\roman B\Lambda_{\partial}[s \fra g]$
and on $EG$;
further, given 
the $(\overline {\Cal B}^*_{(G,C \fra g)}(G))$-module $\Mflat$,
the action of
$$
\overline {\Cal B}^*_{(G,C \fra g)}(G)
=\Cal B^*_{(G,C \fra g)}(\Bobb R,G,\Bobb R)^{(G,C \fra g)}
$$ 
on 
$\Cal B^*_{(G,C \fra g)}(\Bobb R,G,\Mflat)$
preserves the $(G,C \fra g)$-module structure
whence the chain complex of 
$(\overline {\Cal B}^*_{(G,C \fra g)}(G))$-
invariants inherits a $(G,C \fra g)$-module structure.
Since $G$ is of finite homological type, 
the twisted object $h^*(\Mflat)$ 
{\it calculates\/} the differential graded 
$\roman{Tor}^{\overline {\Cal B}^*_{(G,C \fra g)}(G)}
(\Bobb R,\Mflat)$.
The functors $t^*$ and $h^*$ 
are formally different from those 
denoted by $t^*$ and $h^*$ in \cite\duaone\ (4.1); indeed,
since the category of $(G,C \fra g)$-modules is not one of modules 
over a chain algebra
the framework of \cite\duaone\ (4.1) is {\it not\/} directly applicable. 

\proclaim{Proposition 7.3}
The functors $t^*$ and $h^*$ are homotopy inverse to each other.
\endproclaim

\demo{Proof} Given a $(G,C \fra g)$-module $\Nflat$,
the canonical injection $\Nflat@>>> h^*(t^*(\Nflat))$ is a morphism
of $(G,C \fra g)$-modules; given a
$\overline {\Cal B}{^*_{(G,C \fra g)}}(G)$-module $\Mflat$,
the canonical injection 
$\Mflat @>>> t^*(h^*(\Mflat))$ 
is a morphism of
$\overline {\Cal B}{^*_{(G,C \fra g)}}(G)$-modules.
As a $(G,C \fra g)$-module, $h^*(t^*(\Nflat))$ amounts to 
$$
\Cal B^*_{(G,C \fra g)}(\Bobb R,G,\Nflat)
=
\roman{Hom}^{\tau_{|E\fra g|}}
(\roman B\Lambda_{\partial}[s \fra g], \Cal A^0(EG, \Nflat)),
$$
with the diagonal $(G,C \fra g)$-module structure.
This corresponds to the fact that, for a space $Y$ over $BG$,
the Borel construction, applied to the total space $P_Y$ of the
induced $G$-bundle, yields the space $EG \times Y$.
Likewise,
as a $\overline {\Cal B}{^*_{(G,C \fra g)}}(G)$-module,
$t^*(h^*(\Mflat))$ amounts to 
$$
\Cal B^*_{(G,C \fra g)}(\Bobb R,G,\Mflat)
=
\roman{Hom}^{\tau_{|E\fra g|}}
(\roman B\Lambda_{\partial}[s \fra g], \Cal A^0(EG, \Mflat)),
$$
with the diagonal 
$\overline {\Cal B}{^*_{(G,C \fra g)}}(G)$-module structure.
This corresponds to the fact that, for a $G$-space $X$,
the total space $P_Y$ of the induced $G$-bundle
over the Borel construction $Y = EG \times_GX$
amounts to $EG \times X$.
Furthermore, since 
$\roman B\Lambda_{\partial}[s \fra g]$ is contractible 
(cf. Section 3 above) and since,
for any chain complex $V$,
$\Cal A^0(EG, V)$ 
contracts onto $V$, cf. (5.3.3),
for any chain complex $V$, the injection of $V$ into 
$\Cal B^*_{(G,C \fra g)}(\Bobb R,G,V)$ is a 
chain equivalence. Hence the injections
$\Mflat @>>> t^*(h^*(\Mflat))$ and
$\Nflat@>>> h^*(t^*(\Nflat))$ are chain equivalences,
in fact, may be extended to contractions in a canonical way.\qed
\enddemo

We will now apply the duality spelled out in Proposition 7.3 to spaces.
Let $Y$ be a simplicial space $Y$ over the simplicial space $NG$ and
consider the fiber square
$$
\CD
P_Y @>>> EG
\\
@VVV
@VVV
\\
Y @>>> NG
\endCD
\tag7.4
$$
of simplicial spaces, the left-hand arrow being the induced simplicial
principal $G$-bundle over $Y$. Via the induced morphism of differential 
graded algebras from $|\Cal A(NG)|$ to $|\Cal A(Y)|$,
the chain complex
$|\Cal A(Y)|$ inherits a 
differential graded $|\Cal A(NG)|$-module structure and the de Rham cohomology 
of the fiber $P_Y$ is canonically isomorphic to the 
differential graded
$$
\roman {Tor}^{|\Cal A(NG)|}(\Bobb R,|\Cal A(Y)|)
\tag7.5
$$
which, by definition, is the homology of the bar construction
$$
\roman B(\Bobb R,|\Cal A(NG)|,|\Cal A(Y)|)
= \rbar |\Cal A(NG)| \otimes_{\tau^{\rbar}}  |\Cal A(Y)|.
\tag7.6
$$
The duality spelled out in Proposition 7.3 does not apply directly,
since the $|\Cal A(NG)|$-module structure on $|\Cal A(Y)|$ does not factor 
through a $({\overline {\Cal B}^*_{(G,C \fra g)}(G)})$-module structure.
Now the contraction (5.1.7), with $\VV =\Bbb R$, takes the form
$$
\Nsddata {|\Cal A(NG)|} 
{\alpha^*}{\nabla^*} {\overline {\Cal
B}^*_{(G,C \fra g)}(G)} {h^*} 
\tag7.7
$$
and,
in the category of sh-algebras, $|\Cal A(NG)|$ and 
${\overline {\Cal B}^*_{(G,C \fra g)}(G)}$ are isomorphic 
via the contraction (7.7); this notion of isomorphism
is explained in Section 6 of \cite\duaone.  
Hence $|\Cal A(Y)|$ inherits an sh-module structure over
${\overline {\Cal B}^*_{(G,C \fra g)}(G)}$,
unique up to homotopy. We now make this explicit.

To this end, we apply the construction \cite\duaone\ (2.2.1${}^*$) to 
the contraction (7.7), the bar construction twisting cochain
from
$\rbar\, {\overline {\Cal B}^*_{(G,C \fra g)}(G)}$
to
${\overline {\Cal B}
^*_{(G,C \fra g)}(G)}
$
being substituted for the twisting $\sigma$ in
\cite\duaone\ (2.2.1${}^*$). This
yields the acyclic twisting cochain
$$
\xi
\colon
\rbar\, {\overline {\Cal B}
^*_{(G,C \fra g)}(G)}
@>>>
|\Cal A(NG)|.
$$
Via the adjoint
$$
\overline \xi
\colon
\rcob \,\rbar\, {\overline {\Cal B}
^*_{(G,C \fra g)}(G)}
@>>>
|\Cal A(NG)|,
$$
we view henceforth any $|\Cal A(NG)|$-module as an
$\rcob\, \rbar\, {\overline {\Cal B}
^*_{(G,C \fra g)}(G)}$-module, that is, as an
sh-module over $\overline {\Cal B}^*_{(G,C \fra g)}(G)$.

On the category 
${}_{\overline {\Cal B}{^*_{(G,C \fra g)}}(G)}\roman{Mod}^{\infty}$
of sh-modules over $\overline {\Cal B}^*_{(G,C \fra g)}(G)$,
consider the functor
$$
H_{\infty}^* \colon  
{}_{\overline {\Cal B}{^*_{(G,C \fra g)}}(G)}\roman{Mod}^{\infty}
@>>> 
{}_{\rbar\, {\overline {\Cal B}
^*_{(G,C \fra g)}(G)}}
\roman{Comod}
$$
which assigns to an arbitrary sh-module $(M,\tau_{\rcob\, \rbar})$
over $\overline {\Cal B}^*_{(G,C \fra g)}(G)$
the twisted object
$$
H_{\infty}^* (M,\tau_{\rcob\, \rbar})=
\rbar\, {\overline {\Cal B}
^*_{(G,C \fra g)}(G)}
\otimes_{\tau_{\rcob\, \rbar}}  M;
\tag7.8
$$
here $\tau_{\rcob\, \rbar}$ refers to the
universal twisting cochain from  $\rbar$ to $\rcob\, \rbar$.
In particular, the twisting cochain $\xi$
induces an sh-structure on $|\Cal A(Y)|$,
and
$$
H_{\infty}^* (|\Cal A(Y)|,\tau_{\rcob\, \rbar})
=
\rbar\, {\overline {\Cal B}
^*_{(G,C \fra g)}(G)}
\otimes_{\xi}  |\Cal A(Y)|;
$$
we will simplify the notation and write
$$
H_{\infty}^* (|\Cal A(Y)|,\xi)
=
H_{\infty}^* (|\Cal A(Y)|,\tau_{\rcob\, \rbar}).
$$
Since
the induced bundle morphism
$$
\overline \xi \otimes \roman{Id}\colon \rbar\, {\overline {\Cal B}
^*_{(G,C \fra g)}(G)}
\otimes_{\xi}  |\Cal A(Y)|
@>>>
\rbar |\Cal A(NG)| \otimes_{\tau^{\rbar |\Cal A(NG)|}}  |\Cal A(Y)|
\tag7.9
$$
is a chain equivalence, in fact, can be extended to a contraction via the 
perturbation lemma, cf. \cite\duaone\ (2.3),
the twisted object $H_{\infty}^* (|\Cal A(Y)|,\xi)$ 
calculates the 
differential torsion product (7.5).
Thus the coalgebra
$\rbar\, {\overline {\Cal B}^*_{(G,C \fra g)}(G)}$ is a replacement 
for the de Rham complex $\Cal A(G)$ and the missing coalgebra structure 
thereupon. By functoriality, for an arbitrary 
sh-module
$(M,\tau_{\rcob\, \rbar})$
over $\overline {\Cal B}^*_{(G,C \fra g)}(G)$,
the obvious right $(G,C\fra g)$-module structure
on $\overline {\Cal B}^*_{(G,C \fra g)}(G)$ induces a 
$(G,C\fra g)$-module structure on the twisted object 
$\rbar\, {\overline {\Cal B}^*_{(G,C \fra g)}(G)}
\otimes_{\tau_{\rcob\, \rbar}}  M$. In this fashion, 
we view $H_{\infty}^*$ as a functor of the kind
$$
H_{\infty}^* \colon  
{}_{\overline {\Cal B}{^*_{(G,C \fra g)}}(G)}\roman{Mod}^{\infty}
@>>> 
\roman{Mod}_{(G,C \fra g)} .
\tag7.10
$$
The duality between the two functors $t^*$ and $h^*$ spelled out in 
Proposition 7.3 above entails the following:

\proclaim{Theorem 7.11}
On the category of left $G$-manifolds,
the functor
$h^*\circ t^* \circ \Cal A$ is chain-equivalent to the functor
$\Cal A$ as $(G,C \fra g)$-module functors;
and on the category of simplicial manifolds
over $NG$, the functor
$t^*\circ H_{\infty}^* \circ |\Cal A|$ is chain-equivalent to the functor
$|\Cal A|$ as sh-module functors over 
$\overline {\Cal B}^*_{(G,C \fra g)}(G)$.
In particular,  application of the functor $h^*$ to the twisted
object $t^* (\Cal A (X))$ 
($=\overline {\Cal B}^*_{(G,C \fra g)}(\Bbb R,G,\Cal A (X))$)
calculating the $G$-equivariant de Rham cohomology of $X$ 
(in view of Theorem {\rm 2.7.1}) 
reproduces an object calculating the ordinary de Rham
cohomology of $X$;
and application of the functor $t^*$ to the twisted
object $H_{\infty}^* (|\Cal A(Y)|,\xi)$ 
reproduces an object calculating the de Rham cohomology of 
the simplicial space $Y$. \qed
\endproclaim

\noindent
{\smc 7.12. Koszul duality.\/}
Let $\Lambda = \roman H_*G$,
$\SSS = \roman H_*(BG)$, $\Lambda' = \roman H^*G$,
$\Sigm = \roman H^*(BG)$, and let $\tau \colon \SSS \to
\Lambda$ be the transgression twisting cochain.
Ordinary Koszul duality involves the two functors
$$
\align
t^* &\colon 
\roman{Mod}_{\Lambda} @>>> {}_{\Sigm}{\roman{Mod}},
\quad
t^*(N) = \roman{Hom}^{\tau}(\SSS,N)
\\
h^* &\colon  
{}_{\Sigm}{\roman{Mod}} @>>> \roman{Mod}_\Lambda,
\quad
h^*(M) = \roman{Hom}^{\tau}(\Lambda,M).
\endalign
$$
The former assigns to a (right) $\Lambda$-module $N$ the twisted Hom-object 
$t^*(N)$ which calculates the differential graded 
$\roman{Ext}_{\Lambda}(\Bobb R,N)$, and the latter
assigns to a 
(left) $\Sigm$-module $M$ the twisted Hom-object 
$h^*(M)$ which, since $\Sigm$  is of finite type,
calculates the differential graded $\roman{Tor}^{\Sigm}(\Bobb R,M)$. 
These functors are chain homotopy inverse to each other in an 
obvious manner.

Replace $|\Cal A(X)|$ with the ordinary $\Lambda$-module 
$$
\Lambda' \otimes_{\zeta_{\Cal B^*}}t^*(|\Cal A(X)|)
=
\Lambda' \otimes_{\zeta_{\Cal B^*}}
\overline {\Cal B}{^*_{(G,C \fra g)}}
(\Bobb R,G,|\Cal A(X)|),
\tag7.12.1
$$
cf. (5.2) above.
In the same vein, given a simplicial space $Y$ over $NG$, 
we replace $|\Cal A(Y)|$ with an ordinary $\Sigm$-module as follows
where we write
$$
\overline {\Cal B}^*=
\overline {\Cal B}^*_{(G,C \fra g)}
(\Bobb R,G,\Bobb R)
$$
for simplicity:
Extend the adjoint $\overline \zeta_{\Cal B^*}$
of the twisting cochain (5.2.5) to a contraction
$$
\Nsddata 
{\rbar \,\overline {\Cal B^*}}
{\overline \zeta_{\Cal B^*}}{\alpha}
{\Lambda'}{h}.
\tag7.12.2
$$
The construction \cite\duaone\ (2.2.1${}_*$), applied to (7.12.2) and the
(acyclic) transgression twisting cochain $\tau^* \colon \Lambda' \to \Sigm$
(the dual of the
transgression twisting cochain $\tau \colon \SSS \to \Lambda$)
yields the acyclic twisting cochain
$$
\zeta^{\rbar \,\overline {\Cal B^*}}\colon \rbar \overline {\Cal B^*} 
@>>> \Sigm .
\tag7.12.3
$$
This twisting cochain determines the twisted object
$$
\Sigm \otimes_{\zeta^{\rbar \,\overline{\Cal B}^* }} 
H_{\infty}^* (|\Cal A(Y)|,\xi)
\tag7.12.4
$$
which, for our purposes, is the appropriate
replacement for $|\Cal A(Y)|$.
This twisted object is, in particular, an ordinary $\Sigm$-module.
{\sl With these twisted objects, a version of
Koszul duality is given by the functors
$t^*$ and $h^*$ between the categories
${}_{\Sigm}\roman{Mod}$ and $\roman{Mod}_{\Lambda}$:
The functor $h^*$ reconstructs the
ordinary cohomology of a $G$-manifold $X$ from a model
of the kind {\rm (7.12.4)}
for the $G$-equivariant cohomology\/} (where the construction
(7.12.4) is carried out for $Y=N(G,X)$);
and {\sl the functor $t^*$ reconstructs the
equivariant cohomology of a $G$-manifold from a model of the kind 
{\rm (7.12.1)} for the ordinary cohomology.\/}
This corresponds to the procedure employed in \cite\gorkomac\
(cf. e.~g. p. 29) which consists in replacing the naive cochain complexes,
where the $\Lambda$- and $\Sigm$-actions are {\it not\/} defined,
by equivalent cochain complexes where the actions {\it are\/}
defined.

\smallskip\noindent
{\smc 7.13. Koszul duality when $G$ is finite dimensional, 
compact and 
connected,\/} cf. e.~g. \cite\gorkomac.
We recall it, to establish the link with the theory built up above.
As noted in (5.3) above, cf. Corollary 5.3.8,
given $X$, the algebra 
$\Cal A(X)^G$ of invariants inherits now a $(\roman H_*(G))$-module
structure, and the model (6.3) for the $G$-equivariant
de Rham cohomology is exactly
$t^*(\Cal A(X)^G)$.
The functor $h^*$ reconstructs the 
ordinary cohomology of $X$ from $t^*(\Cal A(X)^G)$.
Under our general circumstances
(where $G$ is a general, possibly infinite dimensional Lie group),
{\rm (7.12.1)} is a replacement for (6.3) and
the functor $t^*$ applies, for a general simplicial space over $NG$,
to the model {\rm (7.12.4)}.

\bigskip

\widestnumber\key{999} 
\centerline{References}

\ref \no \allpuone
\by C. Allday and V. Puppe
\paper On a conjecture of Goresky, Kottwitz and MacPherson
\jour Canad. J. of Mathematics
\vol 51
\yr 1999
\pages 3-9
\endref

\ref \no \almei
\by A. Alexeiev and E. Meinrenken
\paper Equivariant cohomology and the Maurer-Cartan equation
\jour Duke Math. J. \vol 130
\yr 2005 \pages  479--521 \finalinfo{\tt math.DG/0406350}
\endref

\ref \no \barrone
\by M.~Barr
\paper Cartan-Eilenberg cohomology and triples
\jour Journal of Pure and Applied Algebra
\vol 112
\yr 1996
\pages  219--238
\endref

\ref \no \bottone \by R. Bott \paper On the Chern-Weil
homomorphism and the continuous cohomology of Lie groups \jour
Advances \vol 11 \yr 1973 \pages  289--303
\endref

\ref \no \botshust \by R. Bott, H. Shulman, and J. Stasheff \paper
On the de Rham theory of certain classifying spaces \jour Advances
\vol 20 \yr 1976 \pages 43--56
\endref

\ref \no \cartanon
\by H. Cartan
\paper Notions d'alg\`ebre diff\'erentielle; 
applications aux groupes 
de Lie et aux vari\'et\'es o\`u op\`ere un groupe de Lie
\jour Coll. Topologie Alg\'ebrique
\paperinfo Bruxelles
\yr 1950
\pages  15--28
\endref

\ref \no \cartantw \bysame \paper La transgression dans un
groupe de Lie et dans un espace fibr\'e principal \jour Coll.
Topologie Alg\'ebrique \paperinfo Bruxelles \yr 1950 \pages
57--72
\endref

\ref \no \cartanse
\bysame
\paper Alg\`ebres d'Eilenberg--Mac Lane et homotopie
\paperinfo expos\'es 2--11
\jour S\'eminaire H. Cartan 1954/55
\publ Ecole Normale Superieure, Paris, 1956
\endref

\ref \no \cartanei \by H. Cartan and S. Eilenberg \book
Homological Algebra \publ Princeton University Press \publaddr
Princeton \yr 1956
\endref

\ref \no \cheveile
\by C. Chevalley and S. Eilenberg
\paper Cohomology theory of Lie groups and Lie algebras
\jour  Trans. Amer. Math. Soc.
\vol 63
\yr 1948
\pages 85--124
\endref

\ref \no \doldpupp \by A. Dold und D. Puppe \paper Homologie
nicht-additiver Funktoren. Anwendungen \jour Annales de l'Institut
Fourier \vol 11 \yr 1961 \pages  201--313
\endref

\ref \no \duponone \by J. L. Dupont \paper Simplicial de Rham
cohomology and characteristic classes of flat bundles \jour
Topology \vol 15 \yr 1976 \pages  233--245
\endref

\ref \no \duskinon
\by J.  Duskin
\paper Simplicial methods and the interpretation of \lq\lq triple\rq\rq\ 
cohomology
\jour Memoirs Amer. Math. Soc.
\vol 163
\yr 1975
\endref

\ref \no \eilmothr
\by S. Eilenberg and J. C. Moore
\paper Foundations of relative homological algebra
\jour Memoirs AMS
\vol 55
\yr 1965
\publ Amer. Math. Soc.
\publaddr Providence, Rhode Island
\endref

\ref \no \franzone \by M. Franz \paper Koszul duality and
equivariant cohomology for tori \jour Int. Math. Res. Not. \vol 42
\yr 2003 \pages 2255--2303 \finalinfo{\tt math.AT/0301083}
\endref

\ref \no \franztwo \bysame \paper Koszul duality and
equivariant cohomology \finalinfo{\tt math.AT/0307115}
\endref

\ref \no \godebook
\by R. Godement
\book Topologie alg\'ebrique et th\'eorie des faisceaux
\publ Hermann
\publaddr Paris
\yr 1958
\endref

\ref \no \gorkomac \by M. Goresky, R. Kottwitz, and R. Mac Pherson
\paper Equivariant cohomology, Koszul duality and the localization
theorem \jour Invent. Math. \vol 131 \yr 1998 \pages 25--83
\endref

\ref \no \gugenhtw \by V.K.A.M. Gugenheim \paper On the chain
complex of a fibration \jour Illinois J. of Mathematics \vol 16
\yr 1972 \pages 398--414
\endref

\ref \no \gugenmay \by V.K.A.M. Gugenheim and J.P. May \paper On
the theory and applications of differential torsion products \jour
Memoirs of the Amer. Math. Soc. \vol 142 \yr 1974
\endref

\ref \no \gugenmun
\by V.K.A.M. Gugenheim and H. J. Munkholm
\paper On the extended functoriality of Tor and Cotor
\jour J. of Pure and Applied Algebra
\vol 4
\yr 1974
\pages  9--29
\endref

\ref \no \hochsone
\by G. Hochschild
\paper Relative homological algebra
\jour  Trans. Amer. Math. Soc.
\vol 82
\yr 1956
\pages 246--269
\endref

\ref \no \hochmost \by G. Hochschild and G. D. Mostow \paper
Cohomology of Lie groups \jour  Illinois J.  of Math. \vol 6 \yr
1962 \pages  367--401
\endref

\ref \no \habili \by J. Huebschmann \paper Perturbation theory and
small models for the chains of certain induced fibre spaces
\paperinfo Habilitationsschrift, Universit\"at Heidelberg, 1984
\finalinfo {\bf Zbl} 576.55012
\endref

\ref \no \perturba \bysame \paper Perturbation theory
and free resolutions for nilpotent groups of class 2 \jour J. of
Algebra \yr 1989 \vol 126 \pages 348--399
\endref

\ref \no \cohomolo \bysame \paper Cohomology of
nilpotent groups of class 2 \jour J. of Algebra \yr 1989 \vol 126
\pages 400--450
\endref

\ref \no \modpcoho \bysame \paper The mod $p$
cohomology rings of metacyclic groups \jour J. of Pure and Applied
Algebra \vol 60 \yr 1989 \pages 53--105
\endref

\ref \no \intecoho \bysame \paper Cohomology of
metacyclic groups \jour Trans. Amer. Math. Soc. \vol 328 \yr 1991
\pages 1-72
\endref

\ref \no \kan \bysame \paper Extended moduli spaces,
the Kan construction, and lattice gauge theory \jour Topology \vol
38 \yr 1999 \pages 555--596 \finalinfo{\tt dg-ga/9505005,
dg-ga/9506006}
\endref

\ref \no \poiscoho
\bysame
\paper Poisson cohomology and quantization
\jour 
J. reine angew. Math.
\vol  408 
\yr 1990
\pages 57--113
\endref

\ref \no \extensta
\bysame
\paper 
Extensions of Lie-Rinehart algebras and the Chern-Weil construction
\paperinfo in: Festschrift in honor of J. Stasheff's 60-th birthday
\jour Cont. Math. 
\vol 227
\yr 1999
\pages 145--176
\publ Amer. Math. Soc.
\publaddr Providence R. I.
\finalinfo{\tt math.DG/9706002}
\endref

\ref \no \lradq
\bysame
\paper Lie-Rinehart algebras, descent, and quantization
\paperinfo in: Galois theory, Hopf algebras, and semiabelian categories
\jour Fields Institute Communications 
\vol 43
\yr 2004
\pages 295--316
\publ Amer. Math. Soc.
\publaddr Providence R. I.
\finalinfo {\tt math.SG/0303016}
\endref

\ref \no \minimult
\bysame 
\paper Minimal free multi models for chain algebras
\paperinfo in: Chogoshvili Memorial
\jour Georgian Math. J.
\vol 11
\yr 2004
\pages 733--752
\finalinfo {\tt math.AT/0405172}
\endref

\ref \no \duaone \bysame \paper Homological
perturbations, equivariant cohomology, and Koszul duality
\finalinfo{\tt math.AT/0401160}
\endref

\ref \no \pertlie \bysame 
\paper 
The Lie algebra perturbation lemma
\paperinfo
in: Festschrift in honor of M. Gerstenhaber's 80-th and
Jim Stasheff's 70-th birthday, Progress in Math., Birkh\"auser-Verlag
(to appear)
\finalinfo{\tt arxiv:0708.3977}
\endref

\ref \no \pertltwo \bysame \paper 
The sh-Lie algebra perturbation lemma
\finalinfo{\tt arxiv:0710.2070}
\endref

\ref \no \huebkade \by J. Huebschmann and T. Kadeishvili \paper
Small models for chain algebras \jour Math. Z. \vol 207 \yr 1991
\pages 245--280
\endref

\ref \no \huebstas
\by J. Huebschmann and J. D. Stasheff
\paper Formal solution of the master equation via HPT and
deformation theory
\paperinfo {\tt math.AG/9906036}
\jour Forum mathematicum
\vol 14
\yr 2002
\pages 847--868
\endref

\ref \no\husmosta \by D. Husemoller, J.~C. Moore, and J.~D.
Stasheff \paper Differential homological algebra and homogeneous
spaces \jour J. of Pure and Applied Algebra \vol 5 \yr 1974 \pages
113--185
\endref

\ref \no \kamtonfo
\by F. W. Kamber and Ph. Tondeur
\paper Alg\`ebres de Weil semi simpliciales
\jour C. R. Acad. Sci. Paris S\'er. A-B
\vol 276
\yr 1973
\pages A1407--1410 
\endref

\ref \no \kamtontw
\bysame 
\paper Characteristic invariants of foliated bundles
\jour Manuscripta Math.
\vol 11
\yr 1974
\pages 51--89  
\endref

\ref \no \kamtonfi
\bysame 
\paper Semi-simplicial Weil algebras and characteristic classes
\jour T\hataccent ohoku Math. J. (2)
\vol 30
\yr 1978
\pages 373-422
\endref

 \ref \no \kostaeig
\by B. Kostant
\paper Clifford algebra analogue of the {H}opf-{K}oszul-{S}amelson
theorem, the $\rho$-decomposition 
$\roman{C}(\frak g)=\roman{End}\,{V}\sb \rho \otimes
\roman{C}(\roman{P})$, 
and the $\frak g$-module structure of $\bigwedge \frak g$ 
\jour Adv. in Math.
\vol 125 no.~2 
\yr 1997 
\pages 275--350
\endref

\ref \no \maclafiv
\by S. Mac Lane
\paper Homologie des anneaux et des modules
\paperinfo in: Colloque de topologie alg\'ebrique, Louvain
\yr 1956
\pages 55--80
\endref

\ref \no \maclaboo \by S. Mac Lane \book Homology \bookinfo Die
Grundlehren der mathematischen Wissenschaften
 No. 114
\publ Springer \publaddr Berlin $\cdot$ G\"ottingen $\cdot$
Heidelberg \yr 1963
\endref

\ref \no \maclbotw
\bysame 
\book Categories for the Working Mathematician
\bookinfo Graduate Texts in Mathematics 
\vol 5
\publ Springer
\publaddr Berlin $\cdot$ G\"ottingen $\cdot$ Heidelberg
\yr 1971
\endref

\ref \no \maszwebe \by T. Maszczyk and A. Weber \paper Koszul
duality for modules over Lie algebras \jour Duke Math. J. \vol 112
\yr 2002 \pages  111--120 \finalinfo{\tt math.AG/0101180}
\endref

\ref \no \mooretwo \by J. C. Moore \paper Cartan's constructions
\paperinfo Colloque analyse et topologie, en l'honneur de Henri
Cartan \jour Ast\'erisque \vol 32--33 \yr 1976 \pages  173--221
\endref

\ref \no \moorefiv
\bysame 
\paper Differential homological algebra
\paperinfo Actes, Congres intern. math. Nice, 1970
\publ Gauthiers-Villars
\publaddr Paris
\yr 1971
\pages 335--339
\endref

\ref \no \quilltwo
\by D. Quillen
\paper Rational homotopy theory
\jour Ann. of Math.
\vol 90
\yr 1969
\pages  205--295
\endref

\ref \no \rinehone
\by G. Rinehart
\paper Differential forms for general commutative algebras
\jour  Trans. Amer. Math. Soc.
\vol 108
\yr 1963
\pages 195--222
\endref

\ref \no \gsegatwo
\by G. B. Segal
\paper Classifying spaces and spectral sequences
\jour Publ. Math. I. H. E. S. 
\vol 34
\yr 1968
\pages 105--112
\endref

\ref \no \shulmone \by H. B. Shulman \book Characteristic classes
and foliations \bookinfo Ph. D. Thesis \publ University of
California \yr 1972
\endref

\ref \no \stashsev
\by J. D. Stasheff
\paper Continuous cohomology of groups and classifying spaces
\jour Bull. Amer. Math. Soc.
\vol 84
\yr 1978
\pages 513--530
\endref

\ref \no \stashalp \by J.~D. Stasheff and S. Halperin \paper
Differential algebra in its own rite \jour Proc. Adv. Study Alg.
Top. August 10--23, 1970, Aarhus, Denmark \pages 567--577
\endref

\ref \no \vanestwo
\by W.~T. Van Est
\paper Une application d'une m\'ethode de Cartan-Leray
\jour Nederl. Akad. Wetensch. Proc. Ser. A 58 =
Indag. Math.
\vol 17
\yr 1955
\pages  542--544
\endref

\ref \no \vanesthr
\bysame  
\paper Alg\`ebres de Maurer-Cartan et
holonomie
\jour
Ann. Fac. Sci. Toulouse Math. 
\vol 5
\year 1989 (suppl.)
\pages 93--134
\endref
\enddocument